\documentclass[preprint,12pt]{elsarticle}
\usepackage[left=0.75in, right=1in, top=1in, bottom=1in]{geometry}
\usepackage{graphicx}
\usepackage{amssymb}
\usepackage{amsthm}
\usepackage[tbtags]{amsmath}
\usepackage{mathtools}
\numberwithin{equation}{section}
\usepackage{lineno}
\usepackage[font=small,labelfont=bf]{caption}
\usepackage{subcaption}
\usepackage{romannum}
\usepackage{array}
\usepackage{booktabs}
\usepackage{float}
\usepackage{cleveref}
\usepackage{bm}
\usepackage{nicefrac}
\usepackage{xcolor}

\def\hf{\frac{1}{2}}

\newcommand{\dt}{\Delta t}
\newcommand{\dx}{\Delta x}
\newcommand{\dy}{\Delta y}
\newcommand{\kph}{{k+\frac{1}{2}}}
\newcommand{\kmh}{{k-\frac{1}{2}}}
\newcommand{\jph}{{j+\frac{1}{2}}}
\newcommand{\jmh}{{j-\frac{1}{2}}}

\newenvironment{acknowledgments}{{\flushleft \bf Acknowledgment:$~$}}{}

\newcommand*\xbar[1]{%
  \hbox{%
    \vbox{%
      \hrule height 0.5pt 
      \kern0.4ex
      \hbox{%
        \kern-0.05em
        \ensuremath{#1}%
        \kern-0.05em
      }%
    }%
  }%
}

\makeatletter
\def\ps@pprintTitle{%
   \let\@oddhead\@empty
   \let\@evenhead\@empty
   \let\@oddfoot\@empty
   \let\@evenfoot\@oddfoot
}
\makeatother

\begin{document}
\begin{frontmatter}
\pagenumbering{arabic}

\title{An adaptive well-balanced positivity preserving central-upwind scheme on quadtree grids for shallow water equations}

\author[a]{Mohammad A. Ghazizadeh\corref{cor1}}
\ead{sghaz023@uottawa.ca}

\author[a]{Abdolmajid Mohammadian}
\ead{majid.mohammadian@uottawa.ca}

\author[b]{Alexander Kurganov}
\ead{alexander@sustech.edu.cn}
 
\cortext[cor1]{Corresponding author.}

\address[a]{Department of Civil Engineering, University of Ottawa, Ottawa, ON K1N 6N5, Canada}
\address[b]{Department of Mathematics and SUSTech International Center for Mathematics, Southern University of Science and Technology,
Shenzhen, 518055, China}

\begin{abstract}
We present an adaptive well-balanced positivity preserving central-upwind scheme on quadtree grids for shallow water equations. The use of
quadtree grids results in a robust, efficient and highly accurate numerical method. The quadtree model is developed based on the
well-balanced positivity preserving central-upwind scheme proposed in [{\sc A. Kurganov and G. Petrova}, {\em Commun. Math. Sci.}, 5 (2007),
pp. 133--160]. The designed scheme is well-balanced in the sense that it is capable of exactly preserving ``lake-at-rest'' steady states. In
order to achieve this as well as to preserve positivity of water depth, a continuous piecewise bilinear interpolation of the bottom
topography function is utilized. This makes the proposed scheme capable of modelling flows over discontinuous bottom topography. Local
gradients are examined to determine new seeding points in grid refinement for the next timestep. Numerical examples demonstrate the
promising performance of the central-upwind quadtree scheme.

\end{abstract}
\begin{keyword}
Shallow water equations, quadtree grids, central-upwind scheme, well-balanced scheme, positivity preserving scheme.
\end{keyword}
\end{frontmatter}
\section{Introduction}
\begin{figure}[ht!]
\centerline{\includegraphics[height=0.30\linewidth]{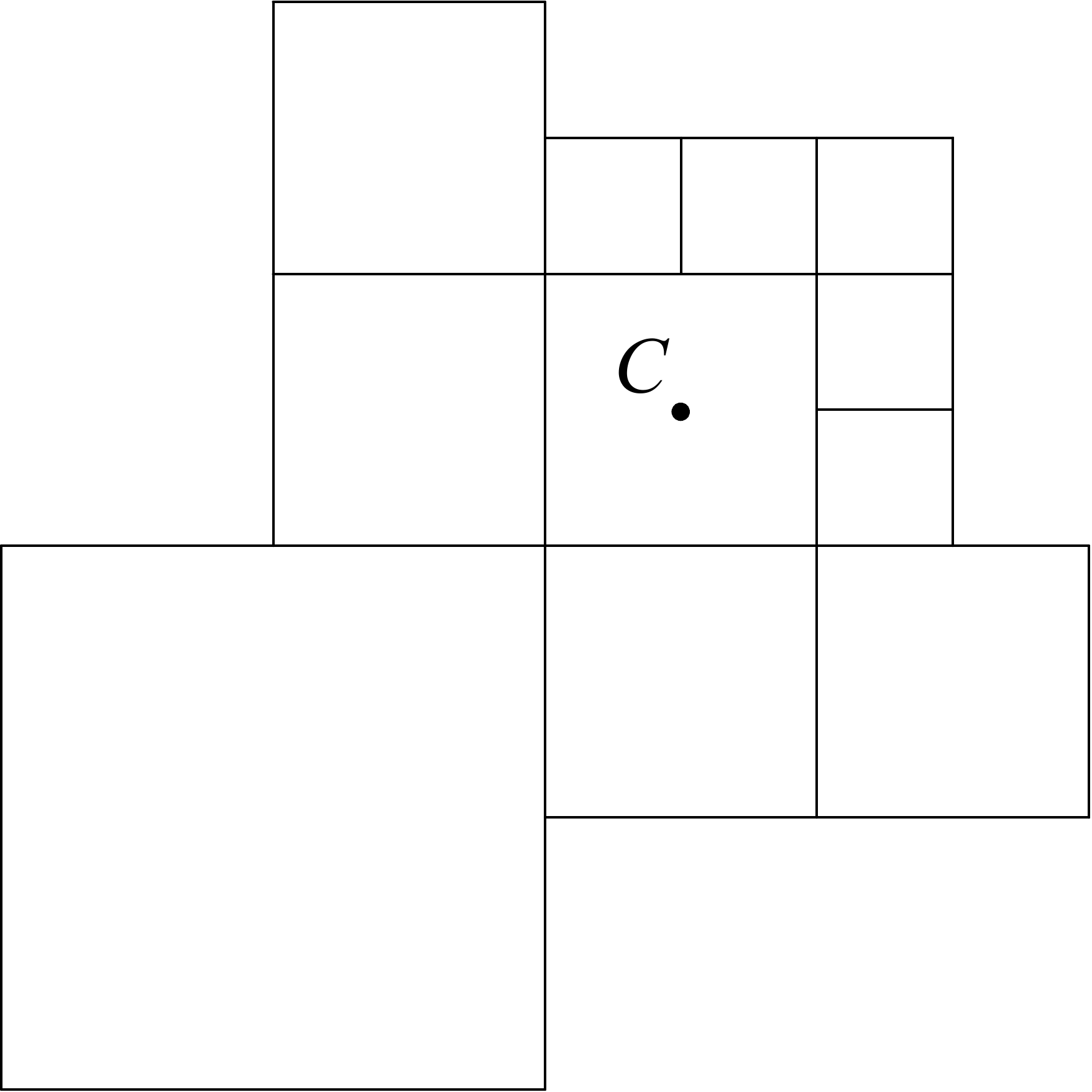}}
\caption{\sf Example of quadtree cells with different level neighboring cell sizes.\label{fig:1}}
\end{figure}

Quadtree grids (Figure \ref{fig:1}), which are two-dimensional (2-D) semi-structured Cartesian grids, are based on hierarchical data
structures, which are widely used in the field of computer science (e.g., image processing and computer graphics), computational geometry,
robotics, video games, and computational fluid dynamics (CFD); see, e.g., \cite{Samet1984,Yiu1996}. Several studies have been conducted on
how to generate quadtree grids; see, e.g., \cite{Aizawa2008,Borthwick2000,Greaves1998,Pascal1998,Popinet2010,Samet1984a,Samet2006}.

Cartesian grids are common in CFD problems because of their efficiency and ability to maintain the simplicity of discretized equations,
which reduces computational cost in comparison to unstructured grids. One of the benefits of quadtree grids over structured grids is grid
coarsening: while the accuracy is maintained, the grid can be coarsened wherever no refinement is needed and thus, the computational cost is
reduced. Note that a disadvantage of Cartesian grids is their inability to adequately represent complex shapes. In such situations, cut-cell
grids become useful; see, e.g., \cite{An2012}. This paper will only focus on quadtree grids.

The main goal of this paper is to develop an adaptive well-balanced positivity preserving scheme on quadtree grids for the Saint-Venant
system of shallow water equations (SWEs). This system was first proposed in \cite{Sai}, but is still extensively used to model flows in
rivers, lakes, coastal areas and estuaries \cite{Sai}. In the 2-D case, the SWEs can be written in terms of the water surface ($w$) and the
unit discharges ($hu$ and $hv$) as follows \cite{Bryson2005}:
\begin{equation}
\begin{dcases}
w_t+(hu)_x+(hv)_y=0,\\
(hu)_t+\bigg[\dfrac{(hu)^2}{w-B}+\dfrac{g}{2}(w-B)^2\bigg]_x+\bigg[\dfrac{(hu)(hv)}{w-B}\bigg]_y=-g\,(w-B)B_x,\\
(hv)_t+\bigg[\dfrac{(hu)(hv)}{w-B}\bigg]_x+\bigg[\dfrac{(hv)^2}{w-B}+\dfrac{g}{2}(w-B)^2\bigg]_y=-g\,(w-B)B_y,
\end{dcases}
\label{eq:SWE}
\end{equation}
where $t$ is time, $g$ is the gravitational constant, $x$ and $y$ are the directions in the 2-D Cartesian coordinate system, $u(x,y,t)$ and
$v(x,y,t)$ are the water velocities in the $x$- and $y$-directions, respectively, $B(x,y)$ is the bottom topography, and
$h(x,y,t)=w(x,y,t)-B(x,y)$ is the water depth.

The system \eqref{eq:SWE} admits ``lake-at-rest'' steady-state solutions,
\begin{equation}
w\equiv{\rm Const},\quad u=v\equiv0,
\label{1.2}
\end{equation}
which are of great practical importance as many waves to be captured are, in fact, small perturbations of these steady states. We would
like to stress that good numerical methods should be able of exactly preserving ``lake-at-rest'' steady states---such methods are called
{\em well-balanced}. Another important property a good numerical method should possess is its ability to preserve non-negativity of water
depth $h$---such methods are called {\em positivity preserving}. We refer the reader to, e.g., a recent review paper \cite{Kur_Acta} for an
extensive discussion on these matters.

Several numerical methods on quadtree grids for SWEs have been developed during the past two decades. For example, an adaptive well-balanced
second-order Godunov-type scheme was proposed in \cite{Rogers2001}. This scheme is able to solve the shallow water system with
discontinuous bottom topography. A
well-balanced scheme on quadtree-cut-cell grids was proposed in \cite{An2012}. This scheme is based on the hydrostatic reconstruction from
\cite{Audusse2005}. In addition, an adaptive quadtree Roe-type scheme for the 2-D two-layer SWEs was introduced in \cite{Lee2011}.
Furthermore, an adaptive quadtree scheme with wet-dry fronts was studied \cite{Liang2009,Liang2004}. For further studies on SWEs over
quadtree grids, we refer the reader to \cite{Borthwick2001,Borthwick2000,Kesserwani2012,Liang2015,Michel-dansac2016}. Besides the aforementioned numerical methods, several well-balanced positivity preserving schemes have been proposed in the past years;
see, e.g., \cite{Audusse2004,Audusse2005,BMK,BM08,BCKN,BNL,Bryson2011,GPC,Kurganov2002,Kurganov2007,LAEK,Ric15,SMSK}, but none of them has
been extended to quadtree grids.

In this paper, we present a central-upwind quadtree scheme which is based on the central-upwind scheme from \cite{Kurganov2007}. Central-upwind schemes are
Godunov-type Riemann-problem-solver-free finite-volume methods, which were proposed in \cite{KLin,KNP,KPW,KTcl} as a ``black-box'' solvers
for general multidimensional systems of hyperbolic systems of conservation laws. Central-upwind schemes were extended to shallow water
models in \cite{Kurganov2002} and many subsequent works; see, e.g., the recent review paper \cite{Kur_Acta} and references therein. The
scheme from \cite{Kurganov2007} is the first well-balanced and at the same time positivity preserving central-upwind scheme, which is
simple, efficient and robust: this is the reason why it was taken as the main building block of the proposed quadtree scheme.

The paper is organized as follows. In \S\ref{S:2}, we briefly describe a quadtree grid generation algorithm. In \S\ref{S:4}, we develop a
well-balanced positivity preserving central-upwind quadtree scheme. The developed scheme is tested on three numerical example in            
\S\ref{S:5}. Finally, some concluding remarks can be found in \S\ref{S:6}.

\section{Quadtree grids}\label{S:2}
Quadtree grids imply recursive spatial decomposition of the computational domain; see Figure \ref{fig:1} for an example of a quadtree cell
$C$ with different level neighboring cells.

Quadtree grids can be generated according to the following algorithm (see \cite{Borthwick2000,Greaves1998}):

\vskip5pt
\noindent\textbf{Step 1.} Choose a domain and generate a set of seeding points considering features of the problem, boundary conditions,
flow characteristics, local gradients and governing equations.

\vskip4pt
\noindent\textbf{Step 2.} Fit the domain within a unit square (root square) by adjusting the size of the square.

\vskip4pt
\noindent\textbf{Step 3.} Determine the level of refinement ($m$) of the quadtree (the size of the smallest cell of the grid is inversely
proportional to $m$).

\vskip4pt
\noindent\textbf{Step 4.} Divide the domain square into four sub-squares. Each sub-square is called a cell (that is, the first level of the
quadtree).

\vskip4pt
\noindent\textbf{Step 5.} Continue dividing each cell into four sub-cells if it contains a seeding point until the maximum refinement level
is reached. If a cell does not include any seeding points, move to the next cell and again implement Step 5. 

\begin{figure}[h!]
\centerline{\includegraphics[height=0.35\linewidth]{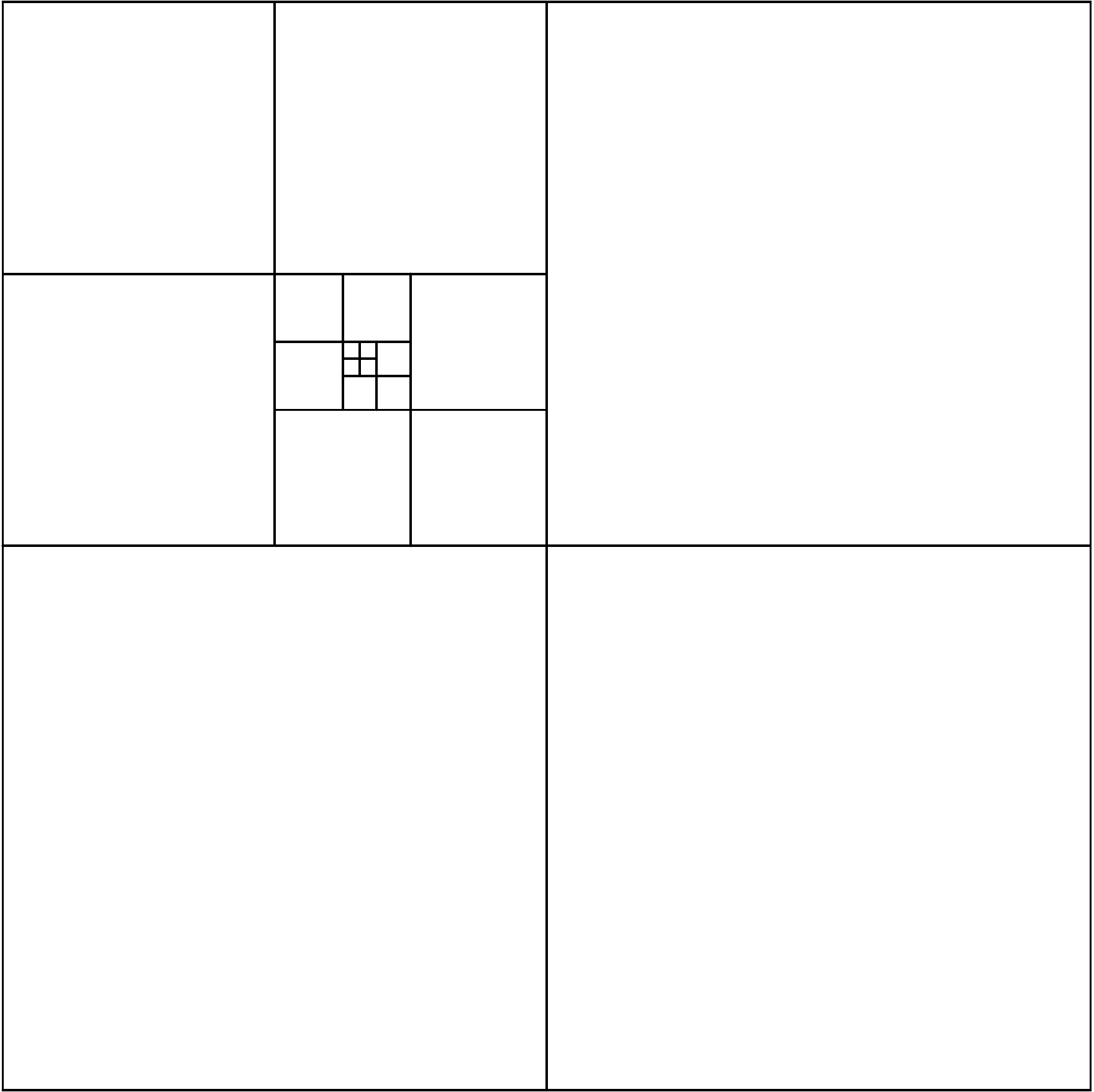}\hspace*{1.5cm}\includegraphics[height=0.35\linewidth]{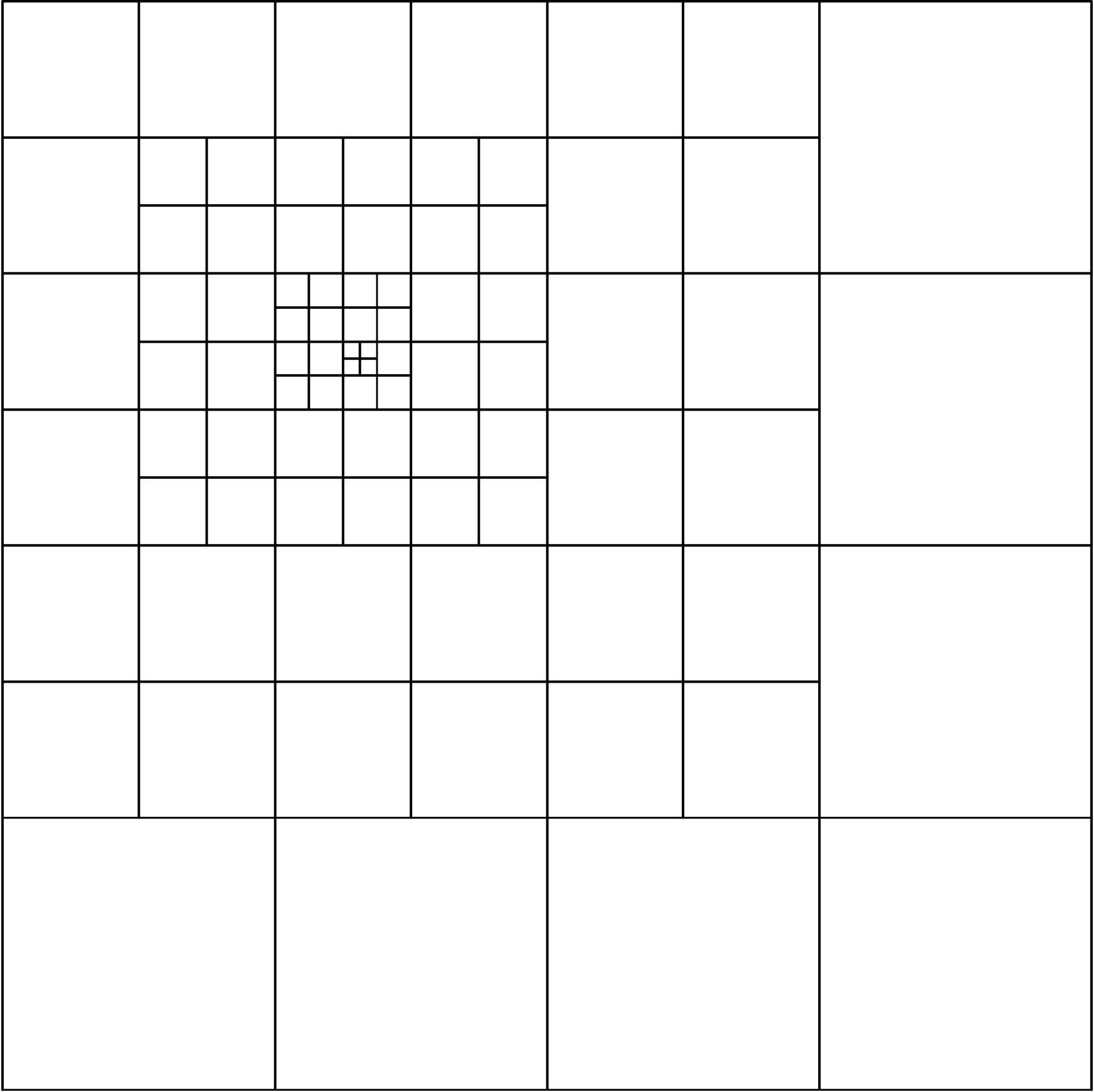}}
\caption{\sf Examples of non-regularized (left) and regularized (right) quadtree grids.\label{fig:2}}
\end{figure}

\bigskip
We note that in order to prevent complicated formulations and enhance the stability of the overall method, no cell can have both an adjacent
neighboring cell and a diagonally neighboring cell with a refinement level difference greater than one; see
\cite{Borthwick2000,Greaves1998,Popinet2003}. This condition is satisfied provided the quadtree is regularized; several regularization algorithms
can be found in \cite{Bern1999,Moore1995,Samet1982,Samet1985,Sankaranarayanan2007,Vaidya1989}. Examples of non-regularized and regularized quadtree
grids are shown in Figure \ref{fig:2}. 
	
\section{Adaptive well-balanced semi-discrete central-upwind scheme}\label{S:4}
In this section, we present an adaptive well-balanced semi-discrete central-upwind scheme for the system \eqref{eq:SWE}, which can be
written in the following vector form:
\begin{equation}
\bm{U}_t+\bm{F}(\bm{U},B)_x+\bm{G}(\bm{U},B)_y=\bm{S}(\bm{U},B),
\label{3.1}          
\end{equation}
where
\begin{equation*}
\bm{U}:=(w,hu,hv)^\top,
\end{equation*}
and the fluxes and source term are:
\begin{align}
&\bm{F}(\bm{U},B)=\left(hu,\dfrac{(hu)^2}{w-B}+\dfrac{g}{2}(w-B)^2,\dfrac{(hu)(hv)}{w-B}\right)^\top,\label{3.2}\\
&\bm{G}(\bm{U},B)=\left(hv,\dfrac{(hu)(hv)}{w-B},\dfrac{(hv)^2}{w-B}+\dfrac{g}{2}(w-B)^2\right)^\top,\label{3.3}\\ 
&\bm{S}(\bm{U},B)=\left(0,-g(w-B)B_x,-g(w-B)B_y\right)^\top.\label{3.4}
\end{align}

The central-upwind quadtree scheme will be designed according to the following algorithm:

\vskip5pt
\noindent\textbf{Step 1.} Generate a non-regularized grid with the seeding points (\S\ref{S:2}).

\vskip4pt
\noindent\textbf{Step 2.} Regularize the non-regularized grid (\S\ref{S:2}).

\vskip4pt
\noindent\textbf{Step 3.} Perform piecewise polynomial reconstructions and obtain the required point values of the bottom topography $B$
(\S\ref{S:4.5}) and conservative quantities $\bm{U}$ (\S\ref{S:4.4}).

\vskip4pt
\noindent\textbf{Step 4.} Calculate the local speeds (\S\ref{S:4.3}) and central-upwind numerical fluxes (\S\ref{S:4.2}).

\vskip4pt
\noindent\textbf{Step 5.} Calculate the well-balanced discrete source term (\S\ref{S:4.8}).

\vskip4pt
\noindent\textbf{Step 6.} Calculate the size of timestep, which can guarantee the positivity and stability (\S\ref{S:4.7}).

\vskip4pt
\noindent\textbf{Step 7.} Calculate local gradients in each cell, which are needed to determine the seeding points at the next timestep
(\S\ref{S:4.10}).

\vskip4pt
\noindent\textbf{Step 8.} Calculate current conservative quantities, which are going to be used as previous timestep data in
the construction of the new quadtree grid (\S\ref{S:4.10}).

\vskip4pt
\noindent\textbf{Step 9.} Evolve the solution by solving the time-dependent system of ODEs, obtained after the semi-discretization of the
studied SWEs over the quadtree grid.

\subsection{Finite-volume semi-discretization over quadtree grids}\label{S:4.1}
Let us consider a typical finite-volume Cartesian cell $C_{j,k}$ of size $\dx_{j,k}\times\dy_{j,k}$ centered at $(x_{j,k},y_{j,k})$. We
assume that at a certain time level $t$, the computed solution is available and represented in terms of its cell averages:
\begin{equation}
\xbar{\bm{U}}_{j,k}(t)\approx
\frac{1}{\dx_{j,k}\dy_{j,k}}\int\limits_{x_\jmh}^{x_\jph}\int\limits_{y_\kmh}^{y_\kph}\bm{U}(x,y,t)\,{\rm d}y\,{\rm d}x,
\label{3.5}          
\end{equation}
where $x_{j\pm\hf}:=x_{j,k}\pm\dx_{j,k}/2$ and $y_{k\pm\hf}:=y_{j,k}\pm\dy_{j,k}/2$.

\begin{figure}[ht!]
\centerline{\subcaptionbox{}{\includegraphics[height=1.5in]{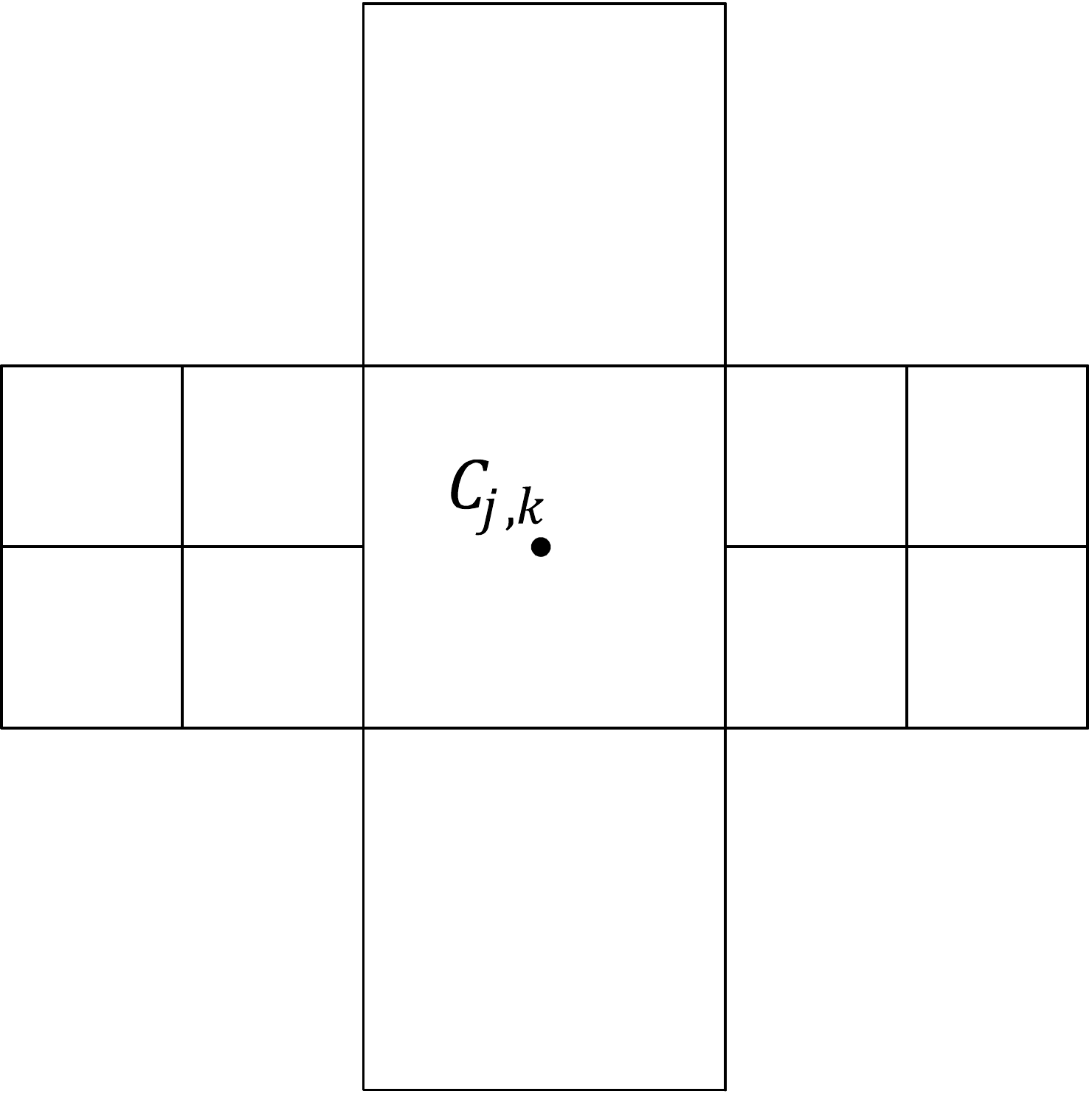}}\hspace{4.5em}
\subcaptionbox{}{\includegraphics[height=1.5in]{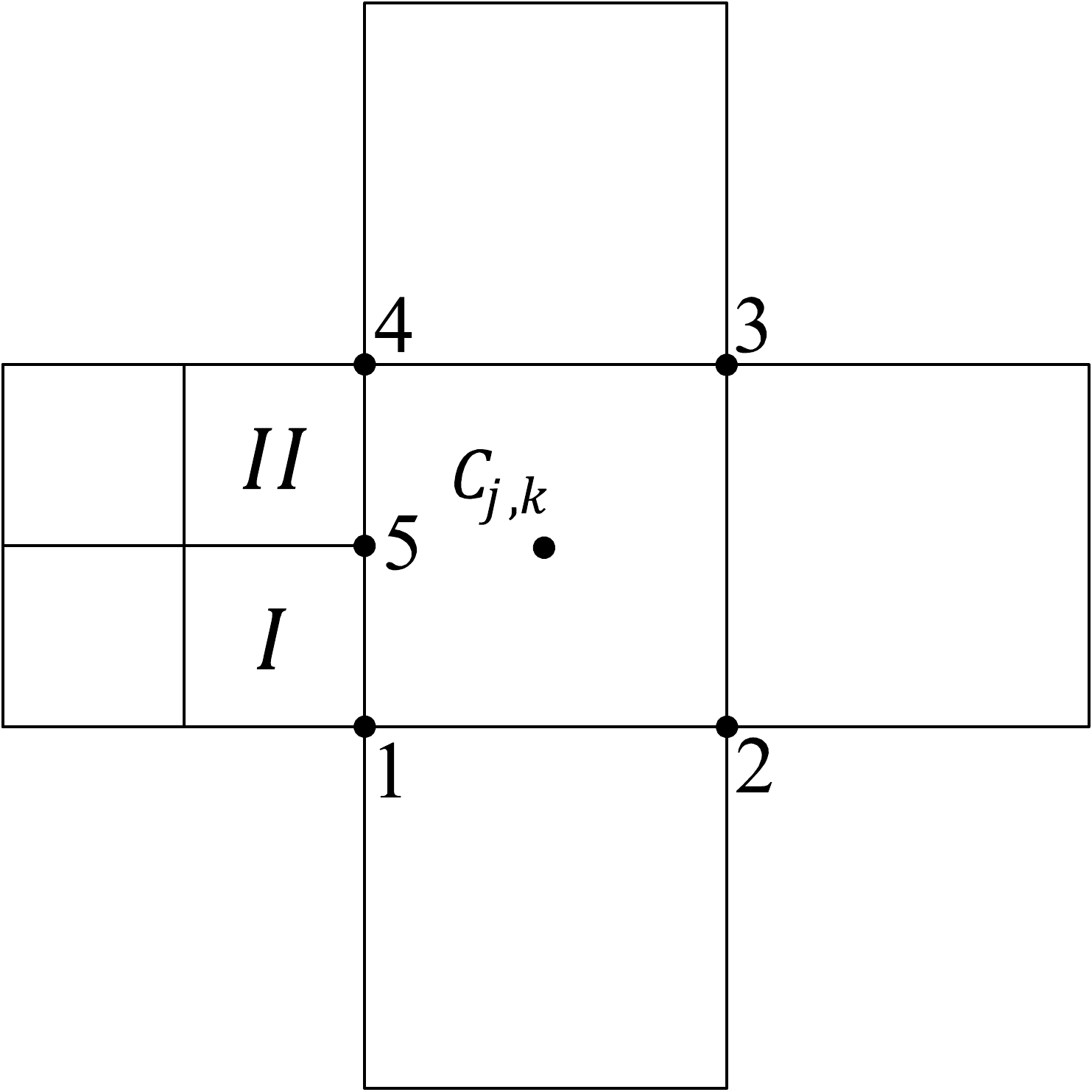}}\hspace{4.5em}\subcaptionbox{}{\includegraphics[height=1.5in]{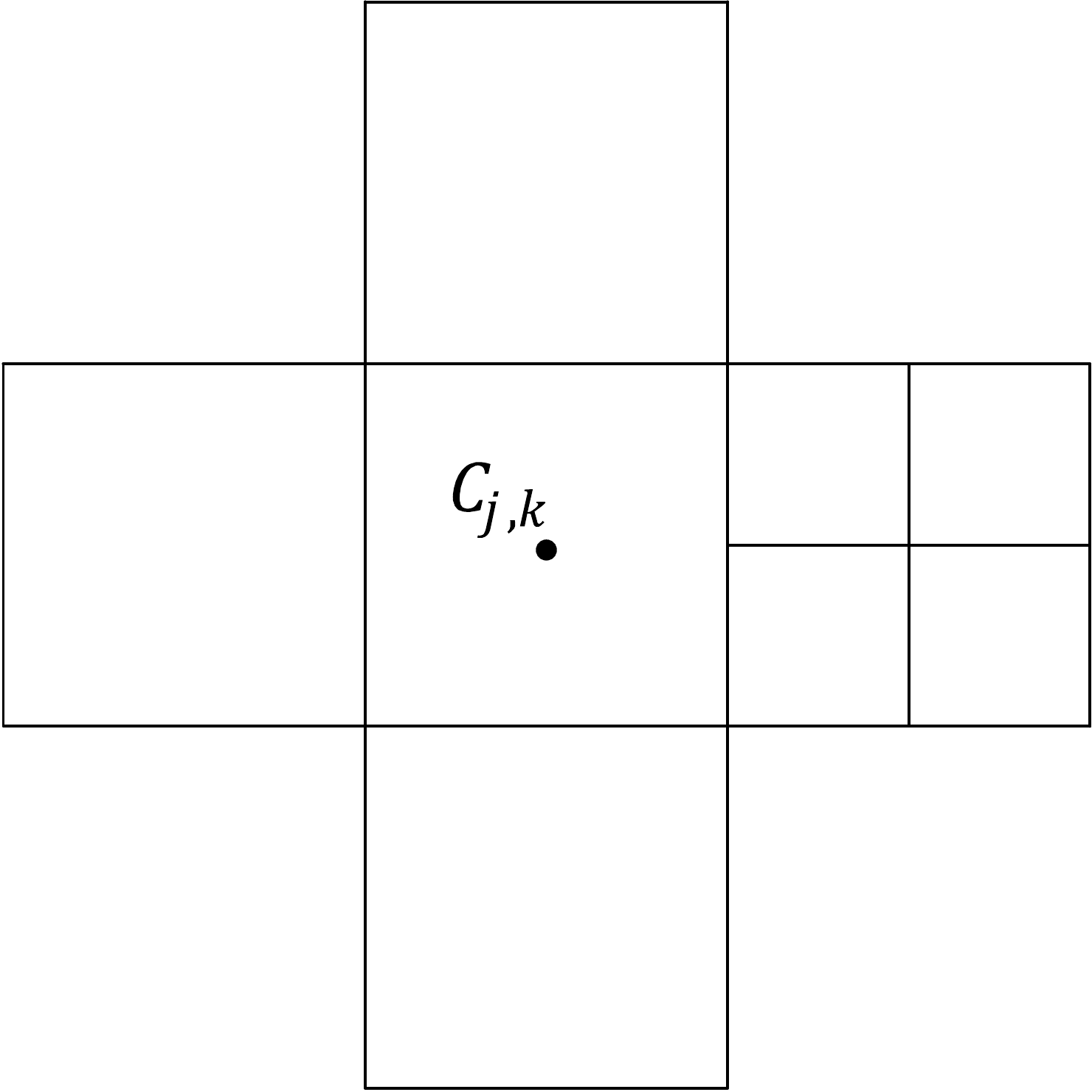}}}
\bigskip
\centerline{\subcaptionbox{}{\includegraphics[height=1.4in]{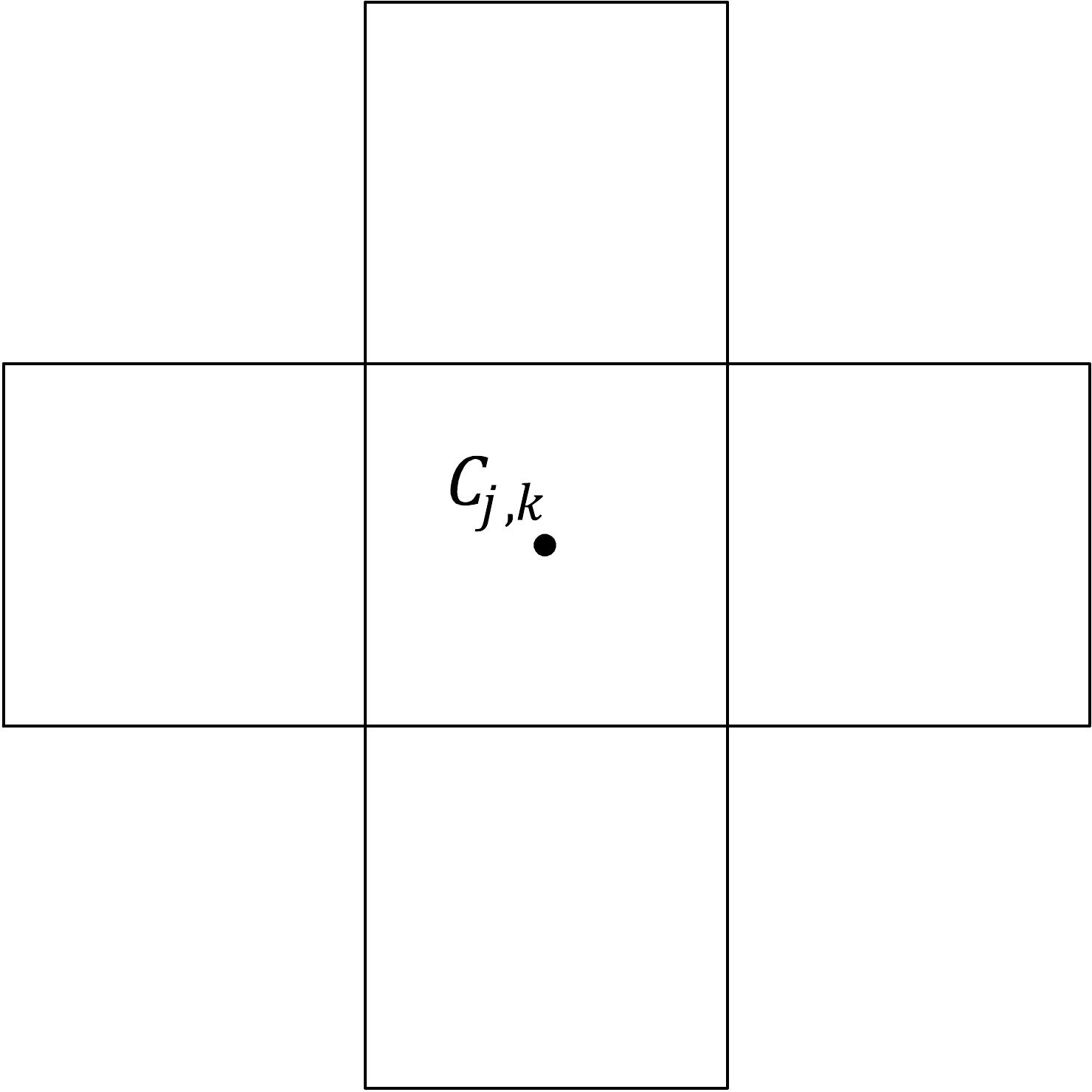}}\hspace{2.5em}
\subcaptionbox{}{\includegraphics[height=1.4in]{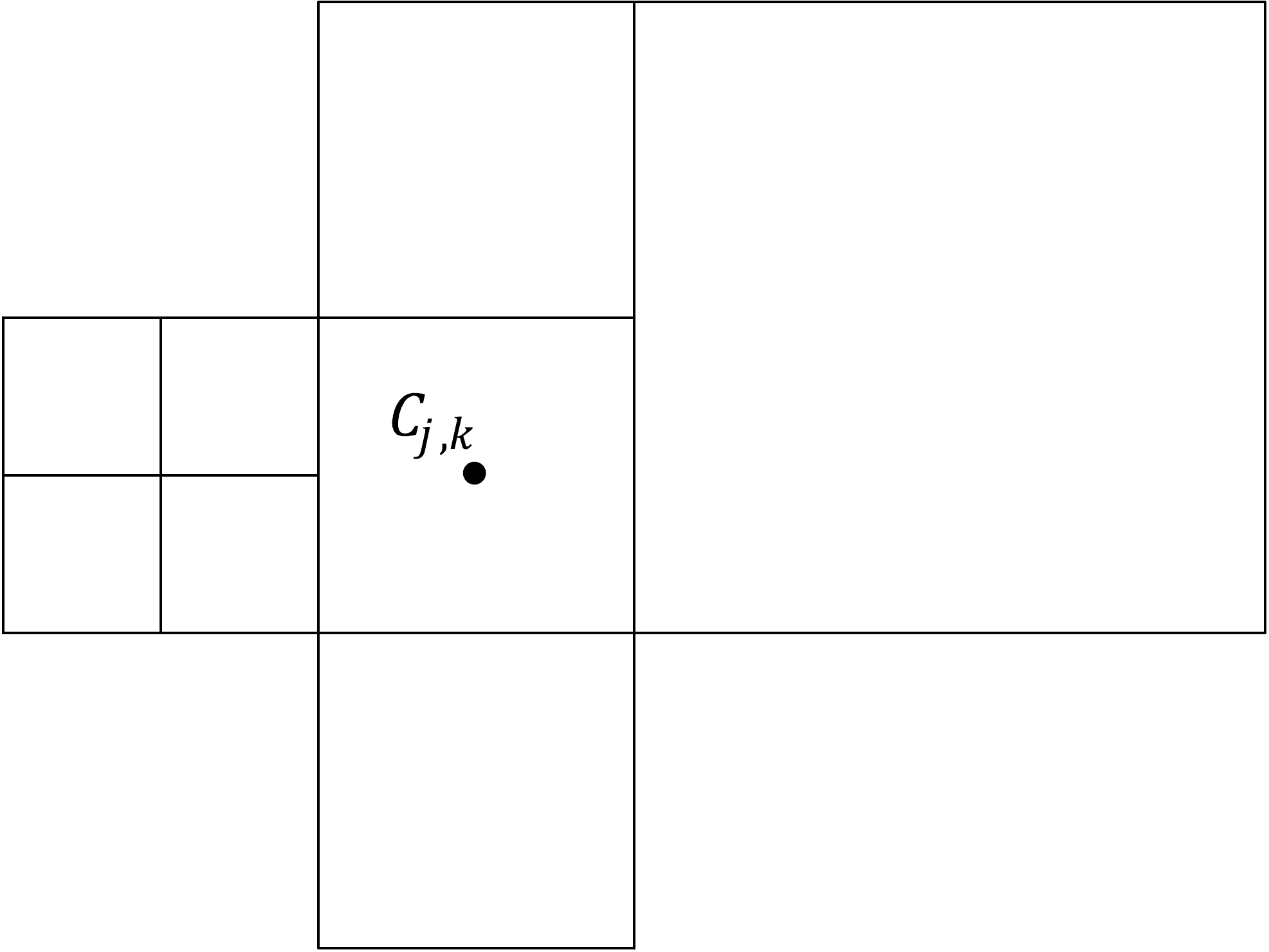}}\hspace{2.5em}\subcaptionbox{}{\includegraphics[height=1.4in]{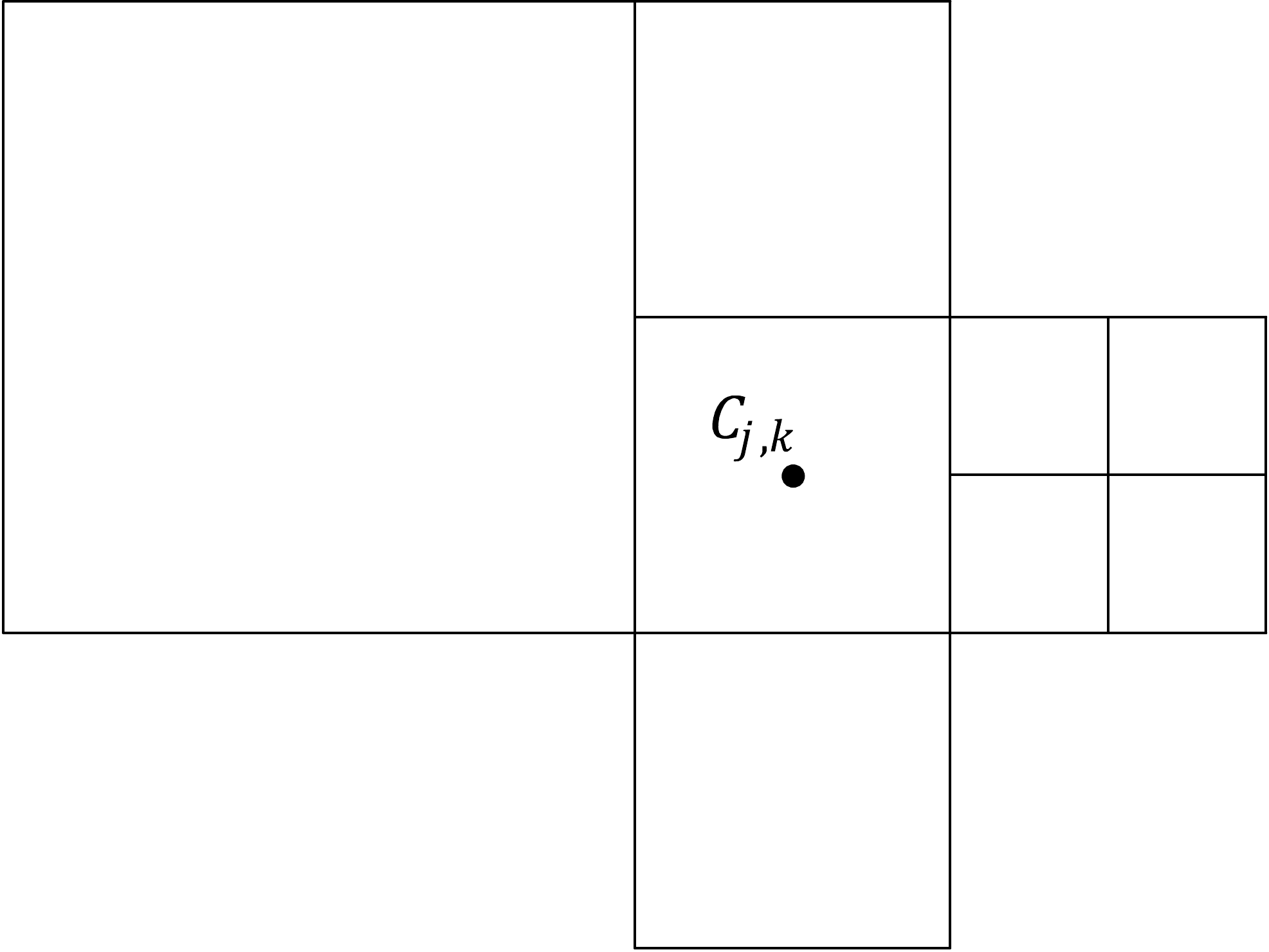}}}
\bigskip
\centerline{\subcaptionbox{}{\includegraphics[height=1.28in]{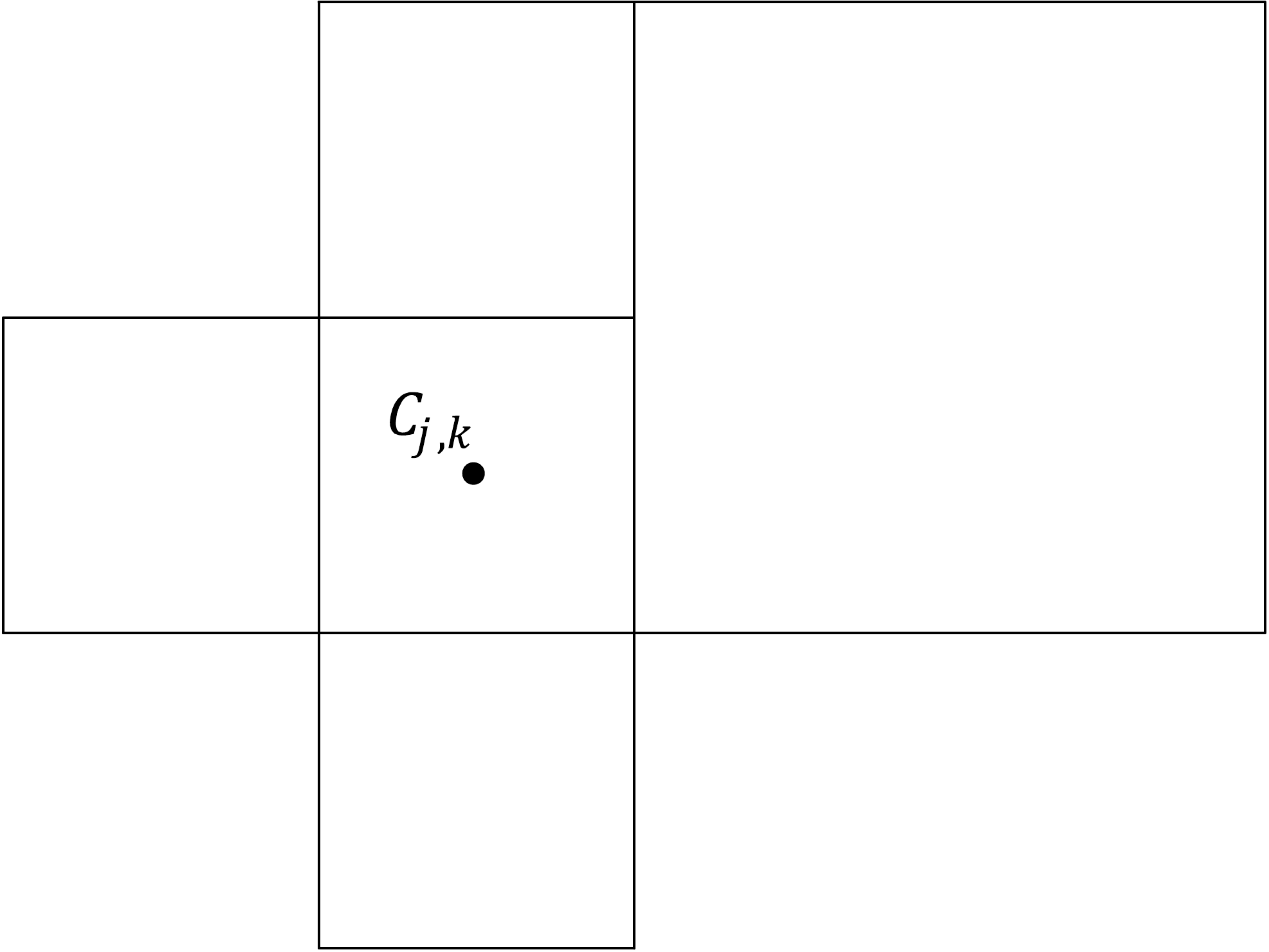}}\hspace{1em}
\subcaptionbox{}{\includegraphics[height=1.28in]{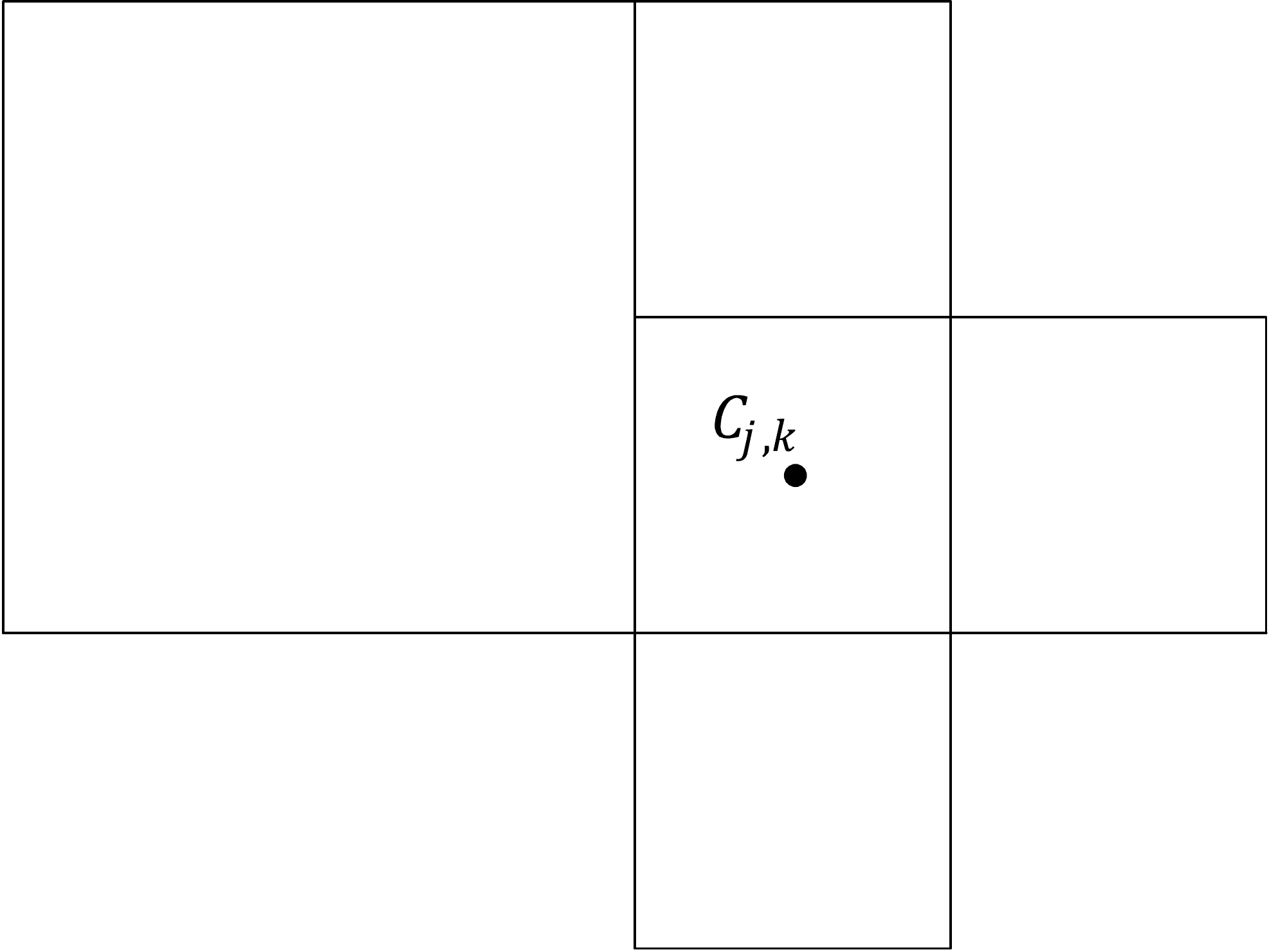}}\hspace{1em}\subcaptionbox{}{\includegraphics[height=1.28in]{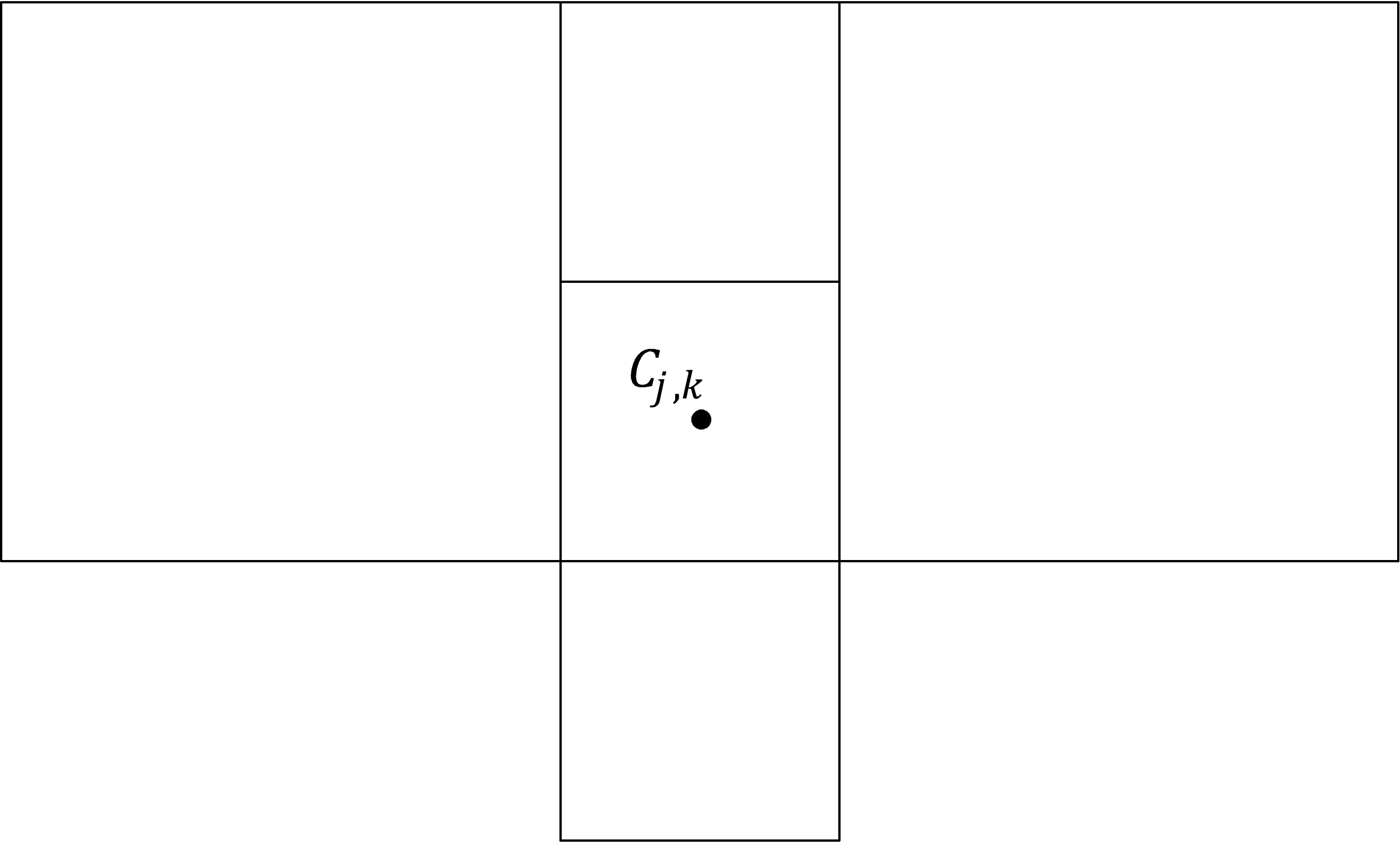}}}
\caption{\sf Permutations of neighboring cells of $C_{j,k}$ in the $x$-direction: (a)--(h) Possible configurations in regularized quadtree
grids; (i) An impossible configuration in regularized quadtree grids.\label{fig:3}}
\end{figure}

Considering the right and left neighbors of cell $C_{j,k}$, there exist nine different permutations of those neighboring cells; see Figure
\ref{fig:3}. We note, however, that only eight of them (configurations (a)--(h) in Figure \ref{fig:3}) are possible in the proposed
regularized quadtree grid. Similar cases are to be considered with respect to the neighboring cells on the top and bottom.

For the sake of brevity, we only present the quadtree scheme for configuration (b) in Figure \ref{fig:3} as an example of a quadtree cell
(other configurations can be treated in a similar manner). We denote left-neighboring cells of $C_{j,k}$ by $\Romannum{1}$ and
$\Romannum{2}$. These two cells centered at $(x_{j,k}-3\dx_{j,k}/4,y_{j,k}\pm\dy_{j,k}/4)$ are of size $\dx_{j,k}/2\times\dy_{j,k}/2$.

The cell averages $\,\xbar{\bm{U}}_{j,k}$ are evolved in time by solving the following system of time-dependent ODEs:
\begin{equation}
\frac{{\rm d}}{{\rm d}t}\,\xbar{\bm{U}}_{j,k}=
-\frac{\bm{H}_{\jph,k}^x-\dfrac{\bm{H}_{\jmh,k-\frac{1}{4}}^x+\bm{H}_{\jmh,k+\frac{1}{4}}^x}{2}}{\dx_{j,k}}
-\frac{\bm{H}_{j,\kph}^y-\bm{H}_{j,\kmh}^y}{\dy_{j,k}}+\,\xbar{\bm{S}}_{j,k},
\label{3.6}          
\end{equation}
obtained after the semi-discretization of the system \eqref{3.1}--\eqref{3.4}. In \eqref{3.6}, $\bm{H}_{\jph,k}^x$,
$\bm{H}_{\jmh,k\pm\frac{1}{4}}^x$, $\bm{H}_{j,\kph}^y$ and $\bm{H}_{j,\kmh}^y$ are the numerical fluxes, which, in general, are
\begin{equation}
\bm{H}^x_{\alpha,\beta}=\bm{H}^x(\bm{U}^-_{\alpha,\beta},\bm{U}^+_{\alpha,\beta};B_{\alpha,\beta})\quad\mbox{and}\quad
\bm{H}^y_{\gamma,\delta}=\bm{H}^y(\bm{U}^-_{\gamma,\delta},\bm{U}^+_{\gamma,\delta};B_{\gamma,\delta})
\label{3.7}
\end{equation}
where
\begin{equation}
\bm{U}^\pm_{\alpha,\beta}=\lim_{x\to x_\alpha\pm0}\widetilde{\bm{U}}(x,y_\beta)\quad\mbox{and}\quad
\bm{U}^\pm_{\gamma,\delta}=\lim_{y\to y_\delta\pm0}\widetilde{\bm{U}}(x_\gamma,y),
\label{3.8}
\end{equation}
and $\widetilde{\bm{U}}$ is a piecewise polynomial interpolation. Second-order schemes employ piecewise linear interpolations,
\begin{equation}
\widetilde{\bm{U}}(x,y)=(\bm{U}_x)_{j,k}[x-x_j]+(\bm{U}_y)_{j,k}[y-y_k],\quad(x,y)\in C_{j,k},
\label{3.8a}
\end{equation}
where the slopes $(\bm{U}_x)_{j,k}$ and $(\bm{U}_y)_{j,k}$ are yet to be determined; see \S\ref{S:4.4} below. Finally,
$\,\xbar{\bm{S}}_{j,k}$ is a cell average of the source term:
\begin{equation}
\xbar{\bm{S}}_{j,k}\approx
\frac{1}{\dx_{j,k}\dy_{j,k}}\int\limits_{x_\jmh}^{x_\jph}\int\limits_{y_\kmh}^{y_\kph}\bm{S}(\bm{U},B)\,{\rm d}y\,{\rm d}x.
\label{3.9}          
\end{equation}

We note that all of the indexed quantities in \eqref{3.6}--\eqref{3.9} are time-dependent, but from now on we omit this dependence for the
sake of brevity.

\subsection{Piecewise bilinear reconstruction of $B$}\label{S:4.5}
The quadtree grid consist of cells of different sizes: $\dx\times\dy$,
$\frac{\dx}{2}\times\frac{\dy}{2},\ldots,\frac{\dx}{2^{m-1}}\times\frac{\dy}{2^{m-1}}$. We denote the set of cells of the corresponding size
by ${\cal C}^{(\ell)}$, that is, ${\cal C}^{(\ell)}=\{C_{j,k}: |C_{j,k}|=\frac{\dx}{2^{\ell-1}}\times\frac{\dy}{2^{\ell-1}}\}$.

We follow the lines of \cite{Kurganov2007} and use a continuous piecewise bilinear reconstruction of the bottom topography
$\widetilde B(x,y)$. We note that on a quadtree grid, the approach from \cite{Kurganov2007} does not directly apply to a quadtree grid since
it contains cells, whose vertex is a midpoint of the edge of the neighboring cell as point 5 in the configuration considered in Figure
\ref{fig:3} (b). We therefore propose the following algorithm for constructing $\widetilde B$.

\vskip5pt
\noindent\textbf{Step 1.} Set $\ell:=1$.

\vskip4pt
\noindent\textbf{Step 2.} Reconstruct bilinear pieces $\widetilde B$ for all $(j,k)$ such that $C_{j,k}\in{\cal C}^{(\ell)}$. This is done as
follows.

We first obtain the point values of $B$ at the vertices of $C_{j,k}$. If one of them lies on the edge of a larger neighboring cell (which
belongs to the set ${\cal C}^{(\ell-1)}$ where the bilinear piece has already been constructed), then there are two possibilities:

(i) either this vertex coincides with a vertex of the neighboring cell and then the point value of $B$ has been already computed there;

(ii) or this vertex is a midpoint of the edge of the neighboring cell and then the point value of $B$ at this vertex is an average of the
point values of $B$ at those two vertices of the neighboring cell that lie on the same edge (for example, in the configuration considered in
Figure \ref{fig:3} (b), the value of $B$ at point 5 will be equal to the average of the values of $B$ at points 1 and 4).

Otherwise, we proceed as in \cite{Kurganov2007} and set
\begin{equation*}
\begin{aligned}
B_{j\pm\hf,k\pm\hf}:=\hf\Big(&\max_{\xi^2+\eta^2=1}\lim_{\ell_x,\ell_y\to0}B(x_{j\pm\hf}+\ell_x\xi,y_{k\pm\hf}+\ell_y\eta)\\
+&\min_{\xi^2+\eta^2=1}\lim_{\ell_x,\ell_y\to0}B(x_{j\pm\hf}+\ell_x\xi,y_{k\pm\hf}+\ell_y\eta)\Big),
\end{aligned}
\end{equation*}
which reduces to $B_{j\pm\hf,k\pm\hf}=B(x_{j\pm\hf},y_{k\pm\hf})$ if the function $B$ is continuous at $(x_{j\pm\hf},y_{k\pm\hf})$.

Equipped with the point values $B_{j\pm\hf,k\pm\hf}$, we construct the following bilinear piece in cell $C_{j,k}\in{\cal C}^{(\ell)}$:
$$
\begin{aligned}
\widetilde B(x,y)=&\,B_{\jmh,\kmh}+(B_{\jph,\kmh}-B_{\jmh,\kmh})\,\frac{x-x_\jmh}{\nicefrac{\dx}{2^{\ell-1}}}+
(B_{\jmh,\kph}-B_{\jmh,\kmh})\,\frac{y-y_\kmh}{\nicefrac{\dy}{2^{\ell-1}}}\\
+&\,(B_{\jph,\kph}-B_{\jph,\kmh}-B_{\jmh,\kph}+B_{\jmh,\kmh})\,\frac{x-x_\jmh}{\nicefrac{\dx}{2^{\ell-1}}}\cdot
\frac{y-y_\kmh}{\nicefrac{\dy}{2^{\ell-1}}},~~(x,y)\in C_{j,k}.
\end{aligned}
$$

\vskip4pt
\noindent\textbf{Step 3.} Set $\ell:=\ell+1$.

\vskip4pt
\noindent\textbf{Step 4.} If $\ell\le m$, then go to Step 2.

\medskip
Note that the restriction of the interpolant $\widetilde B$ along each of the cell is a linear function and the cell average of
$\widetilde B$ over the cell $C_{j,k}$ is equal to its value at the center of the cell and is also equal to the average of the values of
$\widetilde B$ at the midpoints of the edges of $C_{j,k}$, namely, we have
\begin{equation}
B_{j,k}:=\widetilde B(x_j,y_k)=\frac{1}{\dx_{j,k}\dy_{j,k}}\iint\limits_{C_{j,k}}\widetilde B(x,y)\,{\rm d}x\,{\rm d}y=\frac{B_{\jph,k}+
B_{\jmh,k}+B_{j,\kph}+B_{j,\kmh}}{4},
\label{3.12a}
\end{equation}
where
\begin{equation}
B_{j\pm\hf,k}:=\widetilde B(x_{j\pm\hf},y_k)=\hf\left(B_{j\pm\hf,\kph}+B_{j\pm\hf,\kmh}\right),
\label{3.13a}
\end{equation}
and
\begin{equation}
B_{j,k\pm\hf}:=\widetilde B(x_j,y_{k\pm\hf})=\hf\left(B_{\jph,k\pm\hf}+B_{\jmh,k\pm\hf}\right).
\label{3.14a}
\end{equation}
We also note that the values of $B$ at the midpoints of the right edges of cells $I$ and $II$ in configuration considered in Figure
\ref{fig:3} (b) can be obtained in a similar way:
\begin{equation}
\begin{aligned}
&B_{\jmh,k\pm\frac{1}{4}}:=\widetilde B(x_\jmh,y_{k\pm\frac{1}{4}})=\hf\left(B_{\jmh,k}+B_{\jmh,k\pm\hf}\right).
\end{aligned}
\label{3.15a}
\end{equation}

Formulae \eqref{3.12a}--\eqref{3.15a} are crucial for the proof of the positivity preserving property of our well-balanced quadtree
central-upwind scheme; see \S\ref{S:4.7}.

\smallskip
{\em Remark.} We note that the proposed piecewise bilinear reconstruction can be applied to discontinuous bottom functions $B(x,y)$.

\subsection{Piecewise linear reconstruction of $\bm{U}$}\label{S:4.4}
In this paper, we design a second-order scheme, which employs a piecewise linear reconstruction $\widetilde{\bm{U}}$ in each cell. We
then obtain the point values of $\bm{U}$ (required in \eqref{3.7}) using \eqref{3.8}, \eqref{3.8a}, which for cell $C_{j,k}$ from Figure
\ref{fig:3} (b) results in
\begin{equation}
\begin{aligned}
&\bm{U}_{\jph,k}^+=\,\xbar{\bm{U}}_{j+1,k}-\frac{\dx_{j+1,k}}{2}(\bm{U}_x)_{j+1,k}\,,\quad
\bm{U}_{\jph,k}^-=\,\xbar{\bm{U}}_{j,k}+\frac{\dx_{j,k}}{2}(\bm{U}_x)_{j,k}\,,\\
&\bm{U}_{\jmh,k\pm\frac{1}{4}}^+=\,\xbar{\bm{U}}_{j,k}-\frac{\dx_{j,k}}{2}(\bm{U}_x)_{j,k}\pm\frac{\dy_{j,k}}{2}(\bm{U}_y)_{j,k}\,,\\
&\bm{U}_{\jmh,k\pm\frac{1}{4}}^-=\,\xbar{\bm{U}}_{j-\frac{1}{4},k\pm\frac{1}{4}}+\frac{\dx_{j,k}}{4}(\bm{U}_x)_{j-\frac{1}{4},k\pm\frac{1}{4}}
\end{aligned}
\label{3.10} 
\end{equation}
where $\,\xbar{\bm{U}}_{j-\frac{1}{4},k-\frac{1}{4}}$ and $\,\xbar{\bm{U}}_{j-\frac{1}{4},k+\frac{1}{4}}$ denote the cell averages of $\bm{U}$ over the cells $I$ and $II$,
respectively.

In order to achieve the formal second order of accuracy, the slopes $(\bm{U}_x)$ and $(\bm{U}_y)$ in \eqref{3.10} are to be at least
first-order approximations of the corresponding derivatives. In order to minimize oscillations, we compute the slopes using the minmod
limiter (see, e.g., \cite{Berger2005,Kuzmin2006,Sweby1984}), which is implemented in the following way:
\begin{equation}
\begin{aligned}
&(\bm{U}_x)_{j,k}={\rm minmod}\left(\frac{\,\xbar{\bm{U}}_{j,k}-\,\xbar{\bm{U}}_{j-\frac{1}{4},k-\frac{1}{4}}}{3\dx_{j,k}/4},\,
\frac{\,\xbar{\bm{U}}_{j,k}-\,\xbar{\bm{U}}_{j-\frac{1}{4},k+\frac{1}{4}}}{3\dx_{j,k}/4},\,
\frac{\,\xbar{\bm{U}}_{j+1,k}-\,\xbar{\bm{U}}_{j,k}}{\dx_{j,k}}\right),\\
&(\bm{U}_y)_{j,k}={\rm minmod}\left(\frac{\,\xbar{\bm{U}}_{j,k}-\,\xbar{\bm{U}}_{j,k-1}}{\dy_{j,k}},\,
\frac{\,\xbar{\bm{U}}_{j,k+1}-\,\xbar{\bm{U}}_{j,k}}{\dy_{j,k}}\right)
\end{aligned}
\label{3.11} 
\end{equation}
where the minmod function is defined by
\begin{equation*}
min\{z_1,z_2,...\}:=\left\{
\begin{aligned}
&min_j\{z_j\},&&\text{if}~z_j>0\quad\forall{j},\\
&max_j\{z_j\},&&\text{if}~z_j<0\quad\forall{j},\\
&0,&&\text{otherwise}.
\end{aligned}\right.
\end{equation*}

\subsubsection{Positivity preserving correction of $w$}\label{S:4.6}
The piecewise linear reconstruction \eqref{3.11} cannot guarantee the non-negativity of
$$
\widetilde h(x,y):=\widetilde w(x,y)-\widetilde B(x,y).
$$
In fact, for the configuration considered in Figure \ref{fig:3} (b), we only need the following five inequalities to be satisfied: 
$$
\begin{aligned}
&h_{\jmh,k\pm\frac{1}{4}}^+=w_{\jmh,k\pm\frac{1}{4}}^+-B_{\jmh,k\pm\frac{1}{4}}\ge0,\quad h_{\jph,k}^-=w_{\jph,k}^--B_{\jph,k}\ge0,\\
&h_{j,\kmh}^+=w_{j,\kmh}^+-B_{j,\kmh}\ge0,\quad\mbox{and}\quad h_{j,\kph}^-=w_{j,\kph}^--B_{j,\kph}\ge0.
\end{aligned}
$$
If at least one of these inequalities is not satisfied, we need to correct $\widetilde w$ in the cell $C_{j,k}$. We note that the correction
used in \cite{Kurganov2007} will not in general work on quadtree grids. We therefore propose an alternative correction procedure and replace
the linear pieces in the problematic cells with the bilinear one denoted by $\mathring w(x,y)$ and constructed as follows. Let us denote by
$$
\begin{aligned}
&w_{j,k}^{\rm NE}:=\widetilde w(x_\jph-0,y_\kph-0),\quad w_{j,k}^{\rm SE}:=\widetilde w(x_\jph-0,y_\kmh+0),\\
&w_{j,k}^{\rm NW}:=\widetilde w(x_\jmh+0,y_\kph-0),\quad w_{j,k}^{\rm SW}:=\widetilde w(x_\jmh+0,y_\kmh+0)
\end{aligned}
$$
the four corner point values of the linear piece of $\widetilde w$ over the cell $C_{j,k}$. If
\begin{equation}
w_{j,k}^{\rm NE}\ge B_{\jph,\kph},\quad w_{j,k}^{\rm SE}\ge B_{\jph,\kmh},\quad w_{j,k}^{\rm NW}\ge B_{\jmh,\kph}\quad\mbox{and}\quad
w_{j,k}^{\rm SW}\ge B_{\jmh,\kmh},
\label{3.16a}
\end{equation}
then we set
\begin{equation}
\begin{aligned}
\mathring w(x,y)=&\,w_{j,k}^{\rm SW}+(w_{j,k}^{\rm SE}-w_{j,k}^{\rm SW})\,\frac{x-x_\jmh}{\dx_{j,k}}+
(w_{j,k}^{\rm NW}-w_{j,k}^{\rm SW})\,\frac{y-y_\kmh}{\dy_{j,k}}\\
+&\,(w_{j,k}^{\rm NE}-w_{j,k}^{\rm SE}-w_{j,k}^{\rm NW}+w_{j,k}^{\rm SW})\,\frac{x-x_\jmh}{\dx_{j,k}}\cdot\frac{y-y_\kmh}{\dy_{j,k}},~~
(x,y)\in C_{j,k}.
\end{aligned}
\label{3.17a}
\end{equation}

If at least one of the inequalities in \eqref{3.16a} is not satisfied, we would first need to correct the point values of $w$ at the vertices
of $C_{j,k}$. There are three different cases to be considered.

{\em Case 1: only one of the inequalities in \eqref{3.16a} is not satisfied}. Without loss of generality, we assume that
$w_{j,k}^{\rm NE}<B_{\jph,\kph}$. We then replace the point values of $w$ at the vertices of $C_{j,k}$ with
$$
\begin{aligned}
&\mathring w_{j,k}^{\rm NE}=B_{\jph,\kph},\quad\mathring w_{j,k}^{\rm SE}=B_{\jph,\kmh}+\frac{4}{3}(\,\xbar w_{j,k}-B_{j,k}),\\
&\mathring w_{j,k}^{\rm NW}=B_{\jmh,\kph}+\frac{4}{3}(\,\xbar w_{j,k}-B_{j,k}),\quad
\mathring w_{j,k}^{\rm SW}=B_{\jmh,\kmh}+\frac{4}{3}(\,\xbar w_{j,k}-B_{j,k}).
\end{aligned}
$$

{\em Case 2: only two of the inequalities in \eqref{3.16a} are not satisfied}. Without loss of generality, we assume that
$w_{j,k}^{\rm NE}<B_{\jph,\kph}$ and $w_{j,k}^{\rm SE}<B_{\jph,\kmh}$. We then replace the point values of $w$ at the vertices of $C_{j,k}$
with
$$
\begin{aligned}
&\mathring w_{j,k}^{\rm NE}=B_{\jph,\kph},\quad\mathring w_{j,k}^{\rm SE}=B_{\jph,\kmh},\\
&\mathring w_{j,k}^{\rm NW}=B_{\jmh,\kph}+2(\,\xbar w_{j,k}-B_{j,k}),\quad\mathring w_{j,k}^{\rm SW}=B_{\jmh,\kmh}+2(\,\xbar w_{j,k}-B_{j,k}).
\end{aligned}
$$

{\em Case 3: only three of the inequalities in \eqref{3.16a} are not satisfied}. Without loss of generality, we assume that
$w_{j,k}^{\rm NE}<B_{\jph,\kph}$, $w_{j,k}^{\rm SE}<B_{\jph,\kmh}$ and $w_{j,k}^{\rm NW}<B_{\jmh,\kph}$. We then replace the point values of
$w$ at the vertices of $C_{j,k}$ with
$$
\begin{aligned}
&\mathring w_{j,k}^{\rm NE}=B_{\jph,\kph},\quad\mathring w_{j,k}^{\rm SE}=B_{\jph,\kmh},\quad\mathring w_{j,k}^{\rm NW}=B_{\jmh,\kph},\\
&\mathring w_{j,k}^{\rm SW}=4\,\xbar w_{j,k}-B_{\jph,\kph}-B_{\jph,\kmh}-B_{\jmh,\kph}.
\end{aligned}
$$

In all of the above three cases, we use the corrected point values $\mathring w_{j,k}^{\rm NE}$, $\mathring w_{j,k}^{\rm SE}$,
$\mathring w_{j,k}^{\rm NW}$ and $\mathring w_{j,k}^{\rm SW}$ to construct the corrected bilinear approximant (compare with \eqref{3.17a})
$$
\begin{aligned}
\mathring w(x,y)=&\,\mathring w_{j,k}^{\rm SW}+(\mathring w_{j,k}^{\rm SE}-\mathring w_{j,k}^{\rm SW})\,
\frac{x-x_\jmh}{\dx_{j,k}}+(\mathring w_{j,k}^{\rm NW}-\mathring w_{j,k}^{\rm SW})\,\frac{y-y_\kmh}{\dy_{j,k}}\\
+&\,(\mathring w_{j,k}^{\rm NE}-\mathring w_{j,k}^{\rm SE}-\mathring w_{j,k}^{\rm NW}+\mathring w_{j,k}^{\rm SW})\,
\frac{x-x_\jmh}{\dx_{j,k}}\cdot\frac{y-y_\kmh}{\dy_{j,k}},~~(x,y)\in C_{j,k}.
\end{aligned}
$$

It is easy to show that the constructed bilinear piece $\mathring w(x,y)$ is conservative, that is,
$$
\frac{1}{\dx_{j,k}\dy_{j,k}}\int\limits_{C_{j,k}}\mathring w(x,y)\,{\rm d}y\,{\rm d}x=\,\xbar w_{j,k},
$$
and positivity preserving, that is,
$$
\mathring w(x,y)\ge\widetilde B(x,y),~~(x,y)\in C_{j,k}.
$$
We also notice that the point values of $w$ (required in \eqref{3.7}) at the cell $C_{j,k}$ from Figure \ref{fig:3} (b) are
$$
\begin{aligned}
&w_{\jmh,k\pm\frac{1}{4}}^+=\mathring w(x_\jmh+0,y_{k\pm\frac{1}{4}}),&& w_{\jph,k}^-=\mathring w(x_\jph-0,y_k),\\
&w_{j,\kmh}^+=\mathring w(x_j,y_\kmh+0),&& w_{j,\kph}^-=\mathring w(x_j,y_\kph-0)
\end{aligned}
$$
and thus the corresponding corrected values of $h$,
$$
\begin{aligned}
&h_{\jmh,k\pm\frac{1}{4}}^+=w_{\jmh,k\pm\frac{1}{4}}^+-B_{\jmh,k\pm\frac{1}{4}},&&h_{\jph,k}^-=w_{\jph,k}^--B_{\jph,k},\\
&h_{j,\kmh}^+=w_{j,\kmh}^+-B_{j,\kmh},&&h_{j,\kph}^-=w_{j,\kph}^--B_{j,\kph},
\end{aligned}
$$
are nonnegative.

Finally, we would like to point out that the values of $h$ at the boundaries of cell $C_{j,k}$ may be very small or even zero. This will
require the computation of the corresponding point values of $u$ and $v$ to be desingularized. We use the desingularization approach from
\cite{Kurganov2007}:
\begin{equation*}
u:=\frac{\sqrt{2}h(hu)}{\sqrt{h^4+\max\{h^4,\varepsilon\}}},\qquad v:=\frac{\sqrt{2}h(hv)}{\sqrt{h^4+\max\{h^4,\varepsilon\}}},
\end{equation*}
where we take $\varepsilon=\max\{\min_{j,k}\{(\dx_{j,k})^4\},\min_{j,k}\{(\dy_{j,k})^4\}\}$. After recomputing the point values of $h$, $u$
and $v$, the $x$- and $y$-discharges are also recalculated by setting:
\begin{equation*}
(hu):=h\cdot u,\qquad (hv):=h\cdot v.
\end{equation*}
Note that in the above two equations, we have omitted all of the indices for the sake of brevity.

\subsection{Local speeds}\label{S:4.3}
The one-sided local speeds of propagation, denoted at the corresponding cell interfaces by $a_{\alpha,\beta}^\pm$ and
$b_{\gamma,\delta}^\pm$, are calculated using the largest and smallest eigenvalues of the Jacobian matrices
$\frac{\partial\bm{F}}{\partial\bm{U}}$ and $\frac{\partial\bm{G}}{\partial\bm{U}}$ and can be estimated by
\begin{equation}
\begin{aligned}
&a_{\alpha,\beta}^+=\max\left\{u_{\alpha,\beta}^++\sqrt{gh_{\alpha,\beta}^+},\,u_{\alpha,\beta}^-+\sqrt{gh_{\alpha,\beta}^-},\,0\right\},\\
&a_{\alpha,\beta}^-=\min\left\{u_{\alpha,\beta}^+-\sqrt{gh_{\alpha,\beta}^+},\,u_{\alpha,\beta}^--\sqrt{gh_{\alpha,\beta}^-},\,0\right\},\\
&b_{\gamma,\delta}^+=\max\left\{v_{\gamma,\delta}^++\sqrt{gh_{\gamma,\delta}^+},\,v_{\gamma,\delta}^-+\sqrt{gh_{\gamma,\delta}^-},\,
0\right\},\\
&b_{\gamma,\delta}^-=\min\left\{v_{\gamma,\delta}^+-\sqrt{gh_{\gamma,\delta}^+},\,v_{\gamma,\delta}^--\sqrt{gh_{\gamma,\delta}^-},\,
0\right\}.
\end{aligned}
\label{3.12}
\end{equation}

\subsection{Central-upwind numerical fluxes}\label{S:4.2}
We use the central-upwind fluxes from \cite{Kurganov2007}, originally derived in \cite{KTrp}:
\begin{equation}
\begin{aligned}
\bm{H}_{\alpha,\beta}^x=&\frac{a_{\alpha,\beta}^+\bm{F}(\bm{U}_{\alpha,\beta}^-,B_{\alpha,\beta})-
a_{\alpha,\beta}^-\bm{F}(\bm{U}_{\alpha,\beta}^+,B_{\alpha,\beta})}{a_{\alpha,\beta}^+-a_{\alpha,\beta}^-}
+\frac{a_{\alpha,\beta}^+a_{\alpha,\beta}^-}{a_{\alpha,\beta}^+-a_{\alpha,\beta}^-}\left[\bm{U}_{\alpha,\beta}^+-\bm{U}_{\alpha,\beta}^-
\right],\\
\bm{H}_{\gamma,\delta}^y=&\frac{b_{\gamma,\delta}^+\bm{G}(\bm{U}_{\gamma,\delta}^-,B_{\gamma,\delta})-
b_{\gamma,\delta}^-\bm{G}(\bm{U}_{\gamma,\delta}^+,B_{\gamma,\delta})}{b_{\gamma,\delta}^+-b_{\gamma,\delta}^-}
+\frac{b_{\gamma,\delta}^+b_{\gamma,\delta}^-}{b_{\gamma,\delta}^+-b_{\gamma,\delta}^-}
\left[\bm{U}_{\gamma,\delta}^+-\bm{U}_{\gamma,\delta}^-\right],
\end{aligned}
\label{3.13} 
\end{equation}
where $(\alpha,\beta)\in\big\{(\jmh,k-\frac{1}{4}),(\jmh,k+\frac{1}{4}),(\jph,k)\big\}$ and
$(\gamma,\delta)\in\big\{(j,\kmh),(j,\kph)\big\}$ in the configuration considered in Figure \ref{fig:3} (b).

\subsection{Well-balanced discretization of the source term}\label{S:4.8}
A numerical scheme is well-balanced when the discretized cell average of the source term,
$\,\xbar{\bm{S}}_{j,k}=\big(0,\,\xbar S_{j,k}^{\,(2)},\,\xbar S_{j,k}^{\,(3)}\big)^\top$, exactly balances the numerical fluxes in equation
\eqref{3.6} at the ``lake-at-rest'' steady state \eqref{1.2}, that is, when the right-hand side (RHS) of \eqref{3.6} vanishes as long as
$\,\xbar{\bm{U}}_{j,k}\equiv\left(\widehat w,0,0\right)^\top$ for all $(j,k)$, where $\widehat w$ is a constant.

We note that at the ``lake-at-rest data'', all of the reconstructed point values are $w^\pm=\widetilde w$ and $u^\pm=v^\pm=0$ and thus,
$a^+_{\alpha,\beta}=-a^-_{\alpha,\beta},\,\forall(\alpha,\beta)$, $\,b^+_{\gamma,\delta}=-b^-_{\gamma,\delta}\,\forall(\gamma,\delta)$, and
the numerical fluxes \eqref{3.13} reduce to
\begin{equation*}
\bm{H}_{\alpha,\beta}^x=\left(0,\frac{g}{2}\left(\widehat w-B_{\alpha,\beta}\right)^2,0\right)^\top,\quad
\bm{H}_{\gamma,\delta}^y=\left(0,0,\frac{g}{2}\left(\widehat w-B_{\gamma,\delta}\right)^2\right)^\top,
\end{equation*}
and the flux terms on the RHS of \eqref{3.6} then become
\begin{equation}
\begin{aligned}
&-\frac{\bm{H}_{\jph,k}^x-\dfrac{\bm{H}_{\jmh,k-\frac{1}{4}}^x+\bm{H}_{\jmh,k+\frac{1}{4}}^x}{2}}{\dx_{j,k}}
-\frac{\bm{H}_{j,\kph}^y-\bm{H}_{j,\kmh}^y}{\dy_{j,k}}\\
&=-\frac{g}{2}\begin{pmatrix}
\dfrac{\left(\widehat w-B_{\jph,k}\right)^2}{\dx_{j,k}}-\dfrac{\left(\widehat w-B_{\jmh,k-\frac{1}{4}}\right)^2}{2\dx_{j,k}}-
\dfrac{\left(\widehat w-B_{\jmh,k+\frac{1}{4}}\right)^2}{2\dx_{j,k}}\\
\dfrac{\left(\widehat w-B_{j,\kph}\right)^2}{\dy_{j,k}}-\dfrac{\left(\widehat w-B_{j,\kmh}\right)^2}{\dy_{j,k}}
\end{pmatrix}.
\end{aligned}
\label{3.14}
\end{equation}

We now need to approximate the source term in \eqref{3.6} in such a way that $\,\xbar{\bm{S}}_{j,k}$ would cancel \eqref{3.14} at the
``lake-at-rest'' steady states. To this end, we first notice that (at least for smooth solutions)
$$
\begin{aligned}
&-g(w-B)B_x=g(w-B)(w-B)_x-g(w-B)w_x=\frac{g}{2}\left[(w-B)^2\right]_x-g(w-B)w_x,\\
&-g(w-B)B_y=g(w-B)(w-B)_y-g(w-B)w_y=\frac{g}{2}\left[(w-B)^2\right]_y-g(w-B)w_y,
\end{aligned}
$$
and rewrite the cell averages of the second and third components of the integral in \eqref{3.9} as
\begin{equation}
\frac{g}{2}\int\limits_{y_\kmh}^{y_\kph}\left[(w-B)^2\Big|_{x=x_\jph}-(w-B)^2\Big|_{x=x_\jmh}\right]\,{\rm d}y-
g\int\limits_{x_\jmh}^{x_\jph}\int\limits_{y_\kmh}^{y_\kph}(w-B)w_x\,{\rm d}y\,{\rm d}x
\label{3.15}
\end{equation}
and
\begin{equation}
\frac{g}{2}\int\limits_{x_\jmh}^{x_\jph}\left[(w-B)^2\Big|_{y=y_\kph}-(w-B)^2\Big|_{y=y_\kmh}\right]\,{\rm d}x-
g\int\limits_{x_\jmh}^{x_\jph}\int\limits_{y_\kmh}^{y_\kph}(w-B)w_y\,{\rm d}y\,{\rm d}x,
\label{3.16}
\end{equation}
respectively. We then approximate the integrals in \eqref{3.15} and \eqref{3.16} using the second-order midpoint rule (for the
configuration in Figure \ref{fig:3} (b), the integral along the left edge of $C_{j,k}$ is approximated using the composite midpoint rule
as $C_{j,k}$ has two neighboring cells on the left), which results in the following quadrature for the second and third components of the
source term:
\begin{equation}
\begin{aligned}
\xbar{S}_{j,k}^{\,(2)}\approx\frac{g}{2\dx_{j,k}}\Bigg[&\left(w_{\jph,k}^--B_{\jph,k}\right)^2-
\frac{\left(w_{\jmh,k-\frac{1}{4}}^+-B_{\jmh,k-\frac{1}{4}}\right)^2}{2}\\
&-\frac{\left(w_{\jmh,k+\frac{1}{4}}^+-B_{\jmh,k+\frac{1}{4}}\right)^2}{2}\Bigg]-g(w_x)_{j,k}\left(\,\xbar w_{j,k}-B_{j,k}\right),\\
\xbar{S}_{j,k}^{\,(3)}\approx\frac{g}{2\dy_{j,k}}\bigg[&\left(w_{j,\kph}^--B_{j,\kph}\right)^2-\left(w_{j,\kmh}^+-B_{j,\kmh}\right)^2
\bigg]-g(w_y)_{j,k}\left(\,\xbar w_{j,k}-B_{j,k}\right).
\end{aligned}
\label{3.17}
\end{equation}

We finally note that at the ``lake-at-rest'' data, $(w_x)_{j,k}=(w_y)_{j,k}\equiv0,\,\forall(j,k)$ and thus \eqref{3.14} and \eqref{3.17} 
imply that the RHS of \eqref{3.6} vanishes and the resulting scheme is well-balanced.

\subsection{Positivity preserving property and time discretization}\label{S:4.7}
One of the main advantages of the central-upwind scheme is its ability to preserve the positivity of $h$; see \cite{Kur_Acta,Kurganov2007}.
In this section, we extend the positivity proof from \cite{Kurganov2007} to the proposed quadtree scheme. To this end, we integrate
equation \eqref{3.6} in time using a forward Euler method. For the first component, this results in
\begin{equation}
\xbar w_{j,k}^{\,n+1}=\,\xbar w_{j,k}^{\,n}-\lambda_{j,k}^n\left(H_{\jph,k}^{x,(1)}-\frac{H_{\jmh,k-\frac{1}{4}}^{x,(1)}+
H_{\jmh,k+\frac{1}{4}}^{x,(1)}}{2}\right)-\mu_{j,k}^n\left(H_{j,\kph}^{y,(1)}-H_{j,\kmh}^{y,(1)}\right),
\label{3.18}
\end{equation}
where $\,\xbar w_{j,k}^{\,n}:=\,\xbar w_{j,k}(t^n)$ and $\,\xbar w_{j,k}^{\,n+1}:=\,\xbar w_{j,k}(t^{n+1})$ with $t^{n+1}=t^n+\dt^n$,
$\lambda_{j,k}^n:=\dt^n/\dx_{j,k}$, $\mu_{j,k}^n:=\dt^n/\dy_{j,k}$, and the numerical fluxes on the RHS are evaluated at time level $t=t^n$ 
using \eqref{3.13}:
\begin{equation}
\begin{aligned}
H_{\alpha,\beta}^{x,(1)}=&\frac{a_{\alpha,\beta}^+(hu)_{\alpha,\beta}^--a_{\alpha,\beta}^-(hu)_{\alpha,\beta}^+}
{a_{\alpha,\beta}^+-a_{\alpha,\beta}^-}+\frac{a_{\alpha,\beta}^+a_{\alpha,\beta}^-}{a_{\alpha,\beta}^+-a_{\alpha,\beta}^-}
\left[w_{\alpha,\beta}^+-w_{\alpha,\beta}^-\right],\\
H_{\gamma,\delta}^{y,(1)}=&\frac{b_{\gamma,\delta}^+(hv)_{\gamma,\delta}^--b_{\gamma,\delta}^-(hv)_{\gamma,\delta}^+}
{b_{\gamma,\delta}^+-b_{\gamma,\delta}^-}+\frac{b_{\gamma,\delta}^+b_{\gamma,\delta}^-}{b_{\gamma,\delta}^+-b_{\gamma,\delta}^-}
\left[w_{\gamma,\delta}^+-w_{\gamma,\delta}^-\right],
\end{aligned}
\label{3.25}
\end{equation}
where, as before, $(\alpha,\beta)\in\big\{(\jmh,k-\frac{1}{4}),(\jmh,k+\frac{1}{4}),(\jph,k)\big\}$ and
$(\gamma,\delta)\in\big\{(j,\kmh),(j,\kph)\big\}$ for the configuration considered in Figure \ref{fig:3} (b).

If $\,\xbar h_{j,k}^{\,n}\ge0$ for all $(j,k)$, then the point values of $h$ computed using piecewise linear/bilinear reconstructions of $w$
and $B$ presented in \S\ref{S:4.5} and \S\ref{S:4.4}, are nonnegative. Moreover, using \eqref{3.12a}--\eqref{3.15a} and the similar
relationships for the reconstructed point values of $w$, we have
\begin{equation}
\xbar h_{j,k}^{\,n}=\frac{1}{4}\left(\frac{h_{\jmh,k-\frac{1}{4}}^++h_{\jmh,k+\frac{1}{4}}^+}{2}+h_{\jph,k}^-+h_{j,\kmh}^++h_{j,\kph}^-\right)
\label{3.26}
\end{equation}
for the configuration considered in Figure \ref{fig:3} (b).

We now subtract $B_{j,k}$ from both sides of \eqref{3.18} and use \eqref{3.25} and \eqref{3.26} to rewrite \eqref{3.18} as follows:
$$
\begin{aligned}
\xbar h_{j,k}^{\,n+1}&=-\lambda_{j,k}^na_{\jph,k}^-\cdot\frac{a_{\jph,k}^+-u_{\jph,k}^+}{a_{\jph,k}^+-a_{\jph,k}^-}\cdot h_{\jph,k}^++
\left[\frac{1}{4}-\lambda_{j,k}^na_{\jph,k}^+\cdot\frac{u_{\jph,k}^--a_{\jph,k}^-}{a_{\jph,k}^+-a_{\jph,k}^-}\right]h_{\jph,k}^-\\
&+\frac{\lambda_{j,k}^na_{\jmh,k-\frac{1}{4}}^+}{2}\cdot\frac{u_{\jmh,k-\frac{1}{4}}^--a_{\jmh,k-\frac{1}{4}}^-}
{a_{\jmh,k-\frac{1}{4}}^+-a_{\jmh,k-\frac{1}{4}}^-}\cdot h_{\jmh,k-\frac{1}{4}}^+\\
&+\hf\left[\frac{1}{4}-\lambda_{j,k}^na_{\jmh,k-\frac{1}{4}}^-\cdot\frac{a_{\jmh,k-\frac{1}{4}}^+-u_{\jmh,k-\frac{1}{4}}^+}
{a_{\jmh,k-\frac{1}{4}}^+-a_{\jmh,k-\frac{1}{4}}^-}\right]h_{\jmh,k-\frac{1}{4}}^+\\
&+\frac{\lambda_{j,k}^na_{\jmh,k+\frac{1}{4}}^+}{2}\cdot\frac{u_{\jmh,k+\frac{1}{4}}^--a_{\jmh,k+\frac{1}{4}}^-}
{a_{\jmh,k+\frac{1}{4}}^+-a_{\jmh,k+\frac{1}{4}}^-}\cdot h_{\jmh,k+\frac{1}{4}}^+\\
&+\hf\left[\frac{1}{4}-\lambda_{j,k}^na_{\jmh,k+\frac{1}{4}}^-\cdot\frac{a_{\jmh,k+\frac{1}{4}}^+-u_{\jmh,k+\frac{1}{4}}^+}
{a_{\jmh,k+\frac{1}{4}}^+-a_{\jmh,k+\frac{1}{4}}^-}\right]h_{\jmh,k+\frac{1}{4}}^+\\
&-\mu_{j,k}^nb_{j,\kph}^-\cdot\frac{b_{j,\kph}^+-v_{j,\kph}^+}{b_{j,\kph}^+-b_{j,\kph}^-}\cdot h_{j,\kph}^++
\left[\frac{1}{4}-\mu_{j,k}^nb_{j,\kph}^+\cdot\frac{v_{j,\kph}^--b_{j,\kph}^-}{b_{j,\kph}^+-b_{j,\kph}^-}\right]h_{j,\kph}^-\\
&+\mu_{j,k}^nb_{j,\kmh}^+\cdot\frac{v_{j,\kmh}^--b_{j,\kmh}^-}{b_{j,\kmh}^+-b_{j,\kmh}^-}\cdot h_{j,\kmh}^-+
\left[\frac{1}{4}+\mu_{j,k}^nb_{j,\kmh}^-\cdot\frac{b_{j,\kmh}^+-v_{j,\kmh}^+}{b_{j,\kmh}^+-b_{j,\kmh}^-}\right]h_{j,\kmh}^+.
\end{aligned}
\label{3.27} 
$$
This shows that the cell averages of $h$ at the new time level can be written as a linear combination of the reconstructed nonnegative point
values of $h$. Therefore, $\,\xbar h_{j,k}^{\,n+1}\ge0$ provided all of the coefficients in this linear combination are nonnegative, which is,
using the definition of the local speeds of propagation in \eqref{3.12}, true provided the following CFL-type condition are satisfied:
$$
\dt\le\frac{1}{4}\min\left[\min_{j,k}\left\{
\frac{\dx_{j,k}}{\max\limits_{(\alpha,\beta)}\left[\max\left\{a_{\alpha,\beta}^+,-a_{\alpha,\beta}^-\right\}\right]}\right\},\,
\min_{j,k}\left\{
\frac{\dy_{j,k}}{\max\limits_{(\gamma,\delta)}\left[\max\left\{b_{\gamma,\delta}^+,-b_{\gamma,\delta}^-\right\}\right]}\right\}\right].
$$
where, as before, $(\alpha,\beta)\in\big\{(\jmh,k-\frac{1}{4}),(\jmh,k+\frac{1}{4}),(\jph,k)\big\}$ and
$(\gamma,\delta)\in\big\{(j,\kmh),(j,\kph)\big\}$ for the configuration considered in Figure \ref{fig:3} (b).

It should be observed that the above positivity preserving proof is valid not only for the forward Euler time discretization, but for any
strong stability preserving (SSP) ODE solver (see, e.g., \cite{GKS,GST}) as well. In all of our numerical experiments, we have used the
three-stage third-order SSP Runge-Kutta solver.

\subsection{Quadtree grid adaptivity}\label{S:4.10}
After evolving the solution to the new time level $t=t^{n+1}$ the quadtree grid should be adapted (locally either refined or coarsened) to
the new solution structure. To this end, we first compute the slopes $\{(w_x)_{j,k}^{n+1}\}$ and $\{(w_y)_{j,k}^{n+1}\}$ on the old grid
(which we denote by $\{C_{j,k}^{\rm old}\}$) according to \S\ref{S:4.4} and then select the centers of those cells $C_{j,k}^{\rm old}$, at
which either
\begin{equation}
(w_x)_{j,k}^{n+1}\ge C_{\rm seed}\quad\mbox{or}\quad(w_y)_{j,k}^{n+1}\ge C_{\rm seed},
\label{3.28}
\end{equation}
to be the seeding points needed to generate the new grid, which we denote by $\{C_{j,k}^{\rm new}\}$. In \eqref{3.28}, $C_{\rm seed}$ is a
constant that depends on the problem at hand, that is, on such factors as the Froude number, bottom topography function and/or boundary
conditions. 

When the mesh is locally refined or coarsened, the solution realized at the end of the evolution step in terms of the computed cell
averages  $\big\{\big(\,\xbar{\bm{U}}_{j,k}^{\,n+1}\big)_{\rm old}\big\}$ over the grid $\{C_{j,k}^{\rm old}\}$, should be projected onto
the new grid $\{C_{j,k}^{\rm new}\}$ in a conservative manner according to the following three possible cases.

{\em Case 1:} If $C_{j,k}^{\rm new}=C_{j',k'}^{\rm old}$ for some $(j',k')$, that is, if the cell $C_{j',k'}^{\rm old}$ does not need to be
refined/coarsened, then
\begin{equation*}
\big(\,\xbar{\bm{U}}_{j,k}^{\,n+1}\big)_{\rm new}=\big(\,\xbar{\bm{U}}_{j',k'}^{\,n+1}\big)_{\rm old}.
\end{equation*}

{\em Case 2:} If $C_{j,k}^{\rm new}\in{\cal C}^{\ell+p}$ is a ``child'' cell of $C_{j',k'}^{\rm old}\in{\cal C}^{\ell}$ for some $j'$, $k'$
and $p>0$ (that is, if the cell $C_{j',k'}^{\rm old}$ was refined and $C_{j,k}^{\rm new}\subset C_{j',k'}^{\rm old}$), then
\begin{equation*}
\big(\,\xbar{\bm{U}}^{\,n+1}_{j,k}\big)_{\rm new}=\big(\,\xbar{\bm{U}}^{\,n+1}_{j',k'}\big)_{\rm old}+
\big((\bm{U}_x)^{n+1}_{j',k'}\big)_{\rm old}\left[x_j^{\rm new}-x_{j'}^{\rm old}\right]+
\big((\bm{U}_y)^{n+1}_{j',k'}\big)_{\rm old}\left[y_k^{\rm new}-y_{k'}^{\rm old}\right].
\end{equation*}

{\em Case 3:} If $C_{j,k}^{\rm new}\in{\cal C}^{\ell-p}$ is a ``parent'' cell of $C_{j',k'}^{\rm old}\in{\cal C}^{\ell}$ for some $j'$, $k'$
and $p>0$ (that is, if the cell $C_{j',k'}^{\rm old}$ was coarsened and $C_{j,k}^{\rm new}\supset C_{j',k'}^{\rm old}$), then
\begin{equation*}
\big(\,\xbar{\bm{U}}^{\,n+1}_{j,k}\big)_{\rm new}=\frac{1}{4^p}\underset{j'',k'':\,C_{j'',k''}^{\rm old}\subset C_{j,k}^{\rm new}}{\sum\sum}
\big(\,\xbar{\bm{U}}^{\,n+1}_{j'',k''}\big)_{\rm old}.
\end{equation*}

\section{Numerical experiments}\label{S:5}
In this section, we present six numerical examples in which the central-upwind quadtree scheme is tested. In all of the examples (except for
Example 5), we take $g=1$ and obtain the point values of $B$ at the vertices of $C_{j,k}$ using the bottom topography function with
$\ell:=5$ (\S\ref{S:4.5}).

\subsection*{Example 1 --- Accuracy test}
\begin{figure}[ht!]
\centerline{\includegraphics[height=3.6cm]{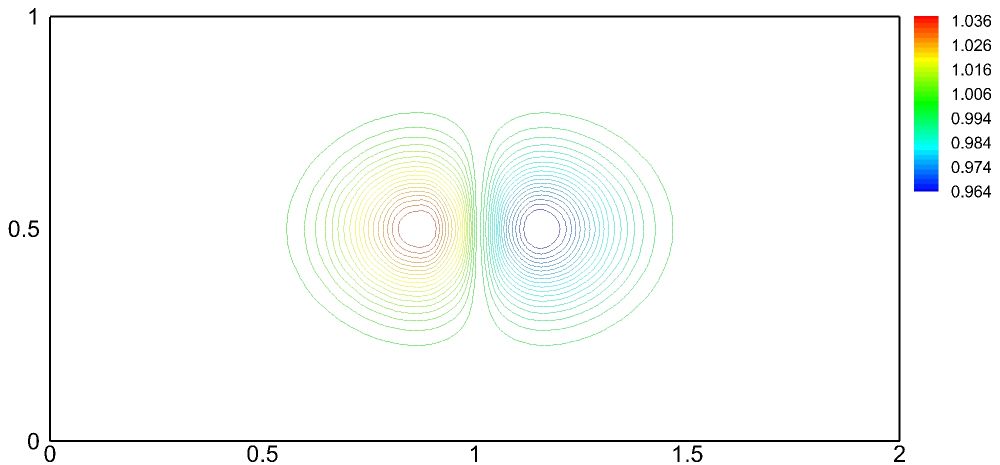}\hspace*{1cm}\includegraphics[height=3.6cm]{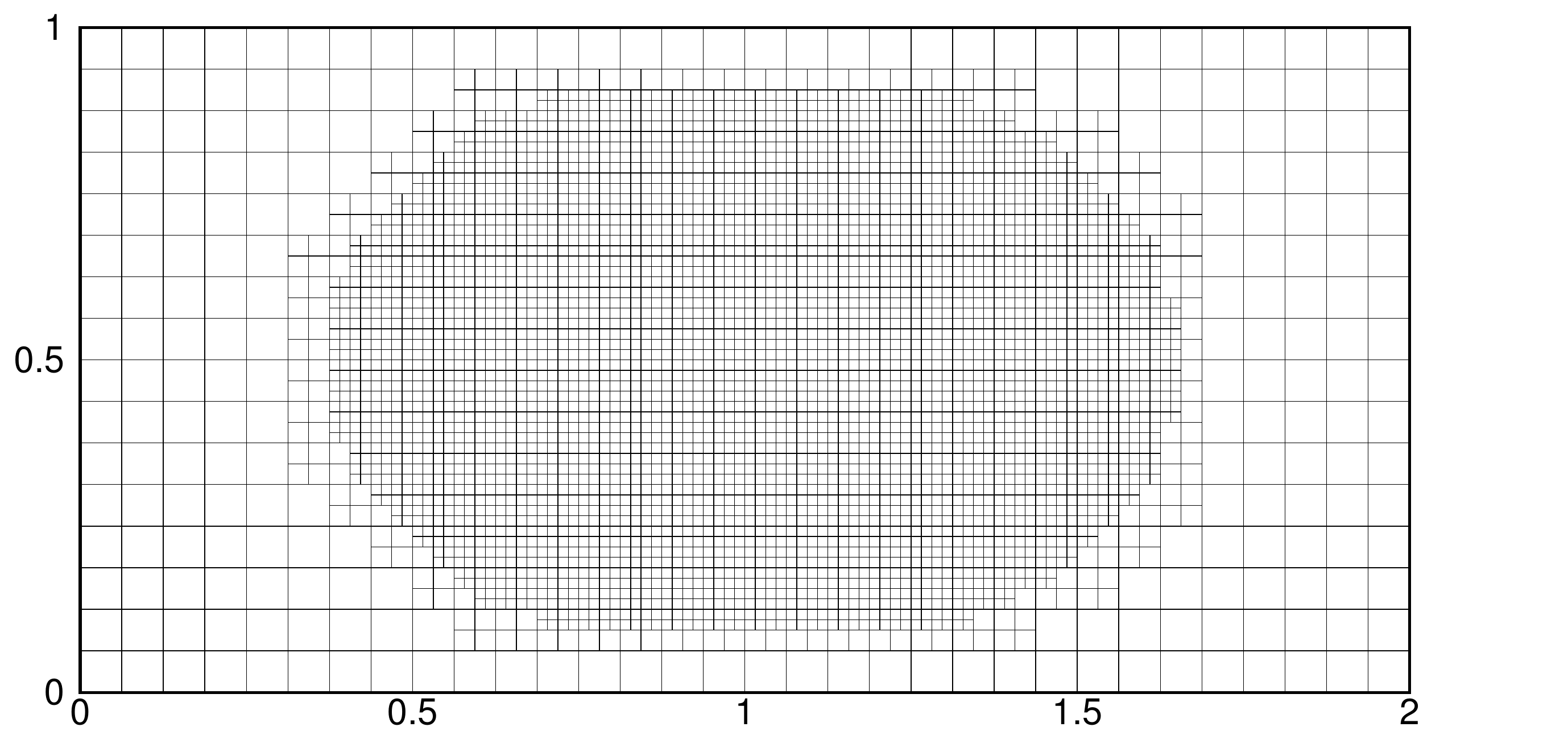}}
\caption{\sf Example 1: Computed water surface $w(x,y,0.07)$ for $m=7$ (left) and its corresponding quadtree grid (right).\label{fig:6.1}}
\end{figure}

In this benchmark, the accuracy of the proposed scheme is tested. We set the computational domain $[0,2]\times[0,1]$ with a zero-order 
extrapolation at all of the boundaries. The following initial data and the bottom topography function are imposed:
$$
w(x,y,0)\equiv1,\qquad u(x,y,0)\equiv0.3,\qquad v(x,y,0)\equiv0,\qquad B(x,y)=0.5e^{-25(x-1)^2-50(y-0.5)^2}.
$$
We generate a structured Cartesian grid $[512\times256]$ for the reference solution. The solution converges to a steady state solution by
$t=0.07$. The $L^1$- and $L^\infty$-errors for $m=5, 6, 7$, and $8$ for $C_{\rm seed}=0.0005$ are presented in Table \ref{Tbl:Tbl1}. The
obtained errors are similar to the ones reported in \cite{Bryson2011,Kurganov2007,SMSK}. The steady state solution computed with $m=7$ and
the corresponding quadtree grid are shown in Figure \ref{fig:6.1}.
\begin{table}[ht!]
\caption{\sf Example 1: $L^1$- and $L^\infty$- errors and numerical orders of accuracy.\label{Tbl:Tbl1}}
\centering
\begin{tabular}{*5c}
\toprule
Quadtree level& $L^1$-error & Order & $L^\infty$-error & Order\\
\midrule
$m=5$&$8.97e-04$&$-$&$5.14e-03$&$-$\\
$m=6$&$4.35e-04$&1.05&$3.22e-03$&0.67\\
$m=7$&$2.80e-04$&1.68&$2.90e-03$&0.83\\
$m=8$&$2.32e-04$&1.95&$2.18e-03$&1.24\\
\bottomrule
\end{tabular}
\end{table}

\subsection*{Example 2 --- Circular dam break}
\begin{figure}[ht!]
\centerline{\includegraphics[height=6cm]{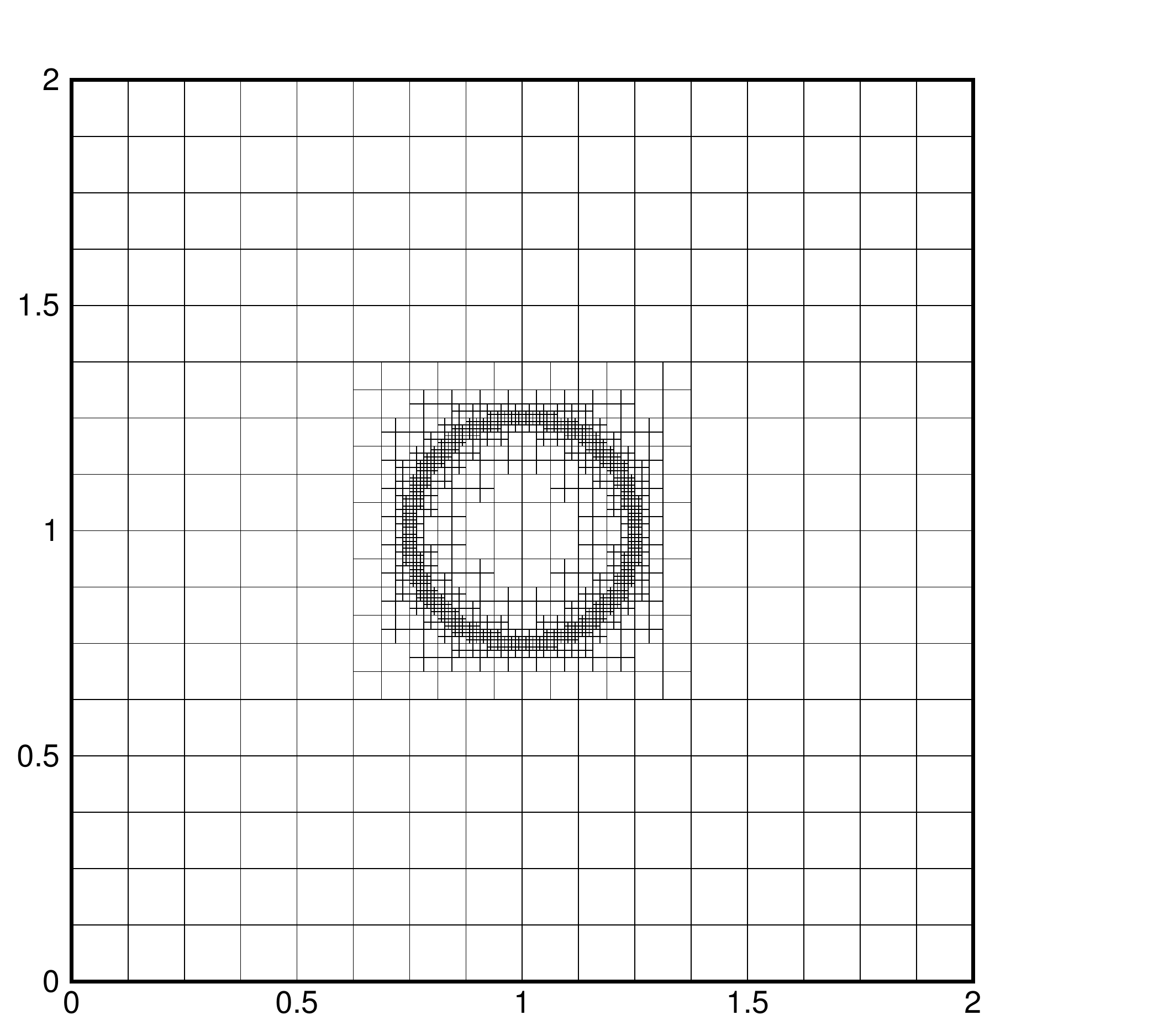}\hspace*{0.1cm}\includegraphics[height=6cm]{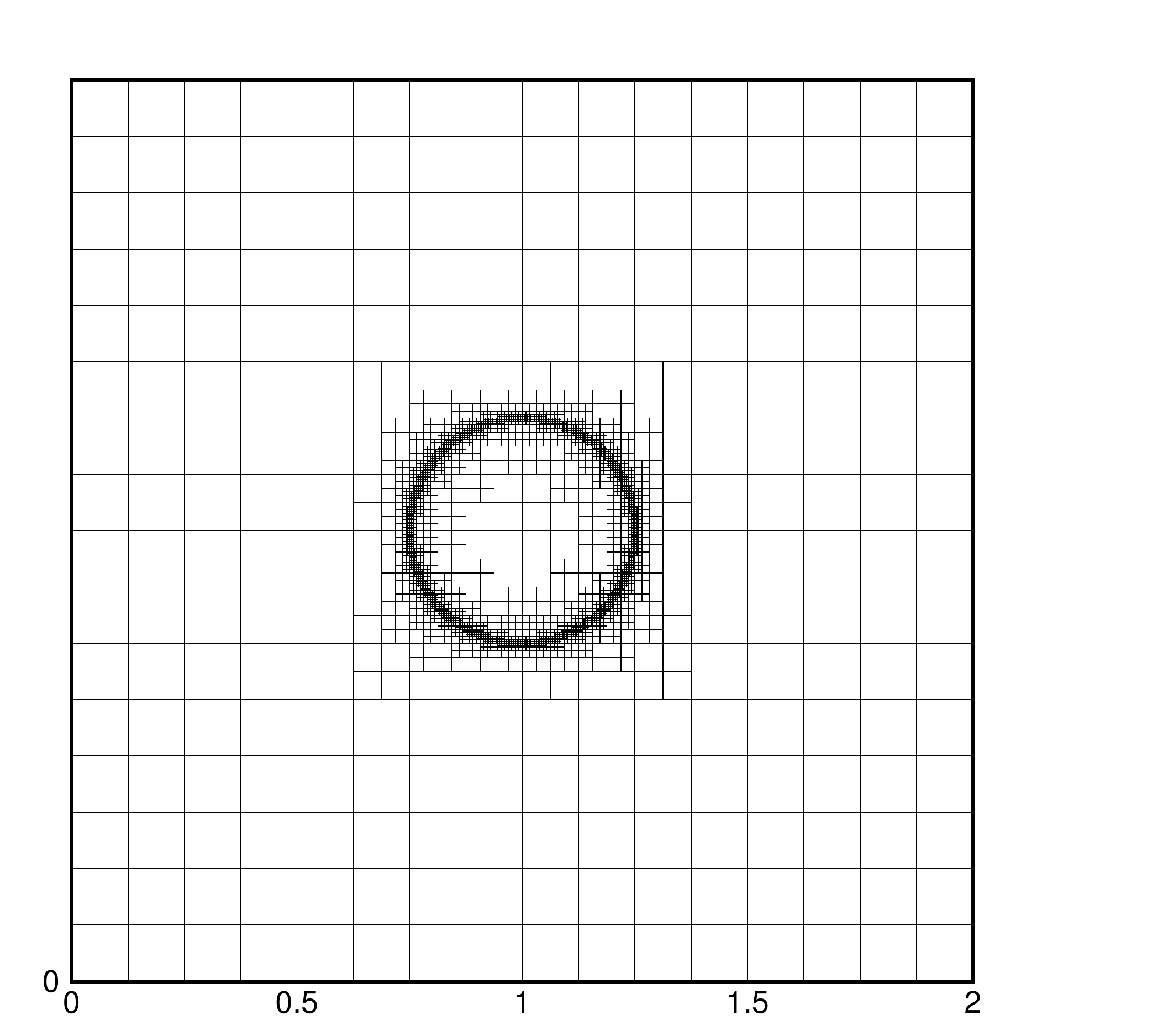}}
\caption{\sf Example 2: Initial quadtree grids for $m=8$ (left) and $m=9$ (right).\label{fig:4}}
\end{figure}

In this example, we demonstrate the ability of the proposed central-upwind quadtree scheme to preserve the positivity of the water surface
and to maintain symmetry. A circular water column, where $w=1$, collapses on a horizontal plane (similar examples were considered in
\cite{Audusse2005,Lin2003,Rogers2001,Rodriguez-Paz2005}), namely,
$$
w(x,y,0)=\left\{\begin{array}{lc}1,&~(x-1)^2+(y-1)^2<0.25,\\10^{-16},&\mbox{otherwise},\end{array}\right.\quad u(x,y,0)=v(x,y,0)\equiv0.
$$
We take the computational domain $[0,2]\times[0,2]$ and and impose zero-order extrapolated boundary conditions at its boundary. In this
example, we take $m=8$ and $m=9$ refinement levels of the quadtree grid and set $C_{\rm seed}=0.1$ in \eqref{3.28}. The initial quadtree
grids are shown in Figure \ref{fig:4}.

\begin{figure}[ht!]
\centerline{\includegraphics[height=6cm]{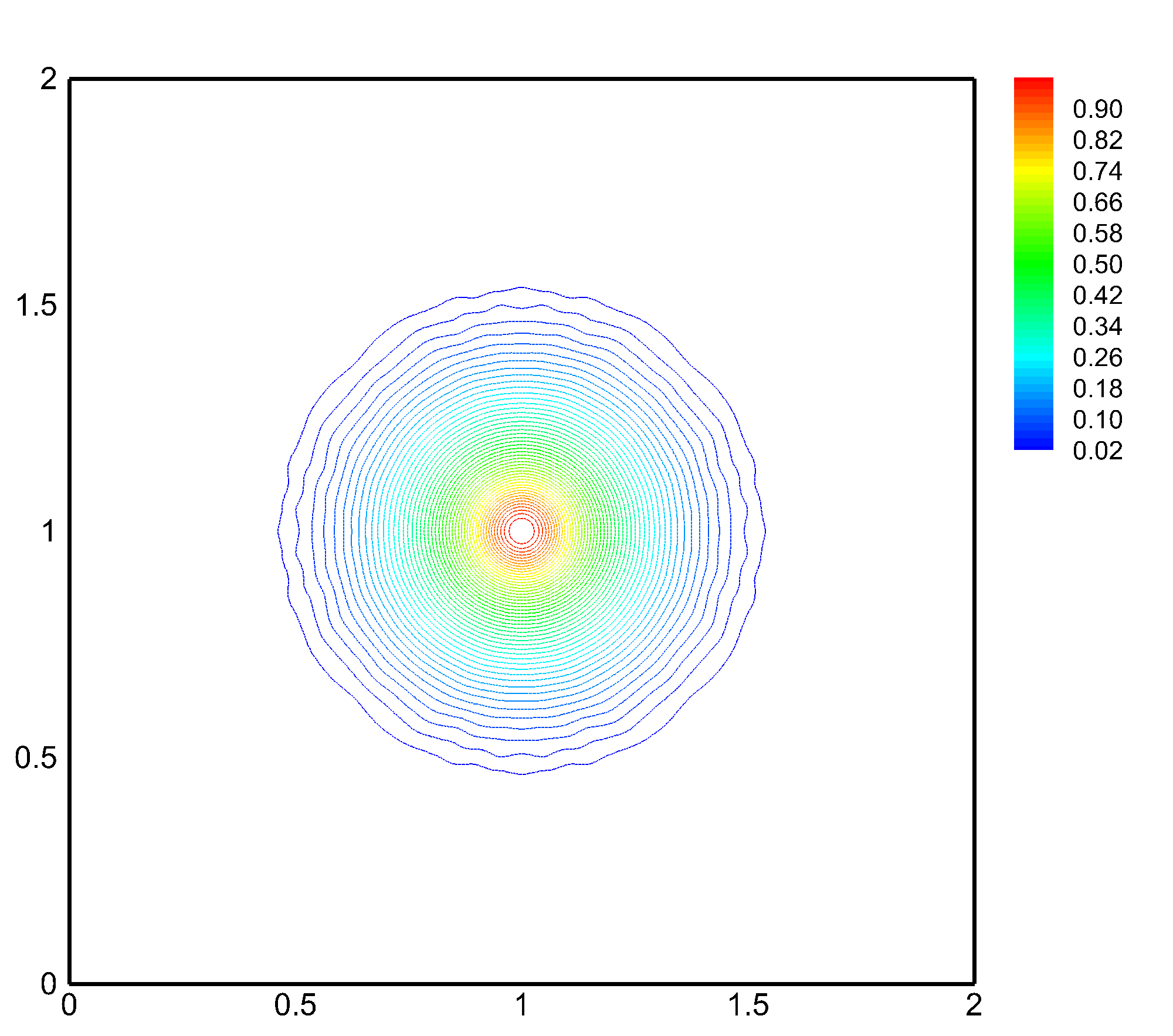}\hspace*{0.1cm}\includegraphics[height=6cm]{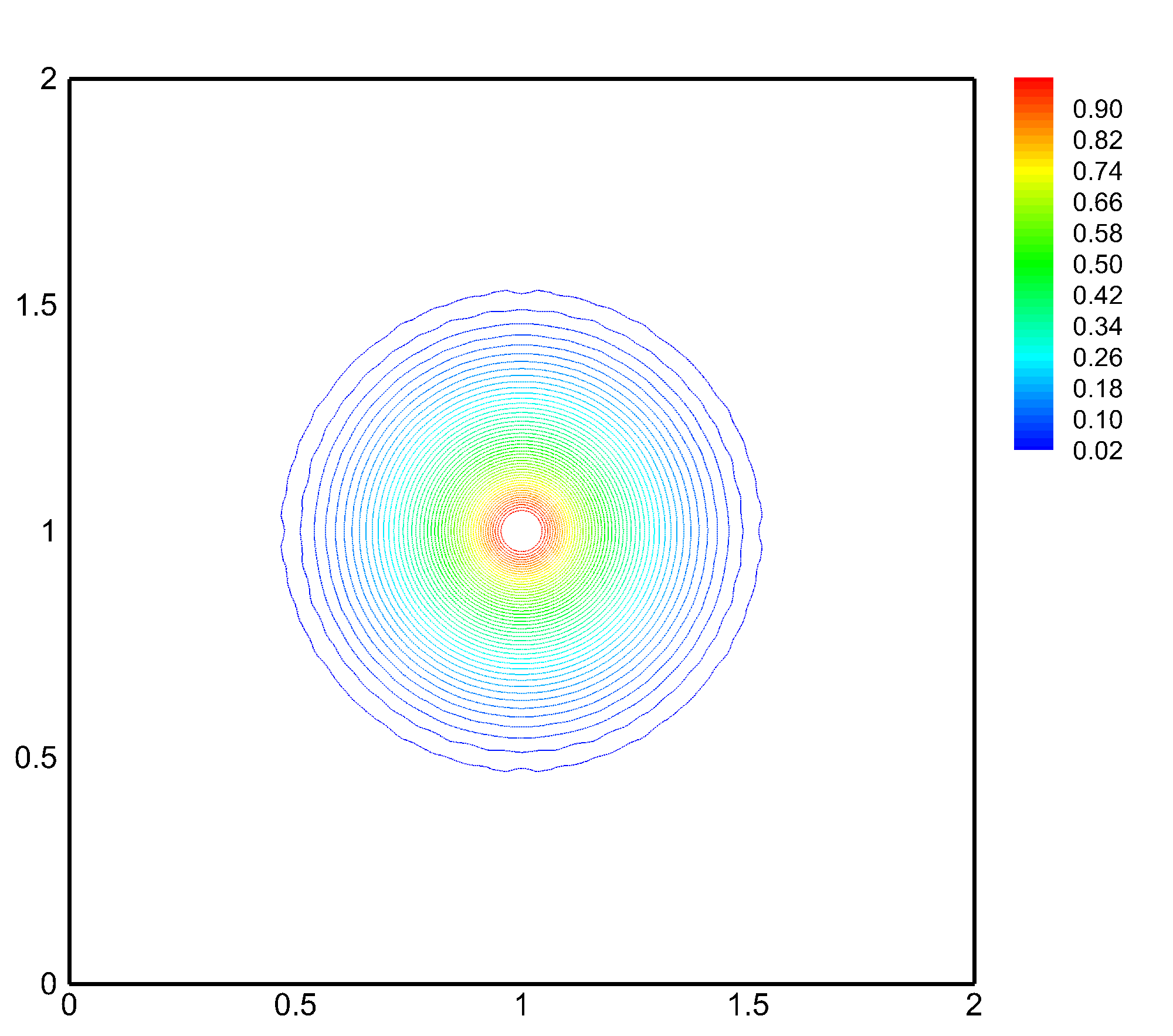}}
\caption{\sf Example 2: Computed water surface $w(x,y,0.2)$ for $m=8$ (left) and $m=9$ (right).\label{fig:4a}}
\end{figure}

We compute the solution until the final time $t=0.2$ and plot the obtained water surface contours in Figure \ref{fig:4a}. As one can see,
the central-upwind quadtree scheme maintains symmetry and preserves positivity. By changing the refinement level from $m=8$ to $m=9$, the
computational cost increases (for $m=8$, the quadtree grid starts with 1852 cells and ends with 12556 cells, whereas for $m=9$, the grid
starts with 3616 cells and ends with 56272 cells), but the results obtained with $m=9$ are clearly sharper and more accurate.

\subsection*{Example 3 --- Small perturbations of a stationary steady-state solution}
\begin{figure}[ht!]
\centerline{\includegraphics[height=5cm]{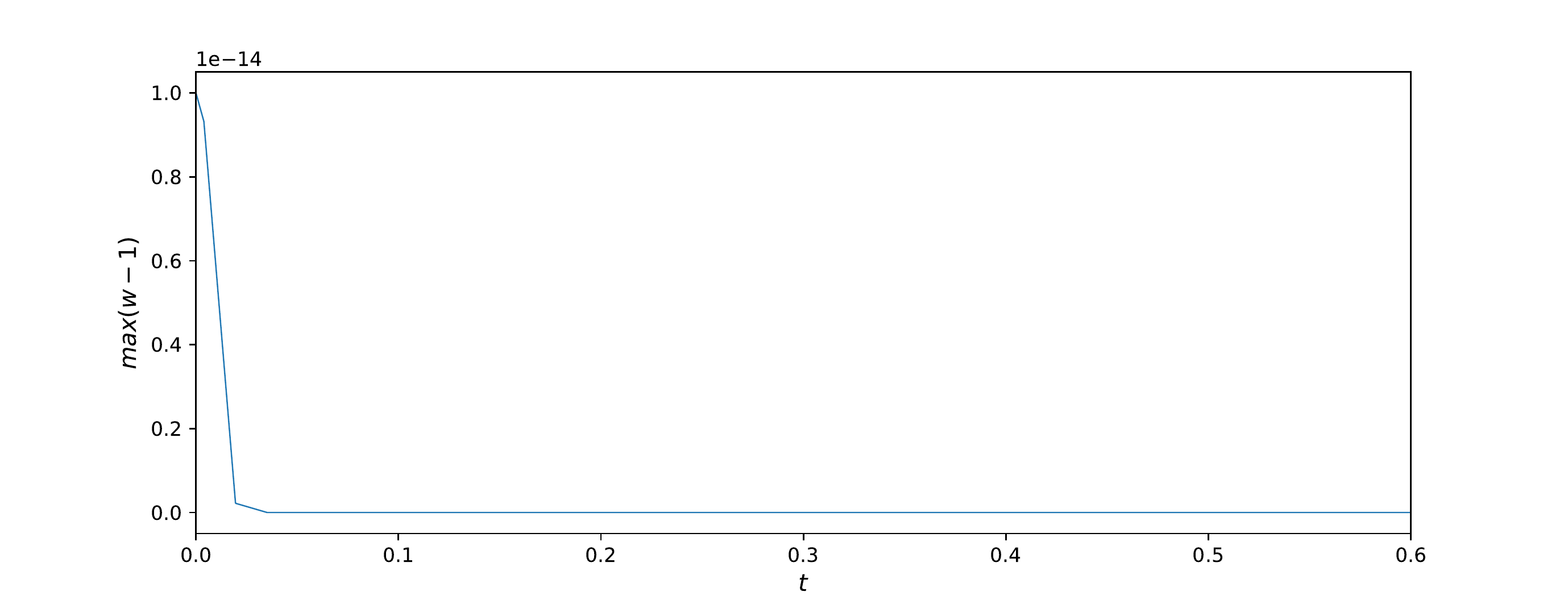}}
\caption{\sf Example 2: Computed $\displaystyle{\max_{x,y}w(x,y,t)}$ as a function of $t$ for $m=5$.\label{fig:4.52}}
\end{figure}

This numerical example is based on the benchmark, which was proposed in \cite{LeVeque1998} to test the ability of studied schemes to
accurately capture small perturbations of a steady state solution (similar examples were considered in, e.g.,
\cite{Bryson2011,Bryson2005,Kurganov2002,Liu2015,SMSK}). The computational domain is $[0,2]\times[0,1]$, the initial conditions are
$$
w(x,y,0)=\left\{\begin{array}{lc}1+\varepsilon,&~0.05<x<0.15,\\1,&~\mbox{otherwise},\end{array}\right.\qquad u(x,y,0)=v(x,y,0)\equiv0,
$$
and the bottom topography is given by
$$
B(x,y)=0.8e^{-5(x-0.9)^2-50(y-0.5)^2}.
$$
A solid wall boundary condition is used at the top and bottom boundaries and zero-order extrapolation is implemented at the left and right
ones. We first consider a very small value $\varepsilon=10^{-14}$ to verify the well-balanced property of the proposed quadtree scheme. The
solution is solved with a coarse quadtree grid for $m=5$. In Figure \ref{fig:4.52}, we plot $\displaystyle{\max_{x,y}(w-1)}$ as a function
of time until $t=0.6$. As one can see, the proposed scheme is stable and the fluxes and source terms balance each other. 

\begin{figure}[ht!]
\centerline{\includegraphics[height=3.6cm]{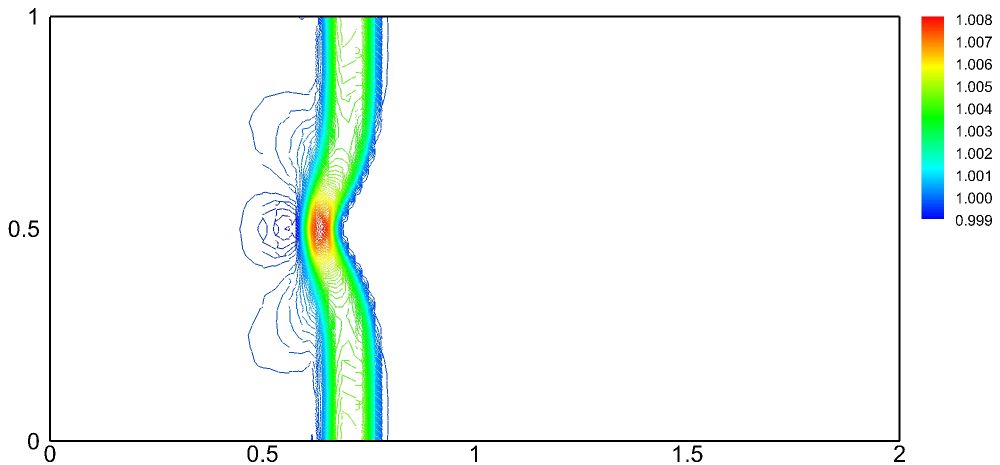}\hspace*{1cm}\includegraphics[height=3.6cm]{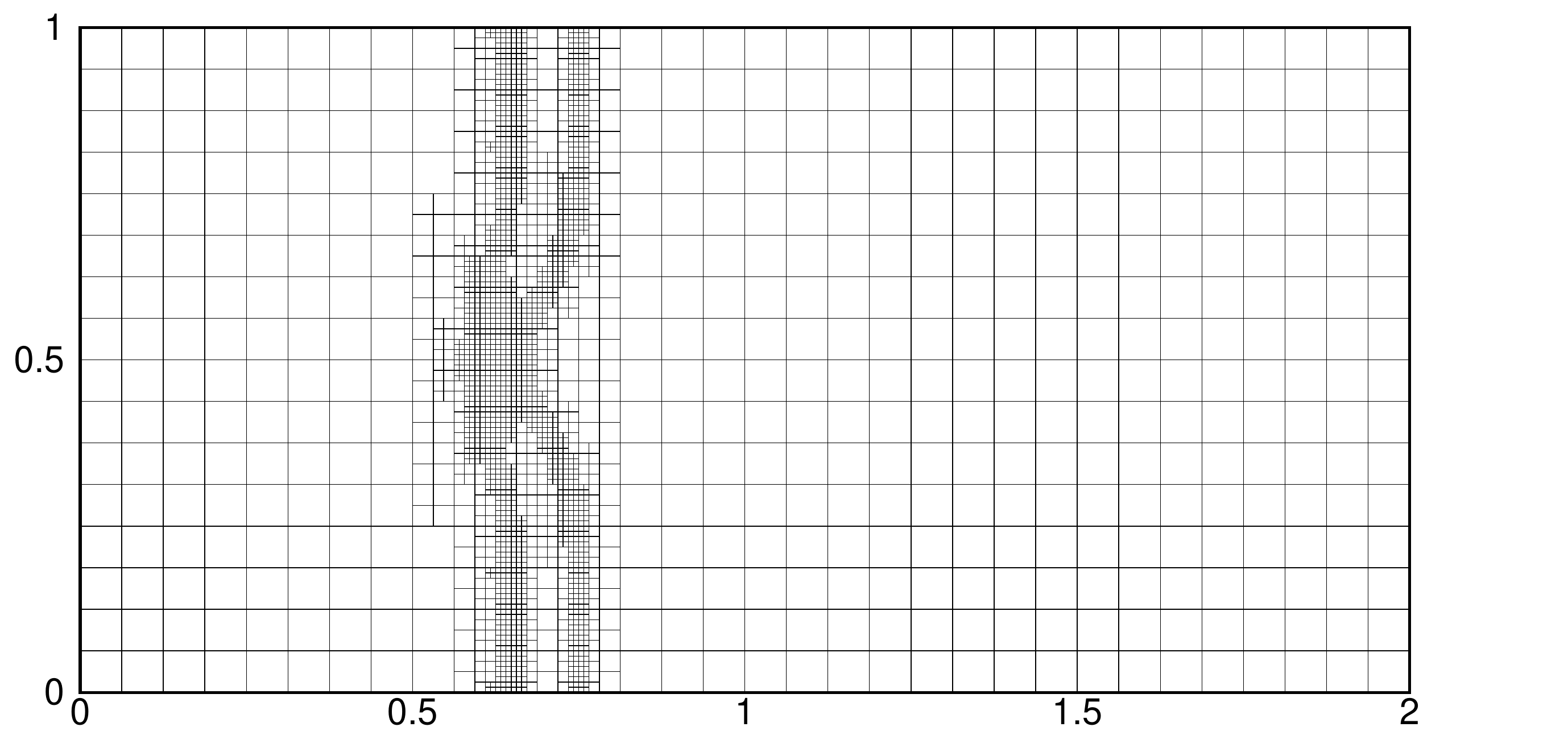}}
\vspace*{0.25cm}
\centerline{\includegraphics[height=3.6cm]{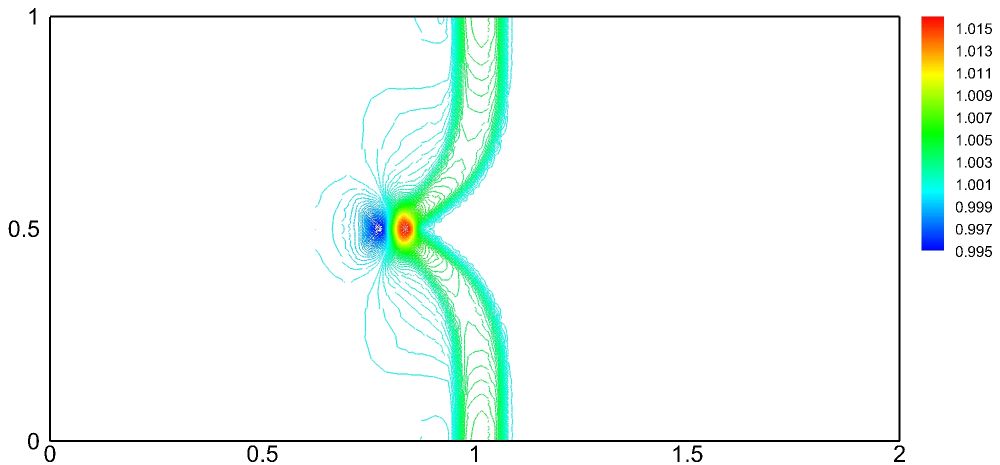}\hspace*{1cm}\includegraphics[height=3.6cm]{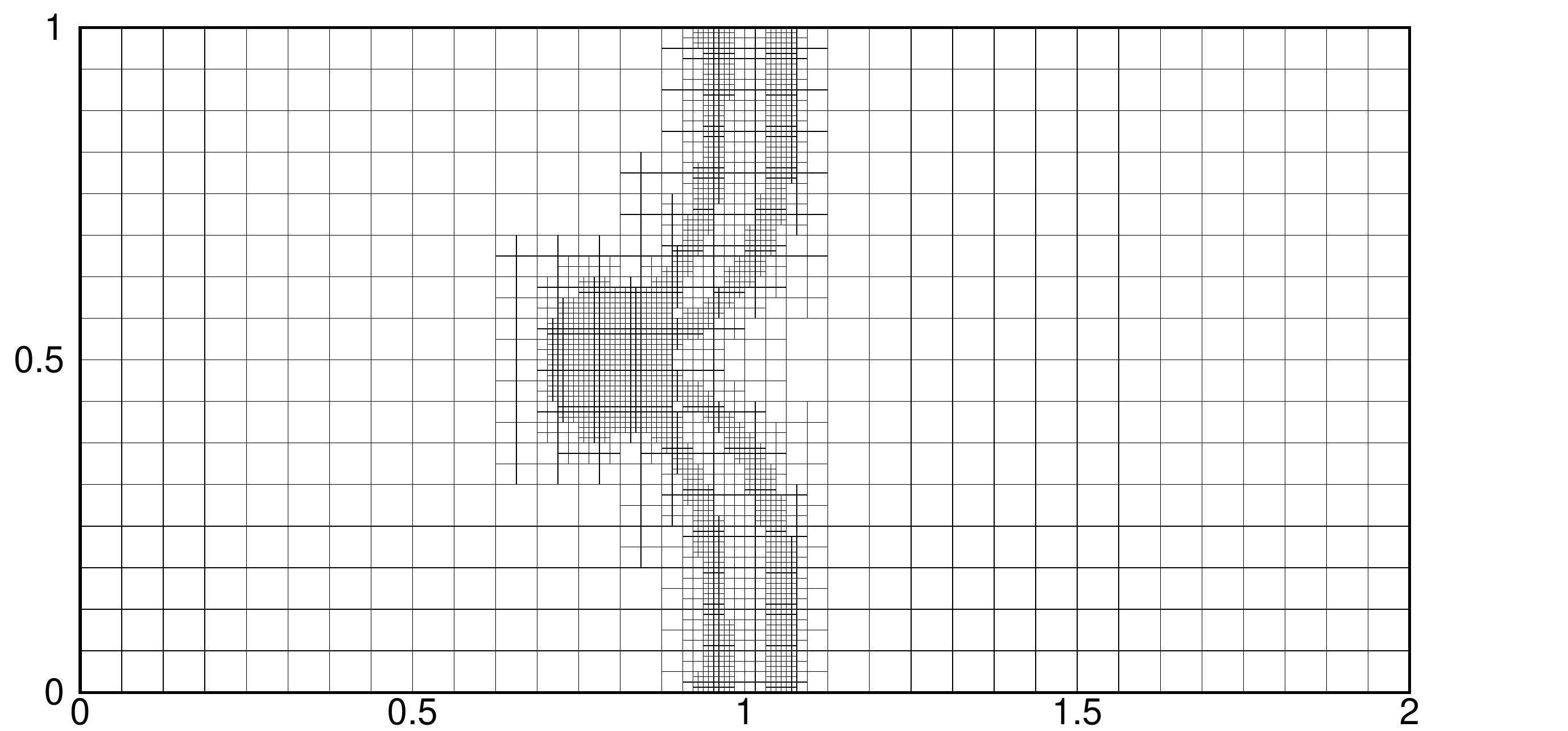}}
\vspace*{0.25cm}
\centerline{\includegraphics[height=3.6cm]{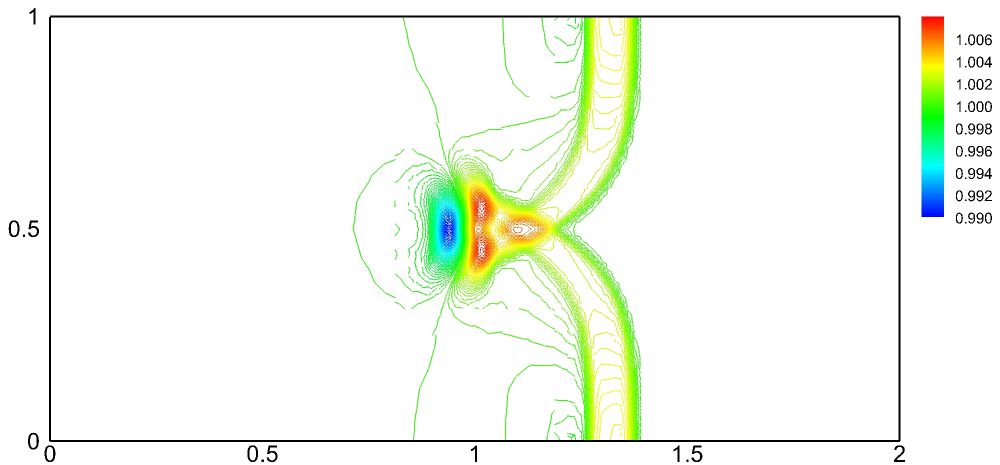}\hspace*{1cm}\includegraphics[height=3.6cm]{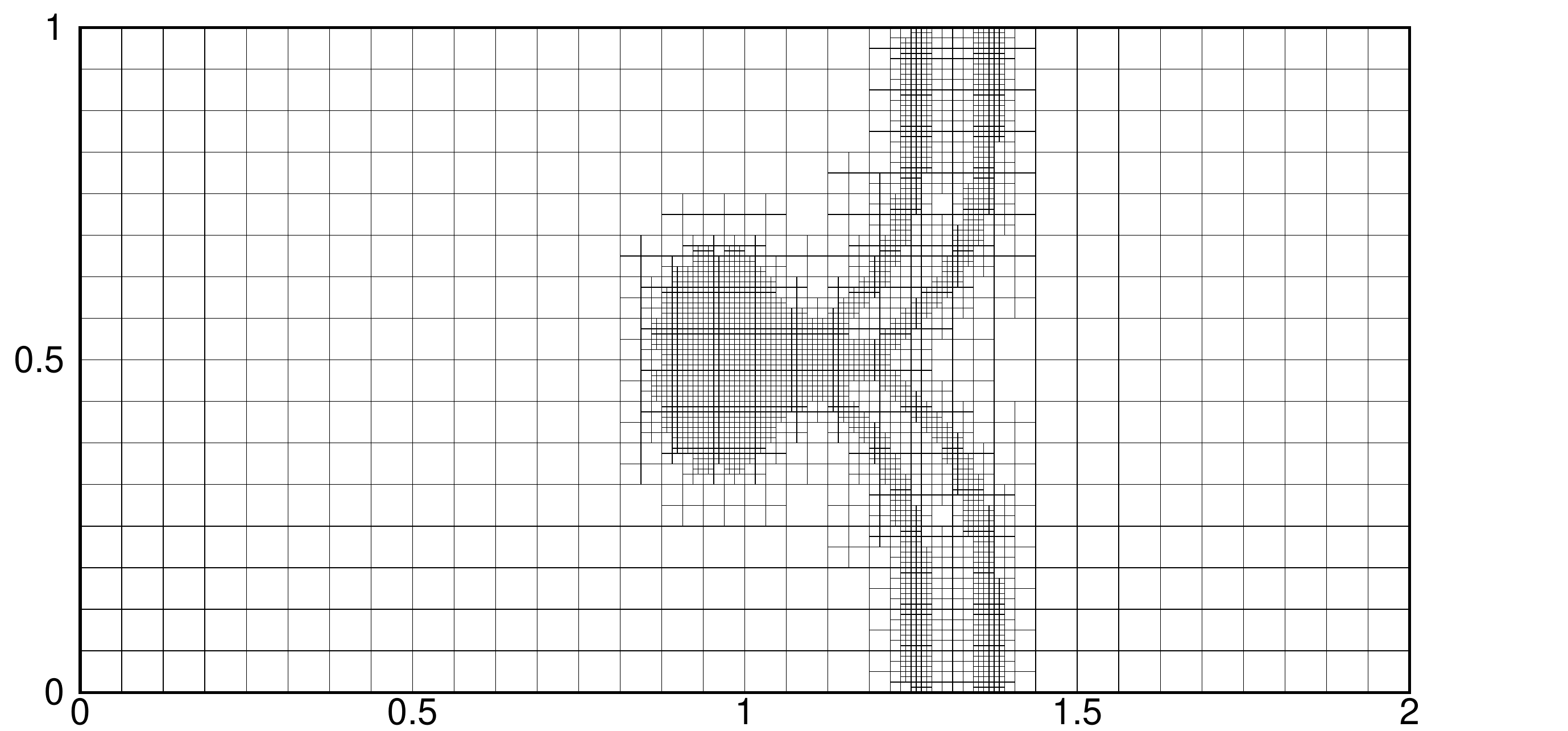}}
\vspace*{0.25cm}
\centerline{\includegraphics[height=3.6cm]{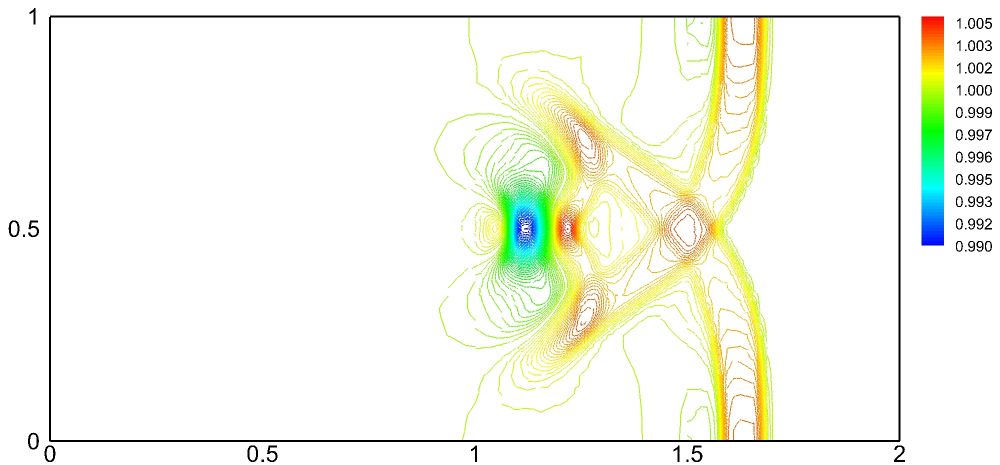}\hspace*{1cm}\includegraphics[height=3.6cm]{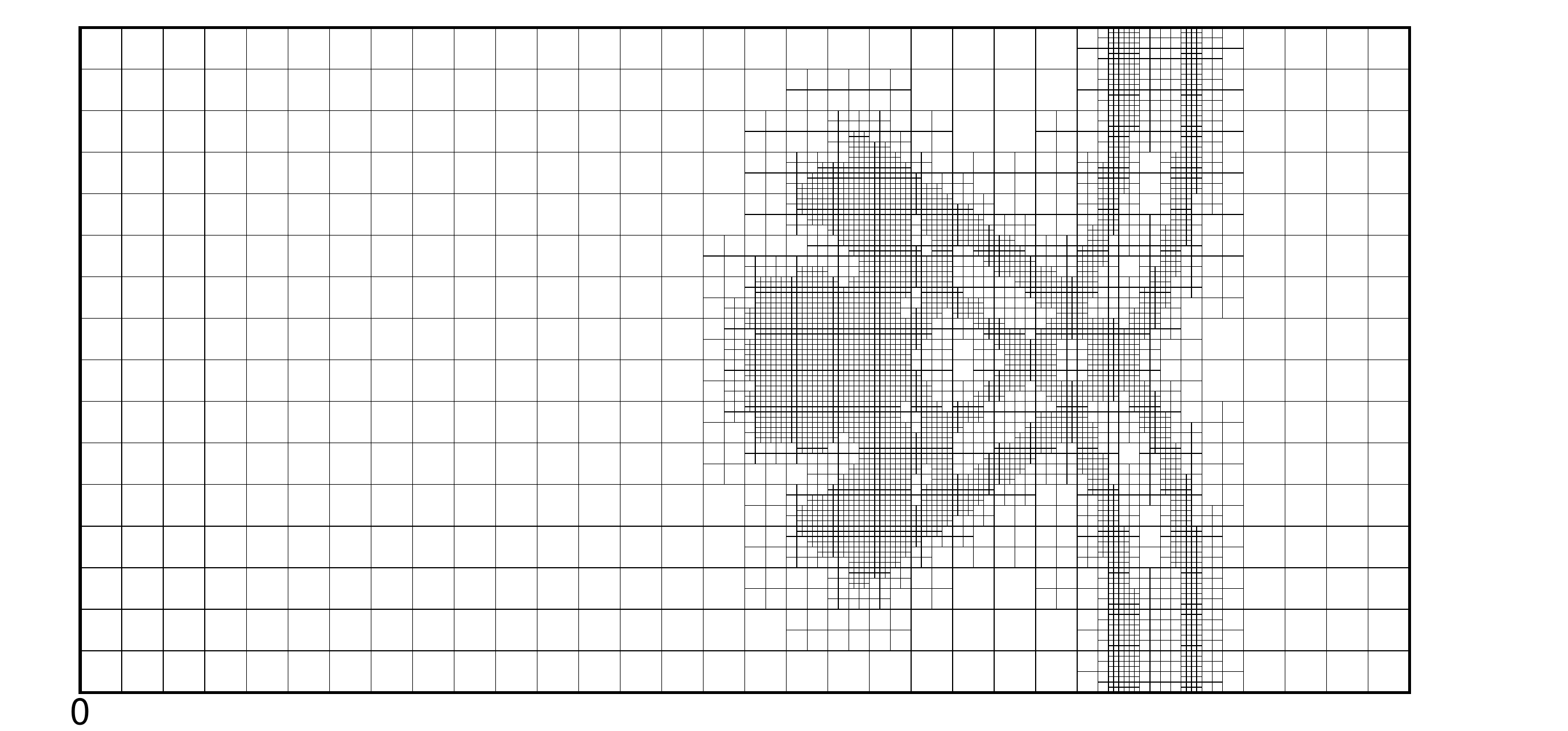}}
\vspace*{0.25cm}
\centerline{\includegraphics[height=3.6cm]{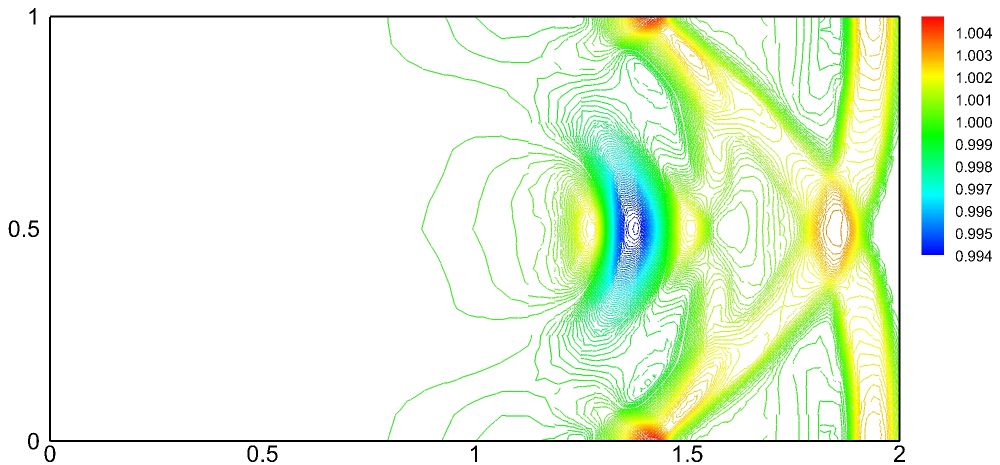}\hspace*{1cm}\includegraphics[height=3.6cm]{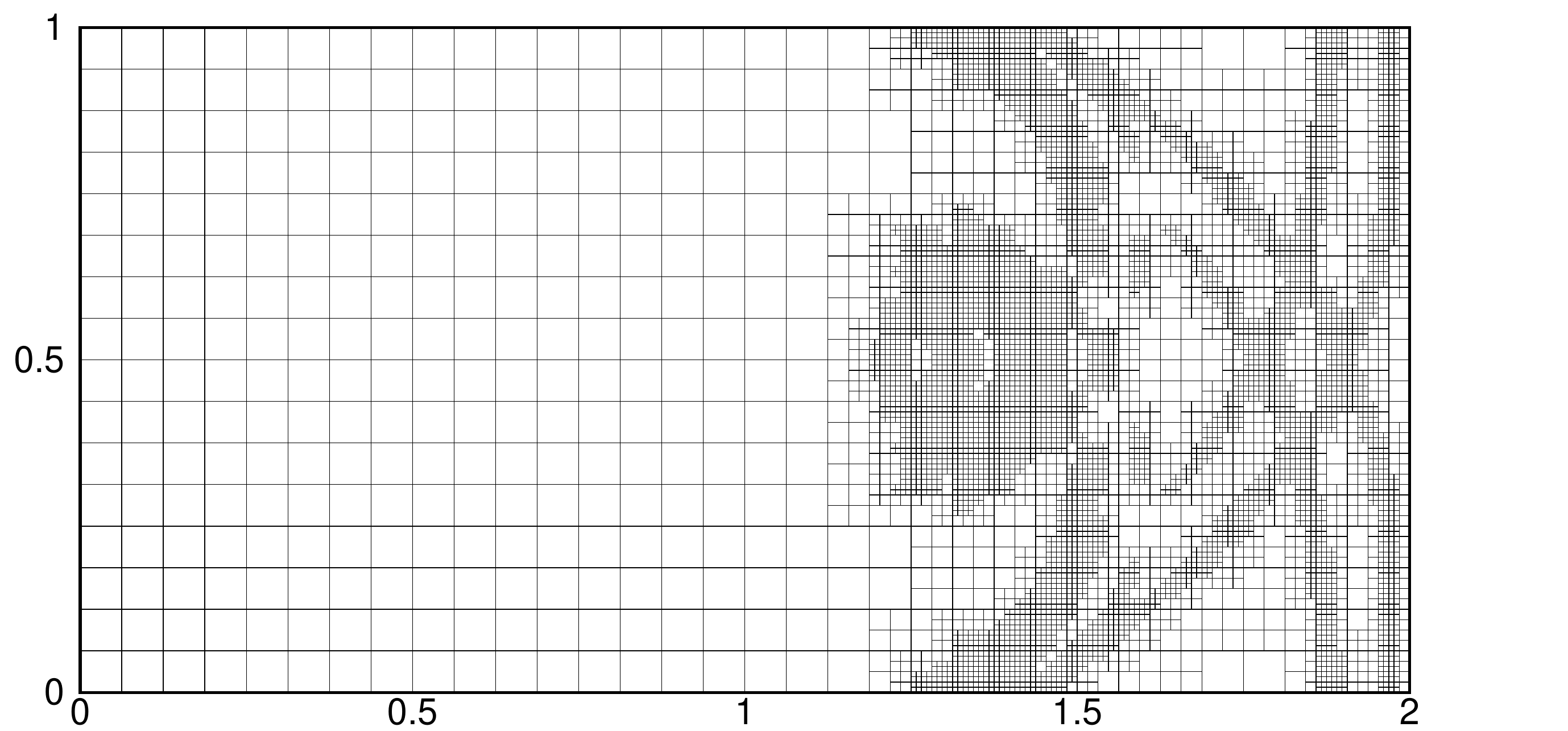}}
\caption{\sf Example 3: Computed water surface $w(x,y,t)$ (left column) and the corresponding quadtree grids (right column) for $t=0.6$,
$0.9$, $1.2$, $1.5$ and $1.8$ (from top to down) obtained using the well-balanced central-upwind quadtree scheme.\label{fig:5}}
\end{figure}

We then take $\varepsilon=0.01$ and $m=8$ refinement levels of the quadtree grid and set small $C_{\rm seed}=0.02$ in \eqref{3.28} in order
to accurately resolve small features of the computed solution. We compute the solution until the final time $t=1.8$ and plot the snapshots
of $w$ at times $t=0.6$, $0.9$, $1.2$, $1.5$ and $1.8$ in Figure \ref{fig:5} (left). The quadtree grid starts with 1970 cells and reaches a
maximum number of 7268 cells during the time evolution. Figure \ref{fig:5} (left) clearly demonstrate that the proposed well-balanced
central-upwind quadtree scheme accurately captures a small perturbation of the ``lake-at-rest'' steady state and that the symmetry of the
solution is preserved. The ability of the scheme to refine grids where local gradients are sharp can be seen in Figure \ref{fig:5} (right),
where the quadtree grids at the same times $t=0.6$, $0.9$, $1.2$, $1.5$ and $1.8$ are presented.

\begin{figure}[ht!]
\centerline{\includegraphics[height=3.6cm]{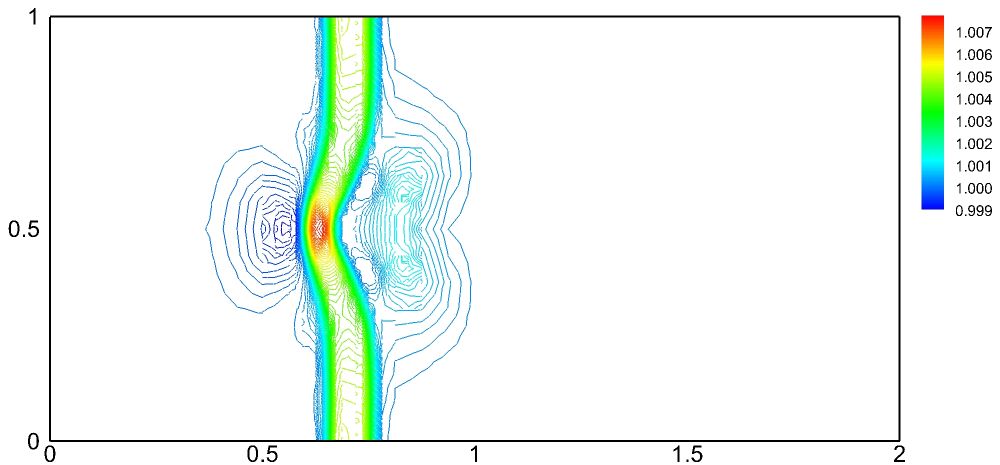}\hspace*{1cm}\includegraphics[height=3.6cm]{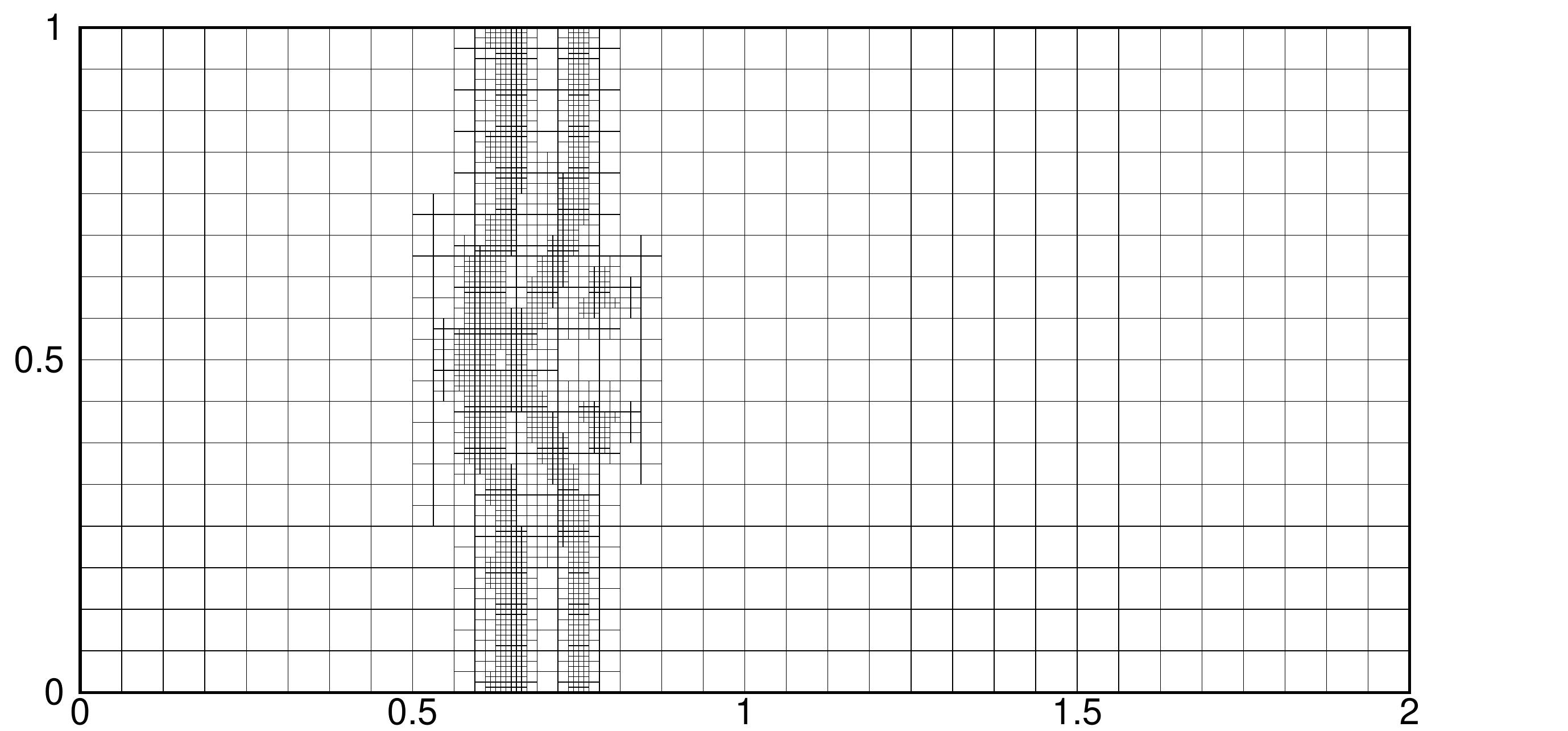}}
\vspace*{0.25cm}
\centerline{\includegraphics[height=3.6cm]{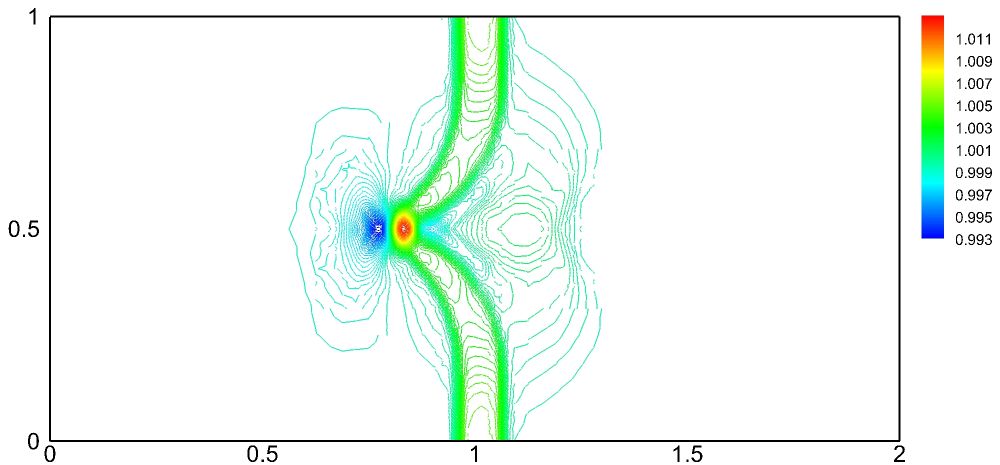}\hspace*{1cm}\includegraphics[height=3.6cm]{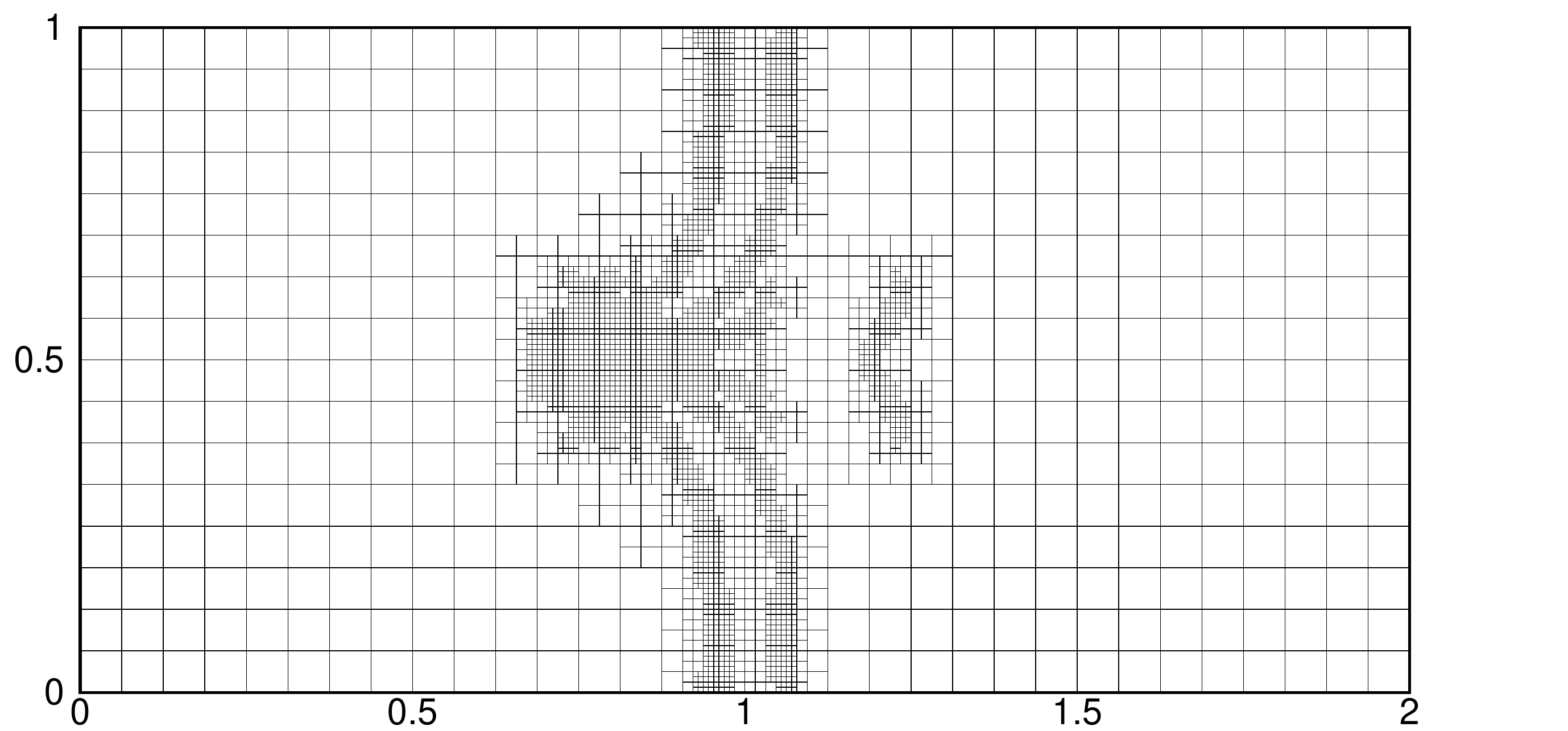}}
\vspace*{0.25cm}
\centerline{\includegraphics[height=3.6cm]{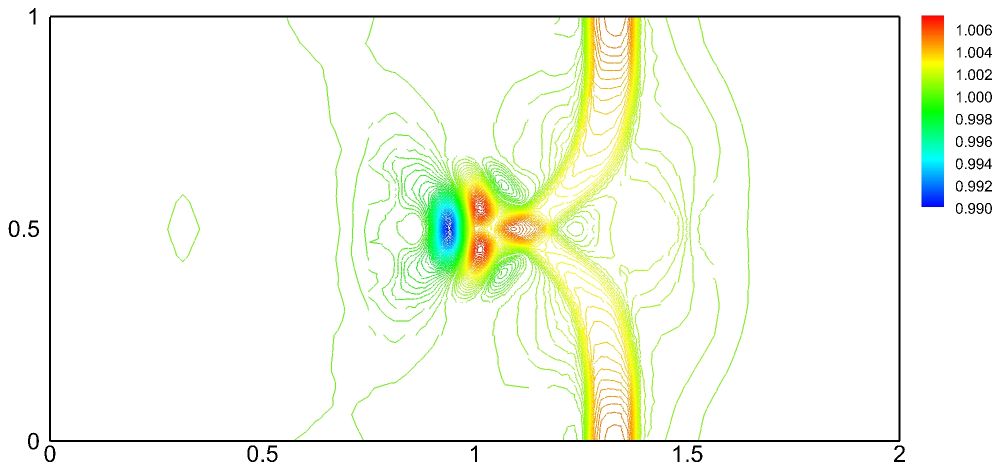}\hspace*{1cm}\includegraphics[height=3.6cm]{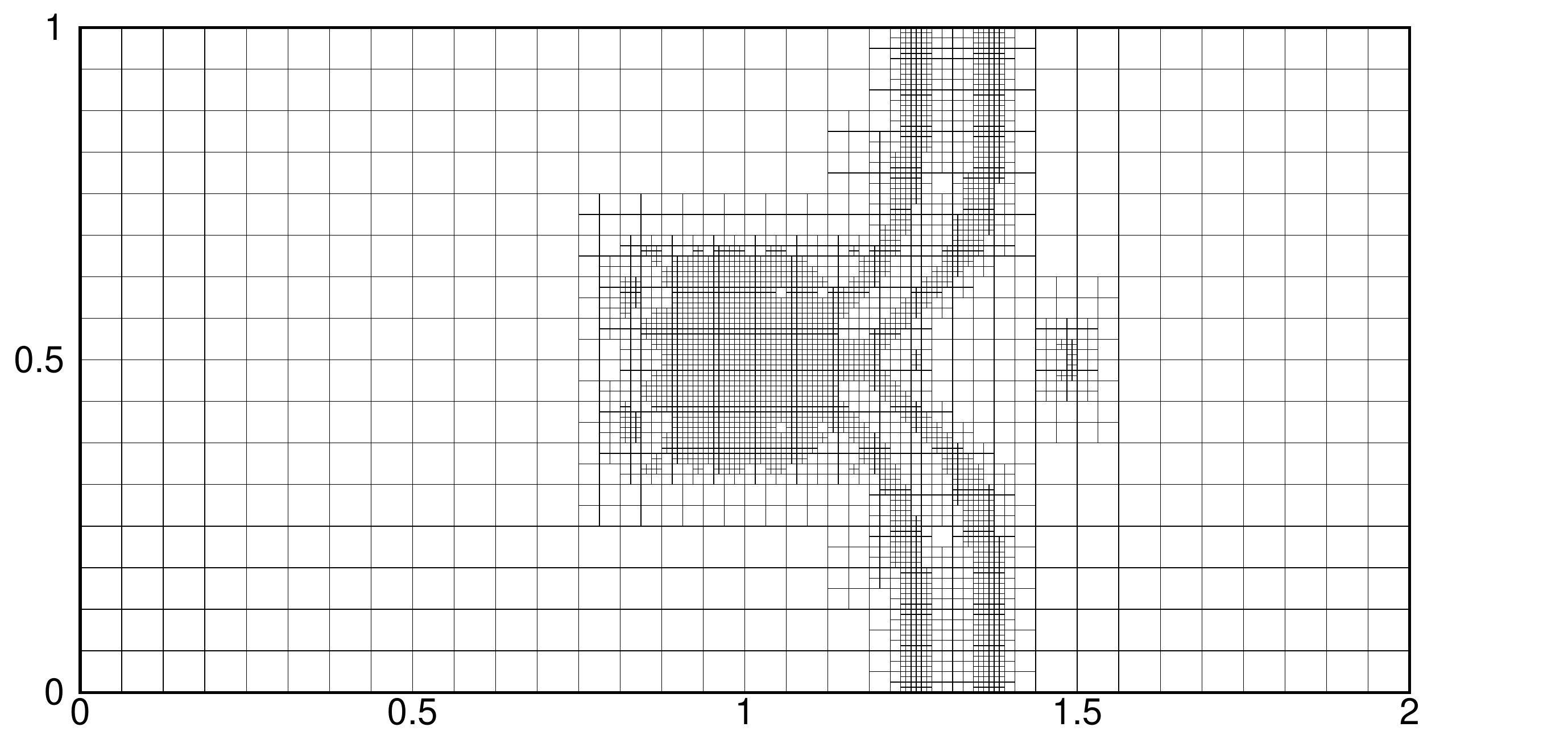}}
\vspace*{0.25cm}
\centerline{\includegraphics[height=3.6cm]{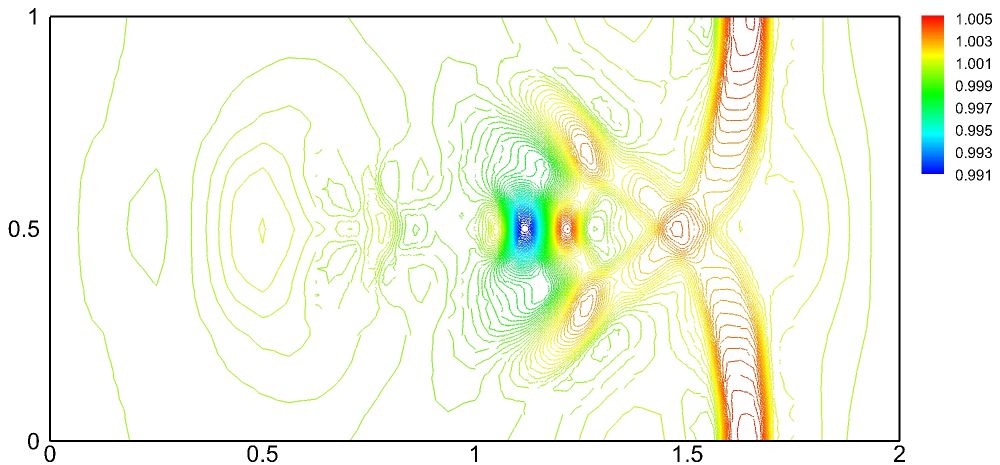}\hspace*{1cm}\includegraphics[height=3.6cm]{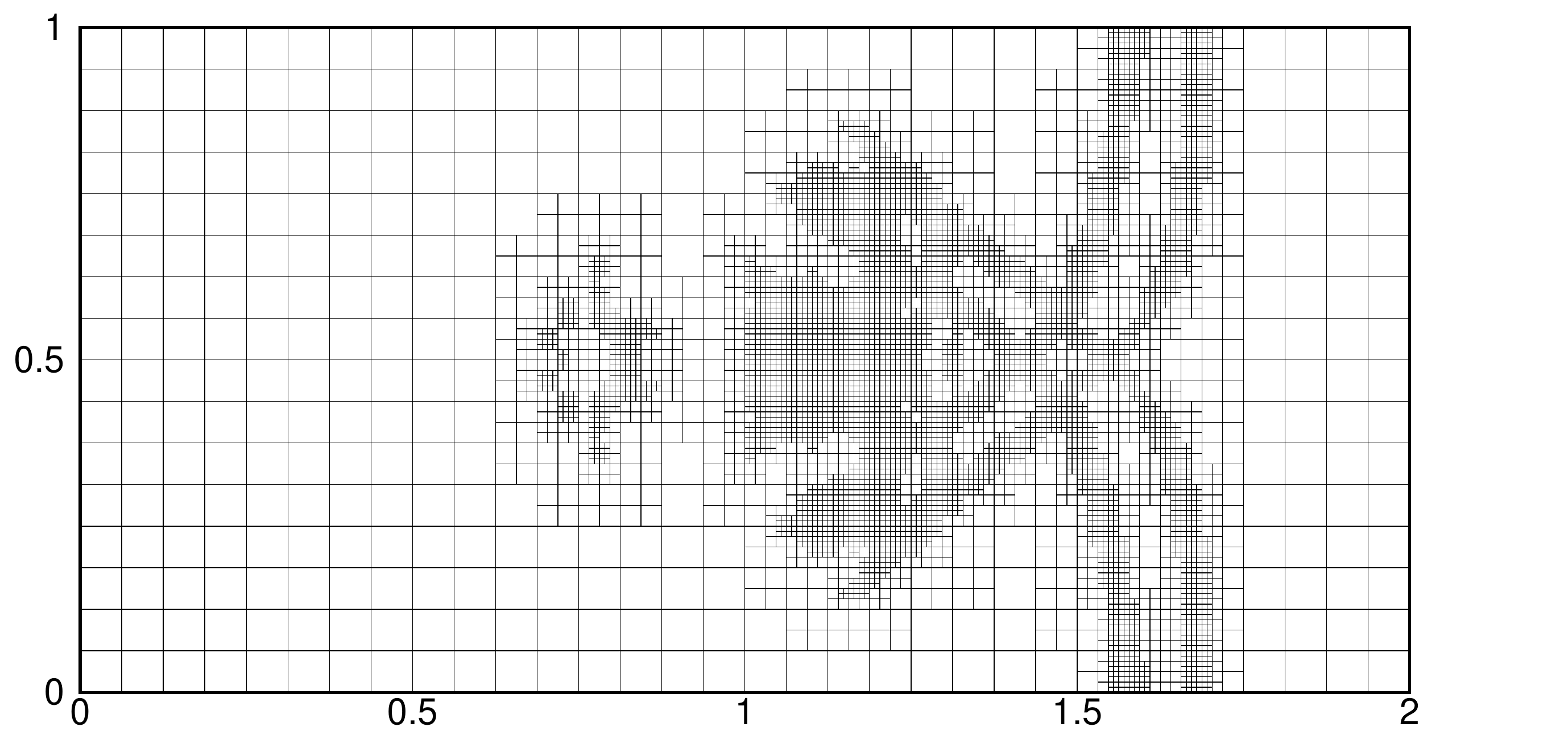}}
\vspace*{0.25cm}
\centerline{\includegraphics[height=3.6cm]{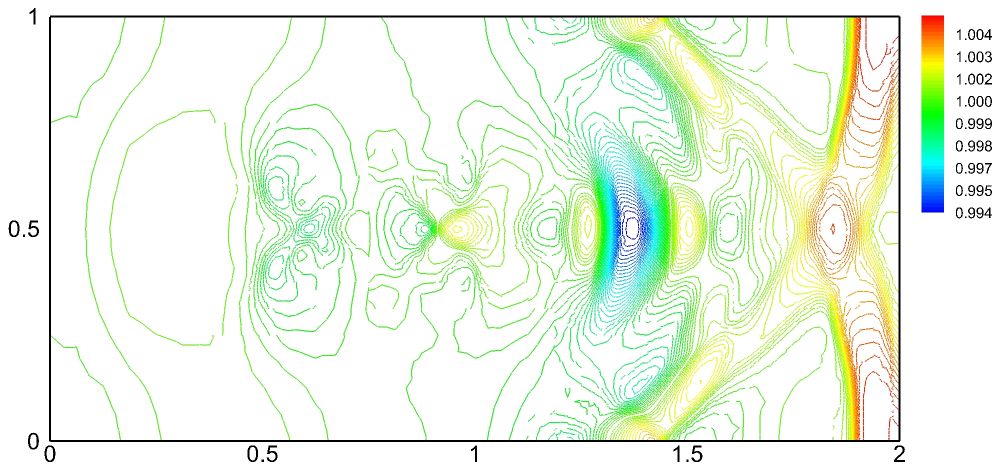}\hspace*{1cm}\includegraphics[height=3.6cm]{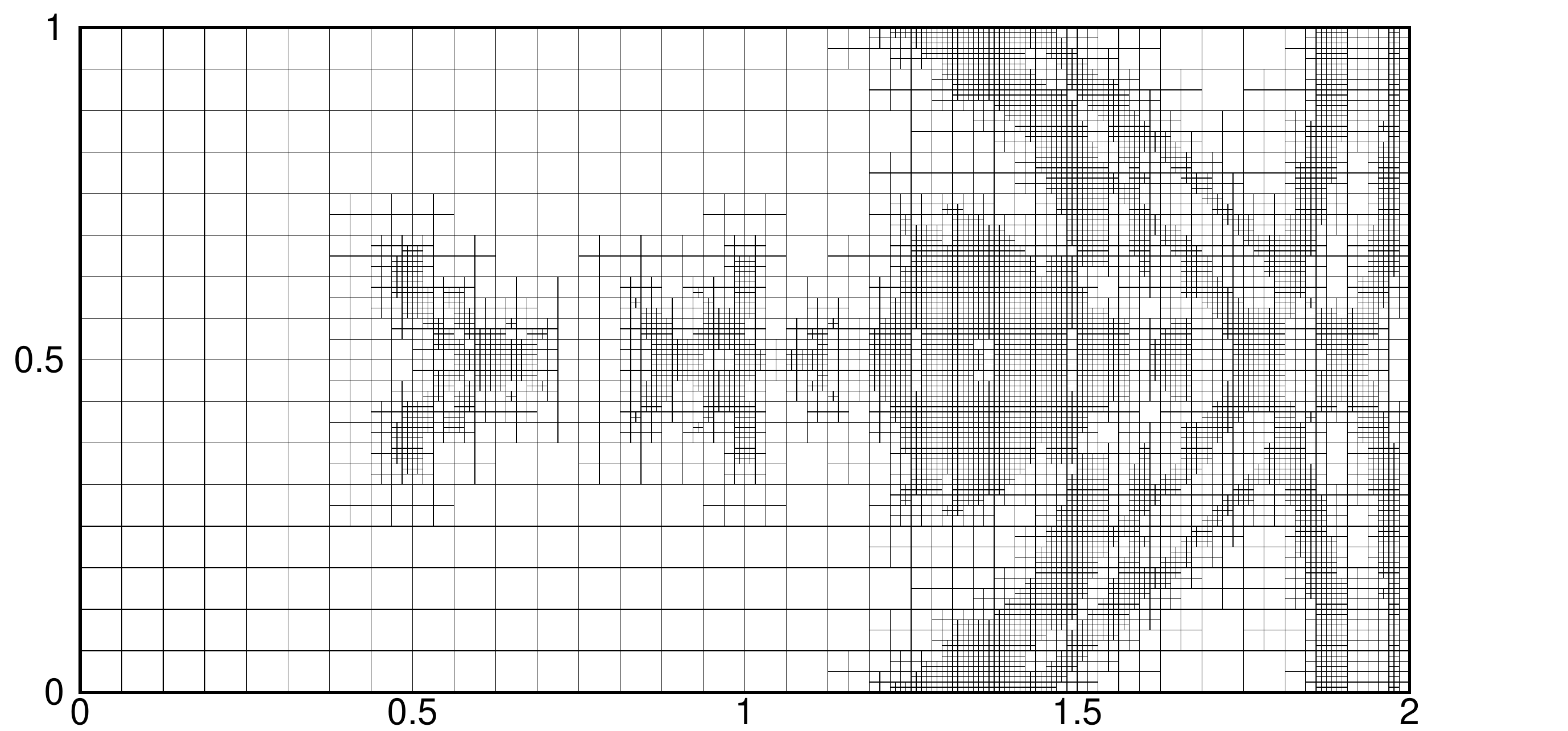}}
\caption{\sf Same as Figure \ref{fig:5}, but for the non-well-balanced central-upwind quadtree scheme.\label{fig:7}}
\end{figure}

We also solve this initial-boundary value problem using a non-well-balanced central-upwind quadtree scheme to stress the importance of the
well-balanced property. In order to design a non-well-balanced scheme, we replace the well-balanced numerical source terms
$\,\xbar S_{j,k}^{\,(2)}$ and $\,\xbar S_{j,k}^{\,(3)}$ given by \eqref{3.17} with the source terms obtained by a straightforward midpoint
rule quadrature. For the configuration considered in Figure \ref{fig:3} (b), the non-well-balanced source term approximations read as
\begin{align*}
&\xbar S_{j,k}^{\,(2)}=
-\frac{g(\,\xbar w_{j,k}-B_{j,k})}{\dx_{j,k}}\left[\frac{B_{\jph,\kph}+B_{\jph,\kmh}}{2}-\frac{B_{\jmh,\kph}+B_{\jmh,\kmh}}{2}\right],\\
&\xbar S_{j,k}^{\,(3)}=
-\frac{g(\,\xbar w_{j,k}-B_{j,k})}{\dx_{j,k}}\left[\frac{B_{\jph,\kph}+B_{\jph,\kmh}}{2}-\frac{B_{\jmh,\kph}+B_{\jmh,\kmh}}{2}\right].
\end{align*}
Figure \ref{fig:7} shows the the snapshots of $w$ at times $t=0.6$, $0.9$, $1.2$, $1.5$ and $1.8$ and the corresponding quadtree grids
obtained using the non-well-balanced computations. As one can see, the use of non-well-balanced numerical source term leads to the
appearance of not small ``parasitic'' waves. Even though these waves are not as large as in the non-well-balanced results presented in,
e.g., \cite{Bryson2011} or \cite{Liu2015}, the unphysical oscillations caused by the non-well-balanced discretization of the source term are
attenuated by adding more seeding points as the quadtree grid reaches a maximum number of 8900 cells during the time evolution. This
demonstrates the importance of the well-balanced property, which eventually reduces the computational cost.

\subsection*{Example 4 --- Small perturbations over a submerged flat plateau}
In this example, which is similar to the examples considered in \cite{Bryson2011,SMSK}, we study small perturbations over a submerged flat
plateau. The computational domain is $[0,1]\times[0,1]$. A solid wall boundary condition is used at the top and  bottom boundaries and
zero-order extrapolation is implemented at the left and right ones. The bottom topography function is given by
$$
B(x,y)=\left\{\begin{array}{lc}1-2\varepsilon,&r\le0.1,\\10(1-2\varepsilon)(0.2-r),&0.1\le r\le0.2,\\0,&~\mbox{otherwise},\end{array}\right.
$$
where $\varepsilon=10^{-4}$ and $r=\sqrt{(x-0.5)^2+(y-0.5)^2}$. The following initial data are imposed:
$$
w(x,y,0)=\left\{\begin{array}{lc}1+\varepsilon,&~0.1\le x \le 0.2,\\1,&\mbox{otherwise},\end{array}\right.\quad u(x,y,0)=v(x,y,0)\equiv0.
$$

\begin{figure}[ht!]
\centerline{\includegraphics[height=5.0cm]{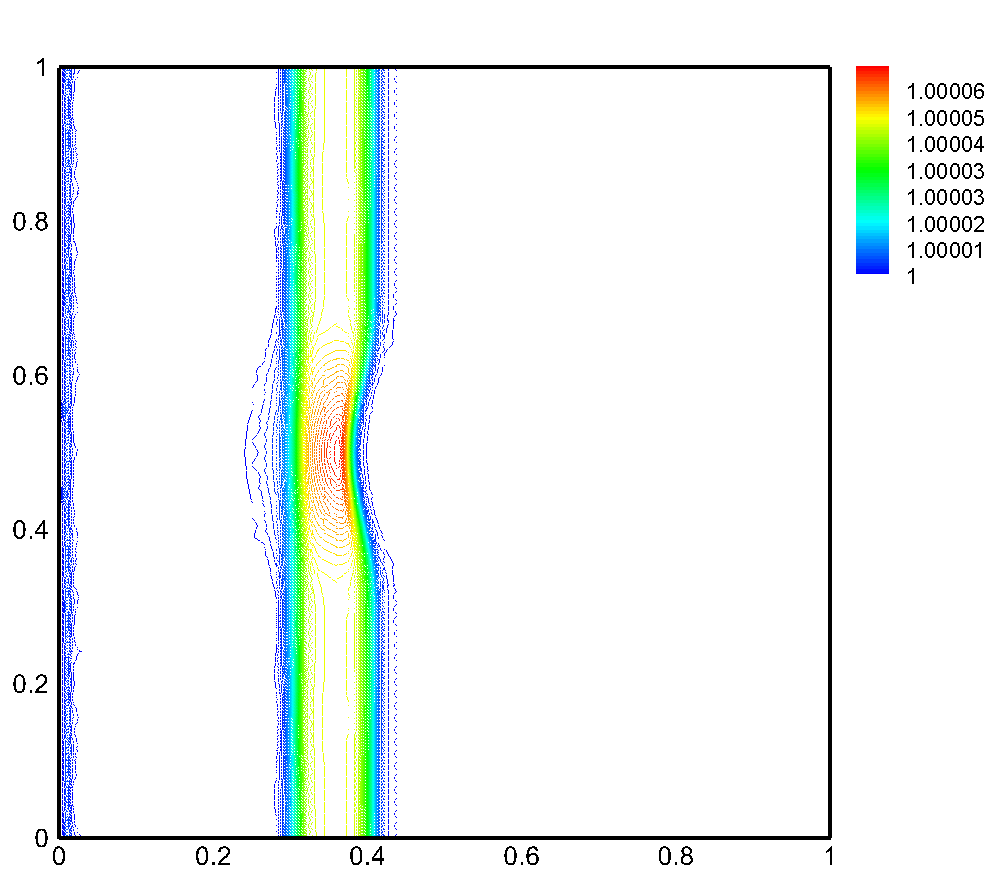}\hspace*{1cm}\includegraphics[height=5.0cm]{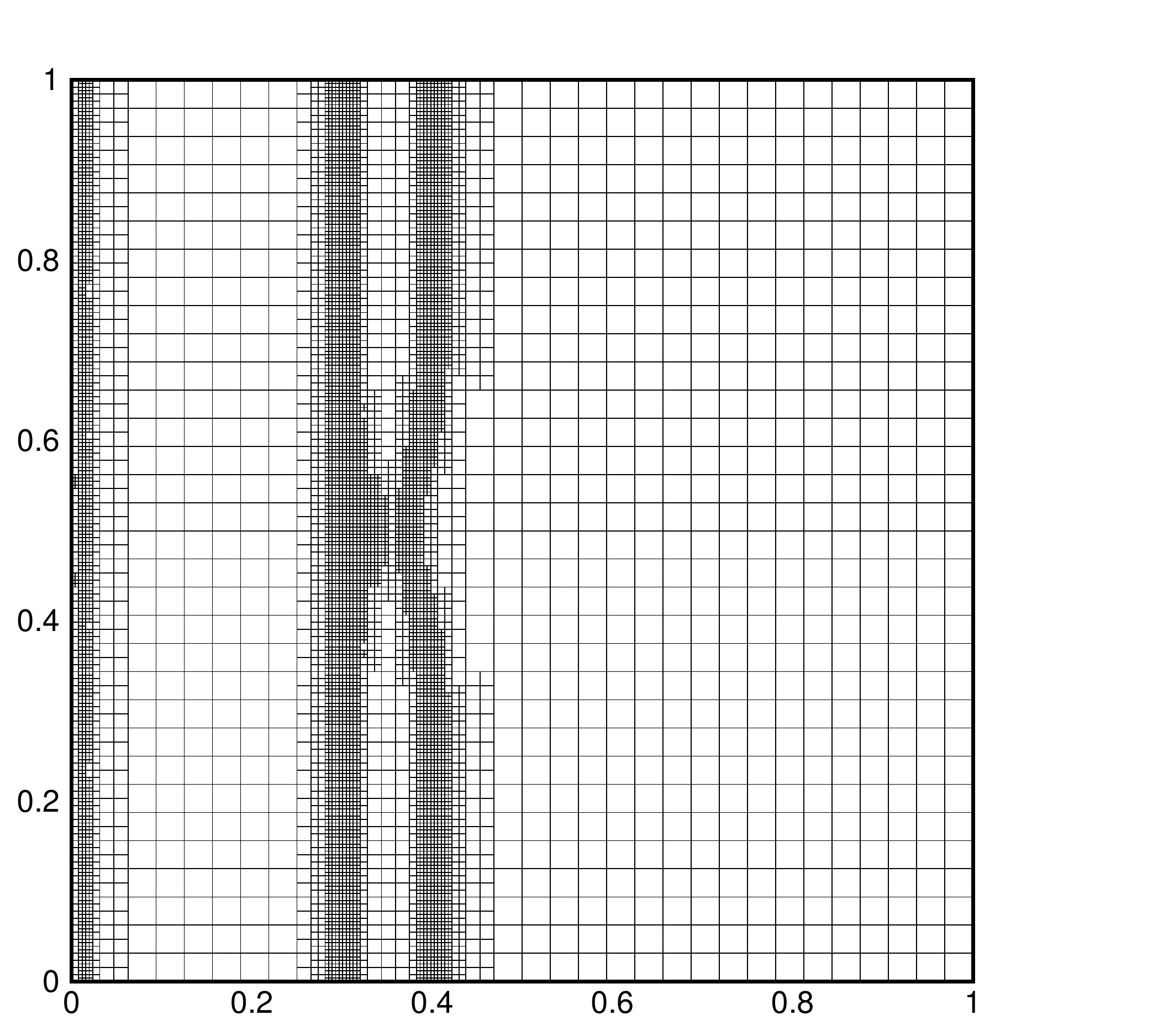}}
\centerline{\includegraphics[height=5.0cm]{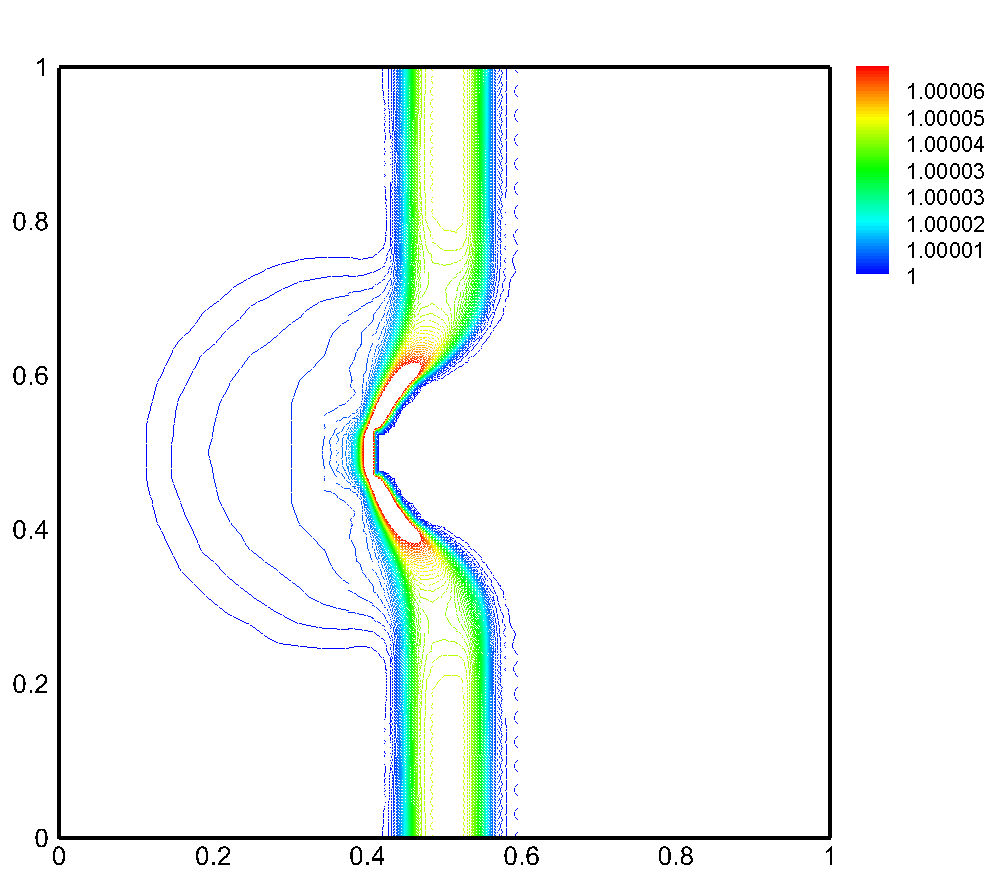}\hspace*{1cm}\includegraphics[height=5.0cm]{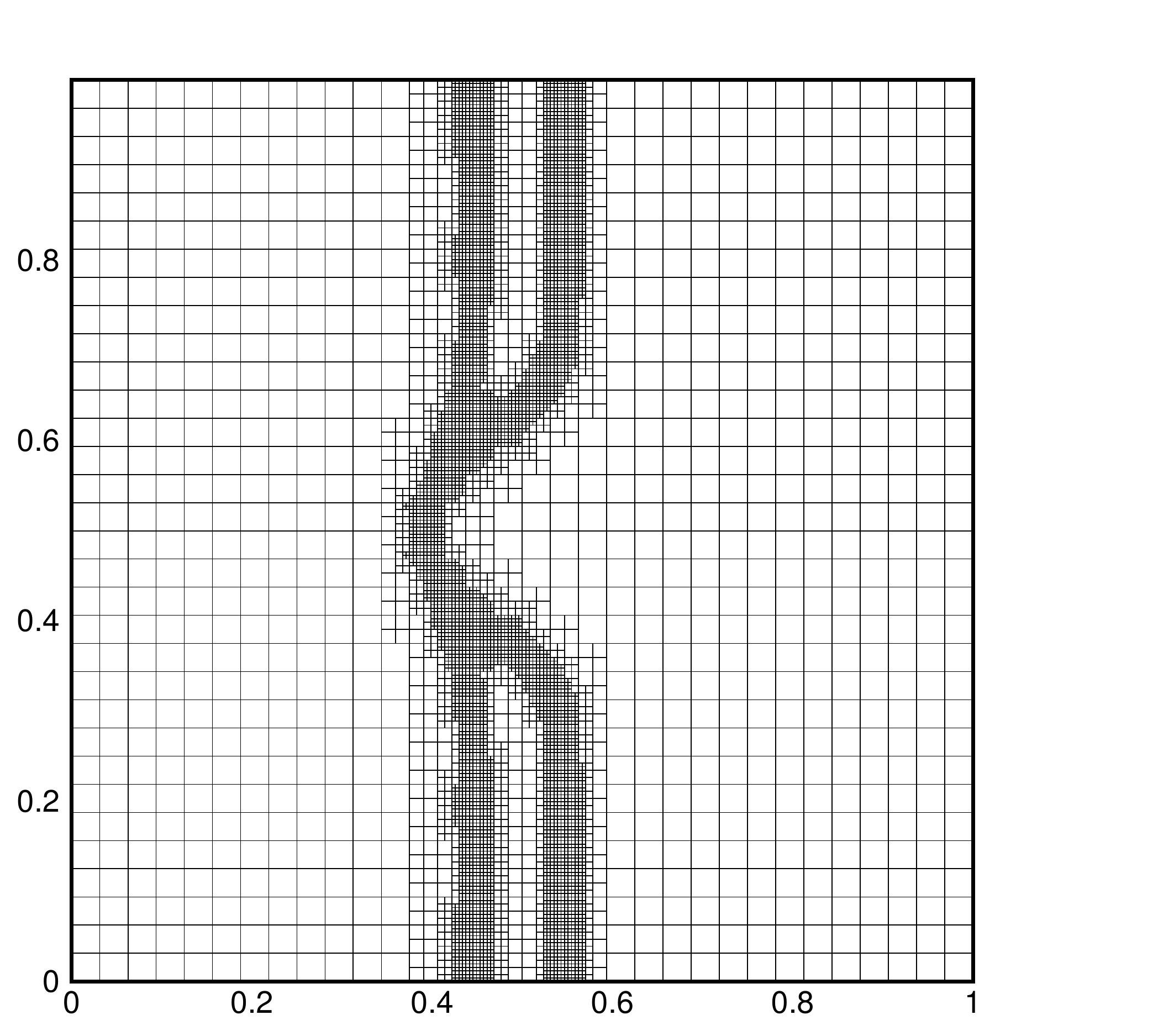}}
\centerline{\includegraphics[height=5.0cm]{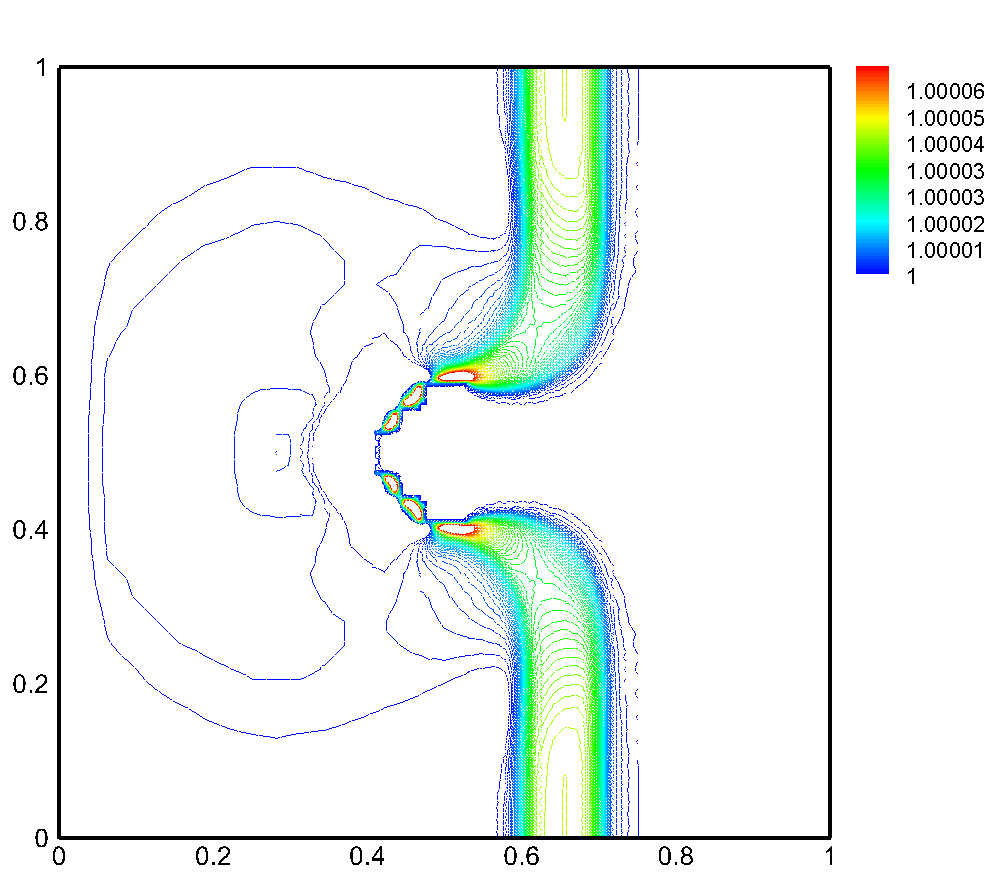}\hspace*{1cm}\includegraphics[height=5.0cm]{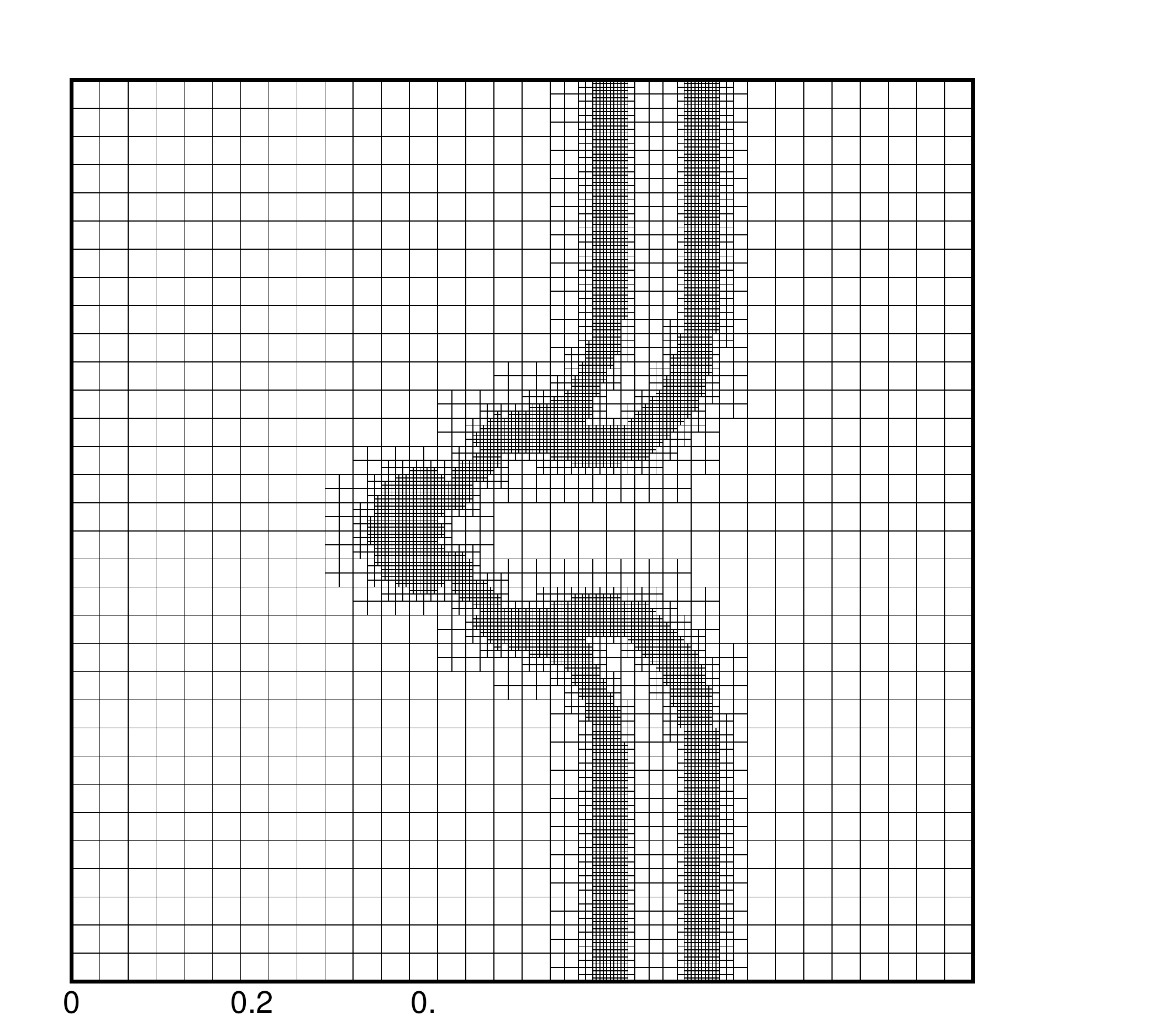}}
\centerline{\includegraphics[height=5.0cm]{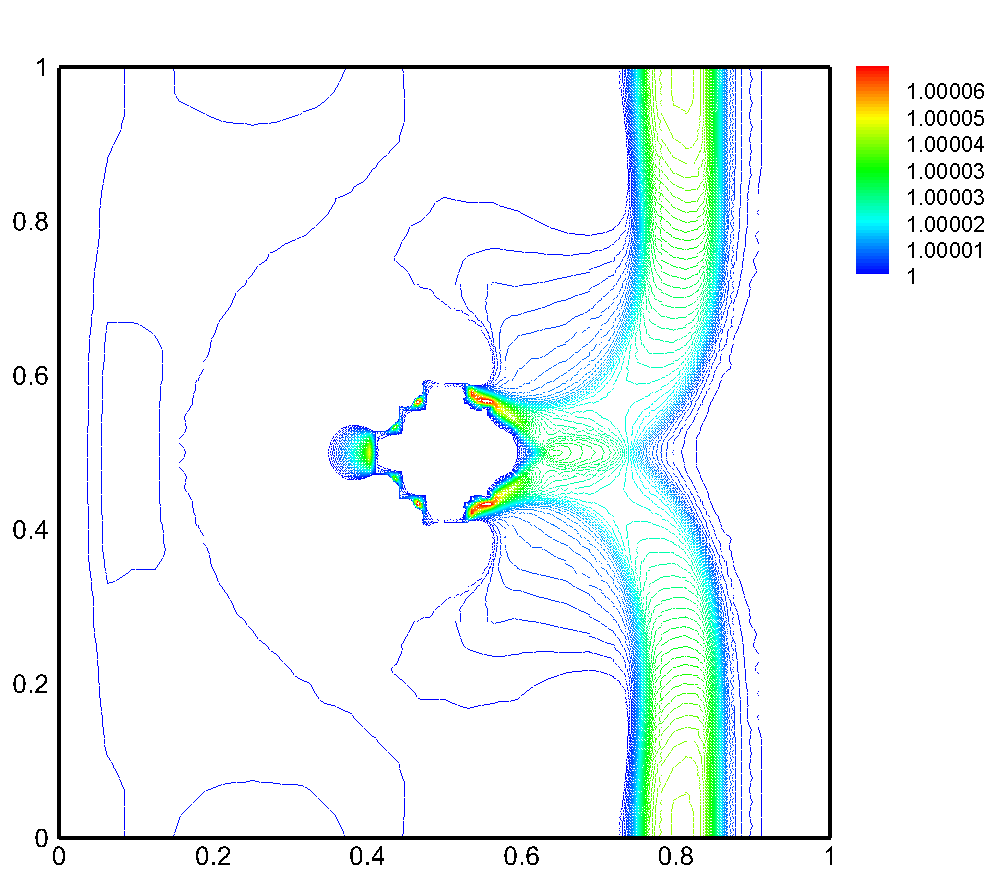}\hspace*{1cm}\includegraphics[height=5.0cm]{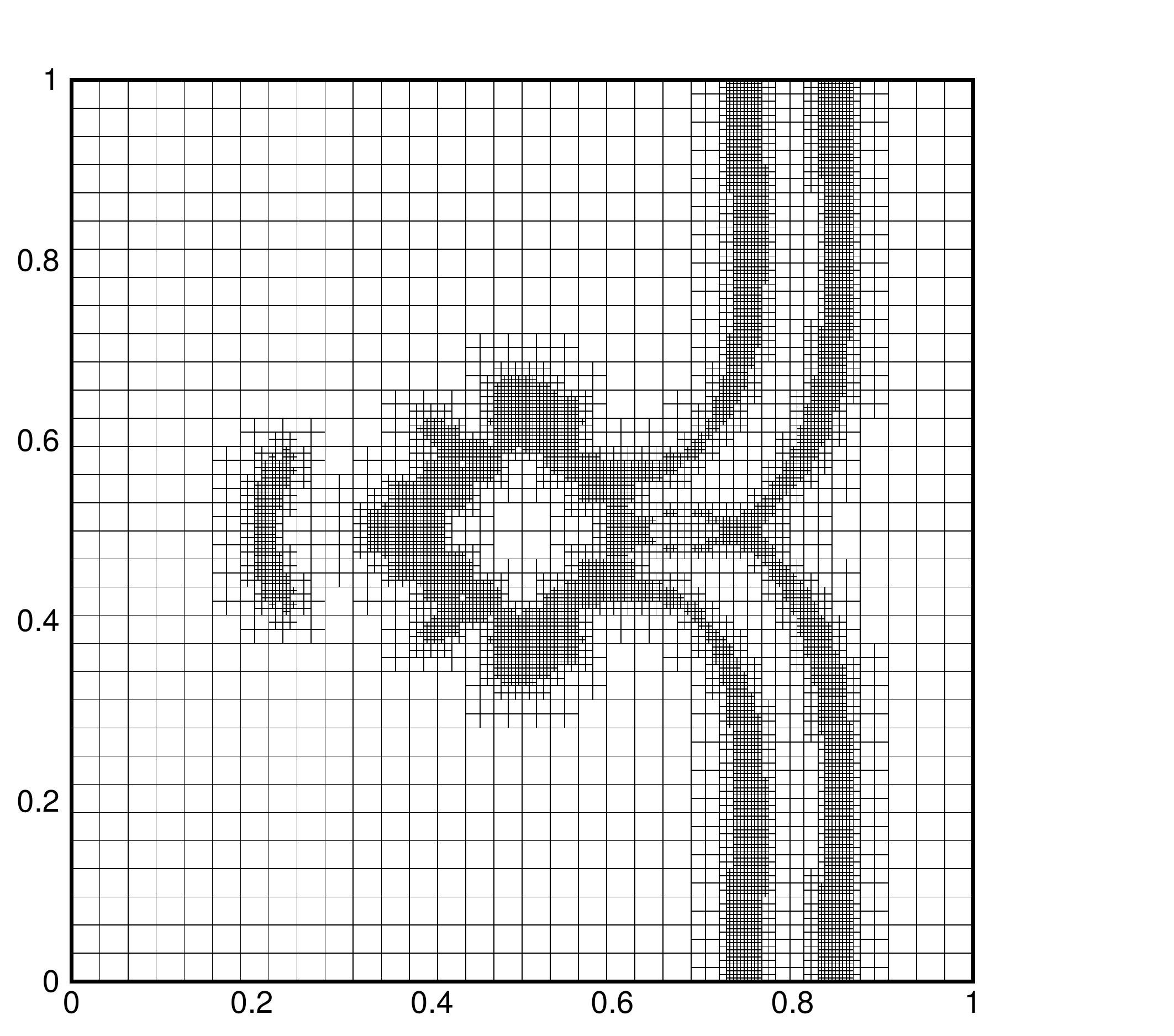}}
\caption{\sf Example 4: Computed water surface $w(x,y,t)$ (left column) and the corresponding quadtree grids (right column) for $t=0.2$,
$0.35$, $0.5$, and $0.65$ (from top to down) obtained using the well-balanced central-upwind quadtree scheme.\label{fig:10.5}}
\end{figure}

\begin{figure}[ht!]
\centerline{\includegraphics[height=5.0cm]{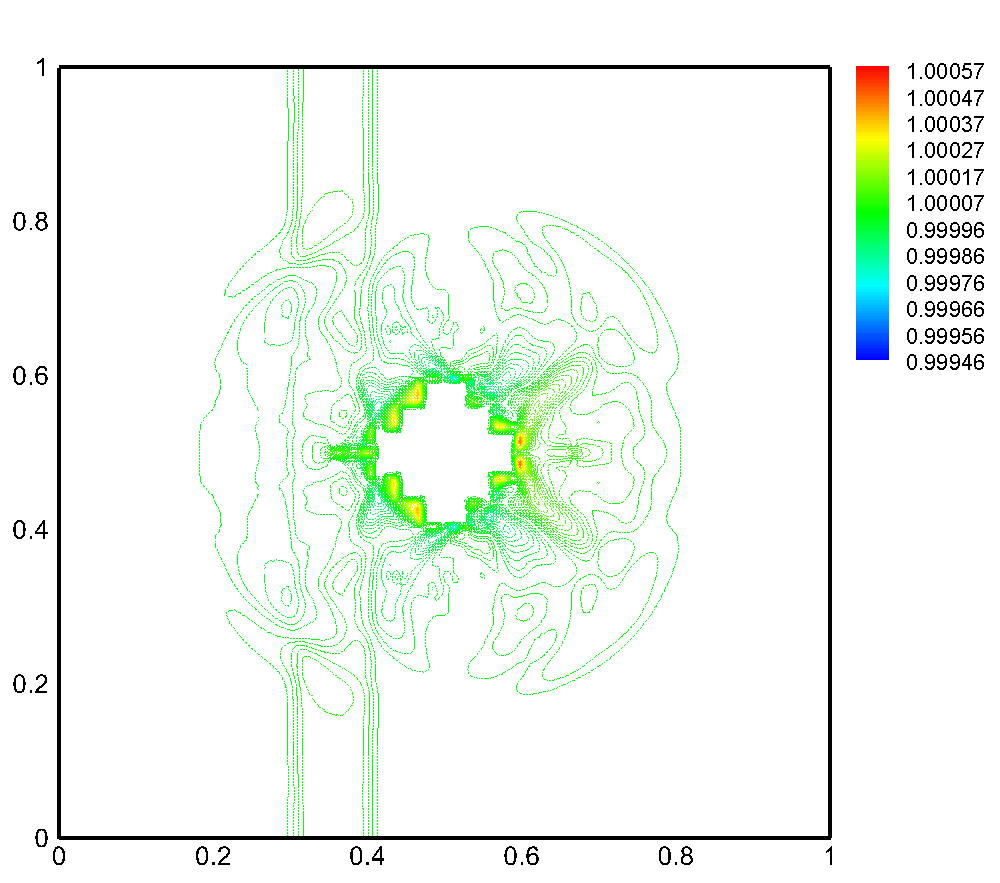}\hspace*{1cm}\includegraphics[height=5.0cm]{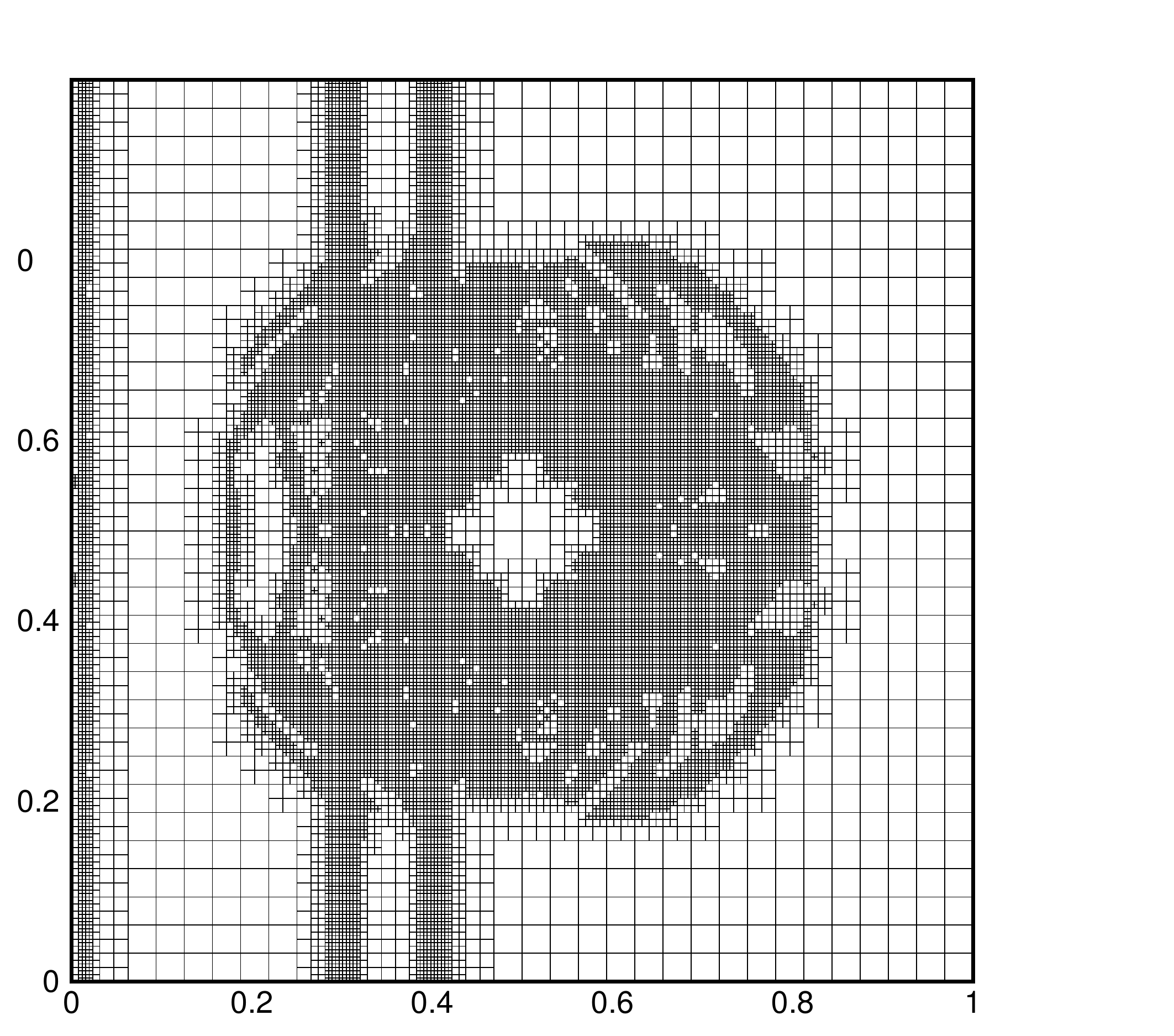}}
\centerline{\includegraphics[height=5.0cm]{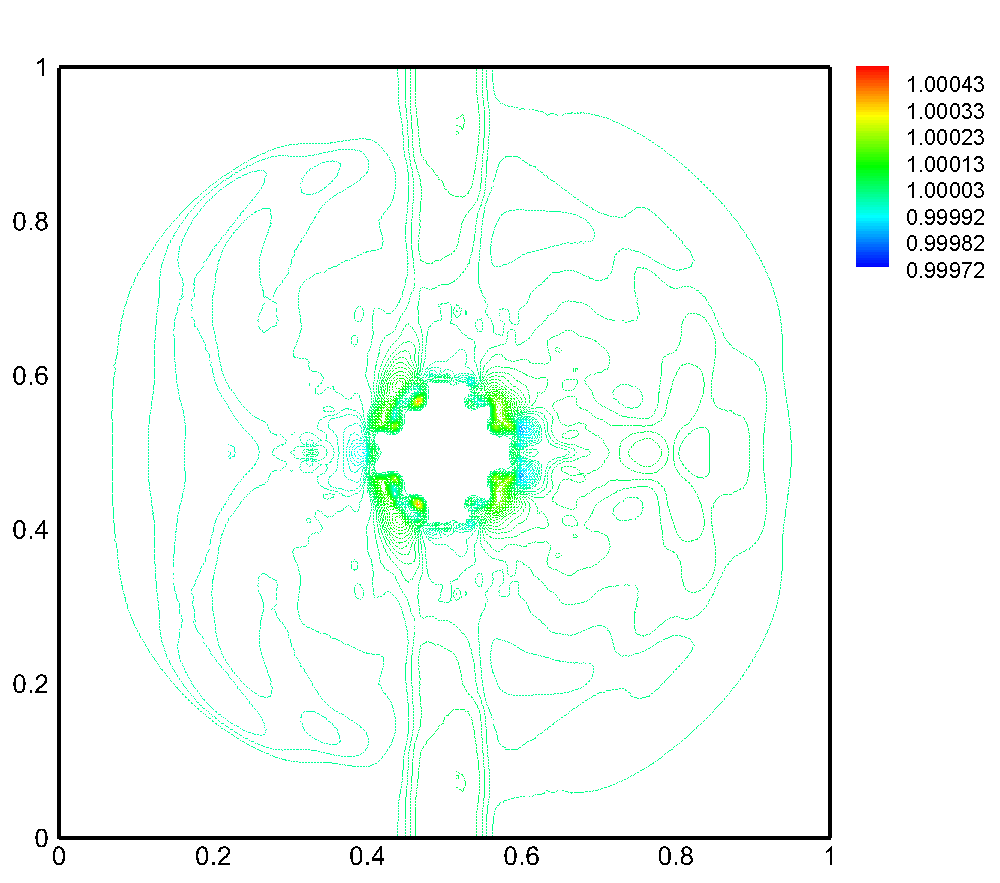}\hspace*{1cm}\includegraphics[height=5.0cm]{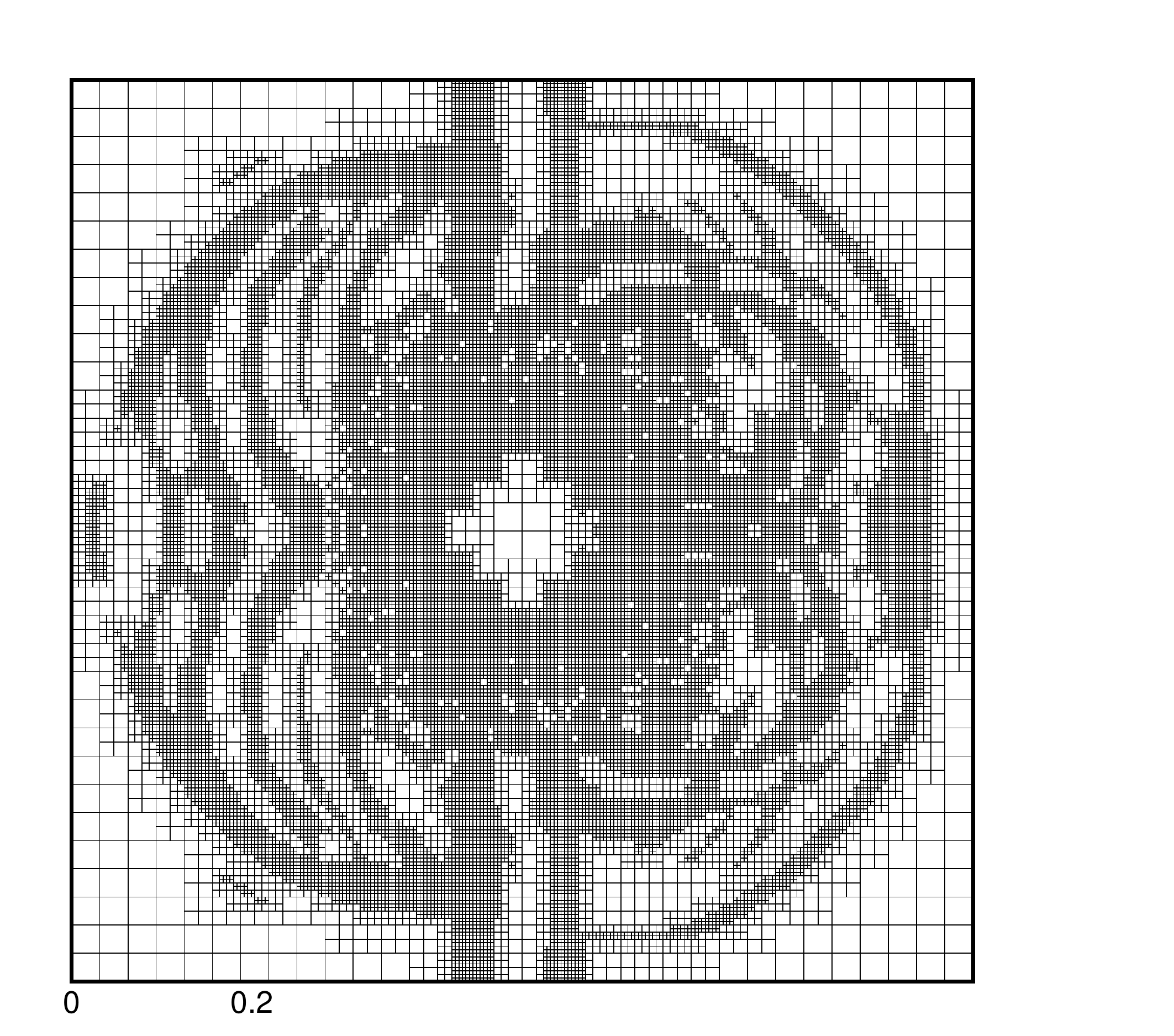}}
\centerline{\includegraphics[height=5.0cm]{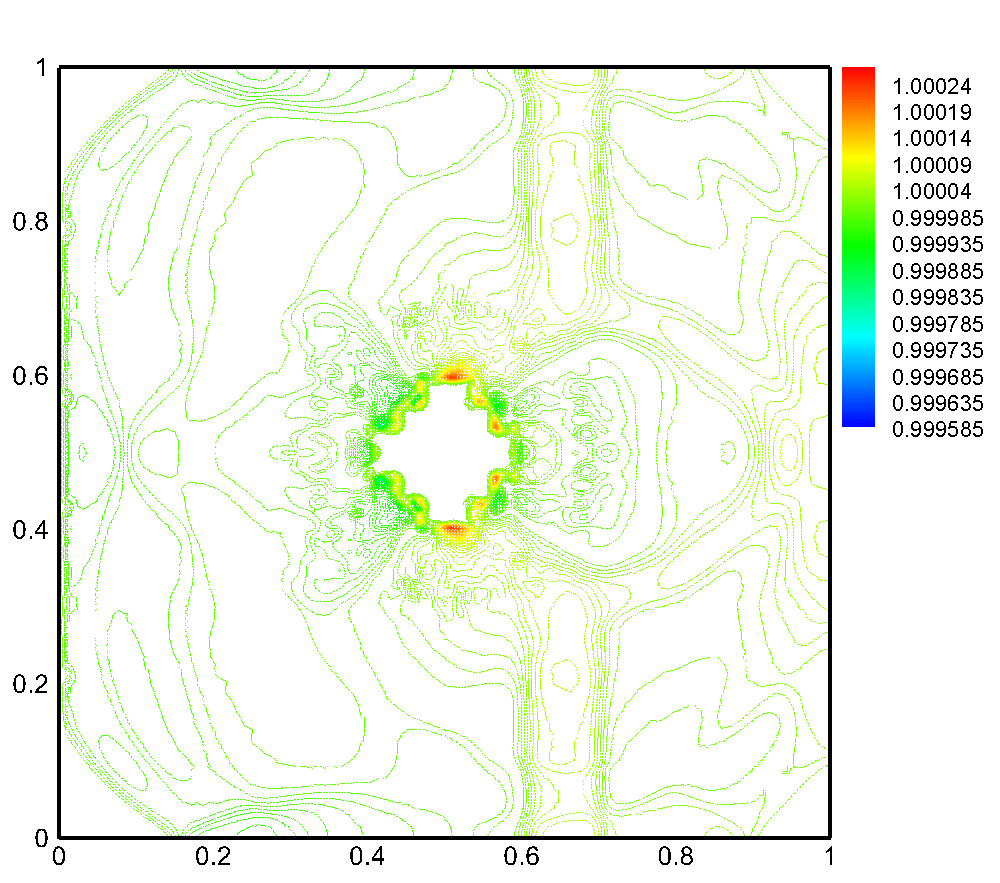}\hspace*{1cm}\includegraphics[height=5.0cm]{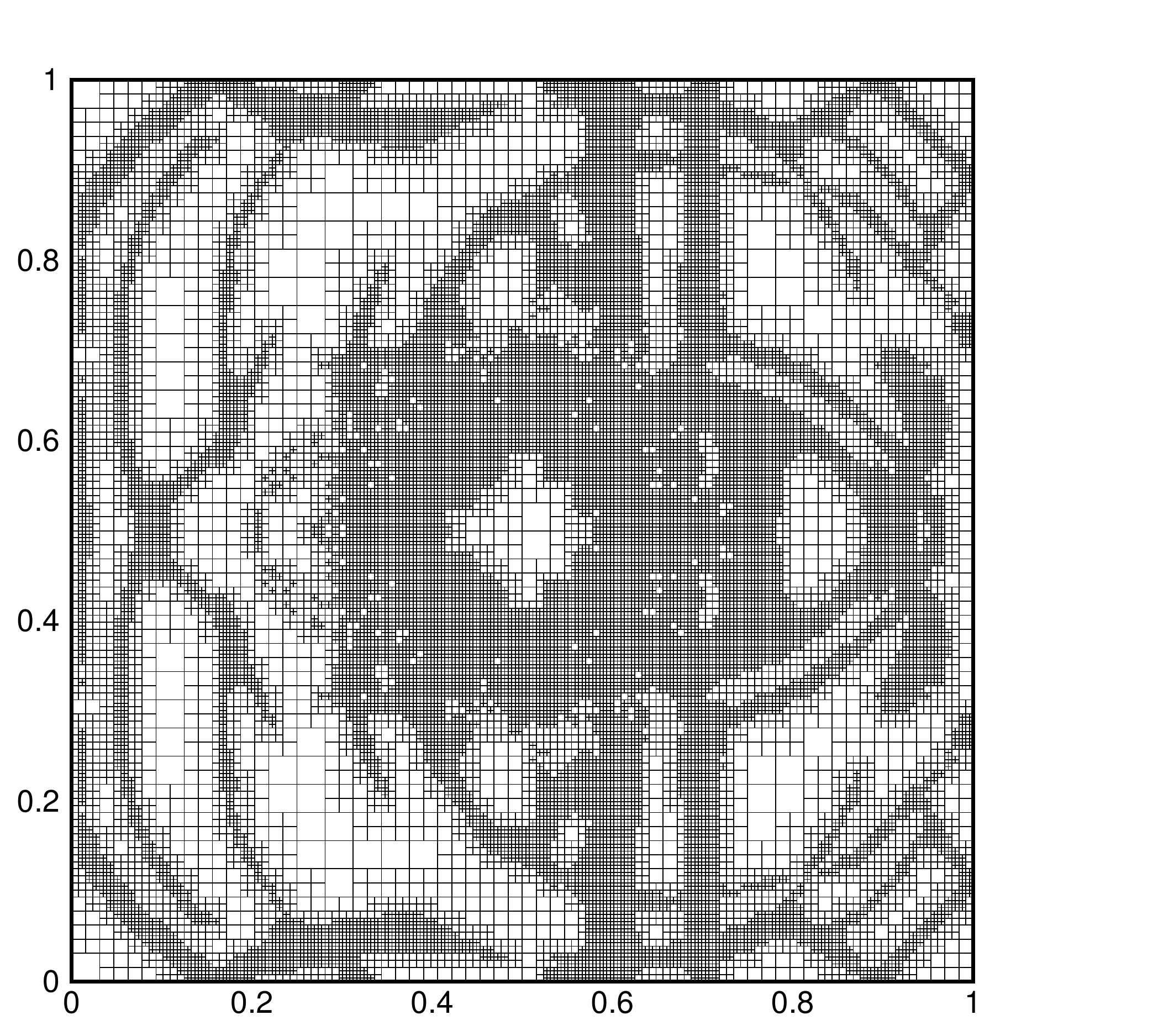}}
\centerline{\includegraphics[height=5.0cm]{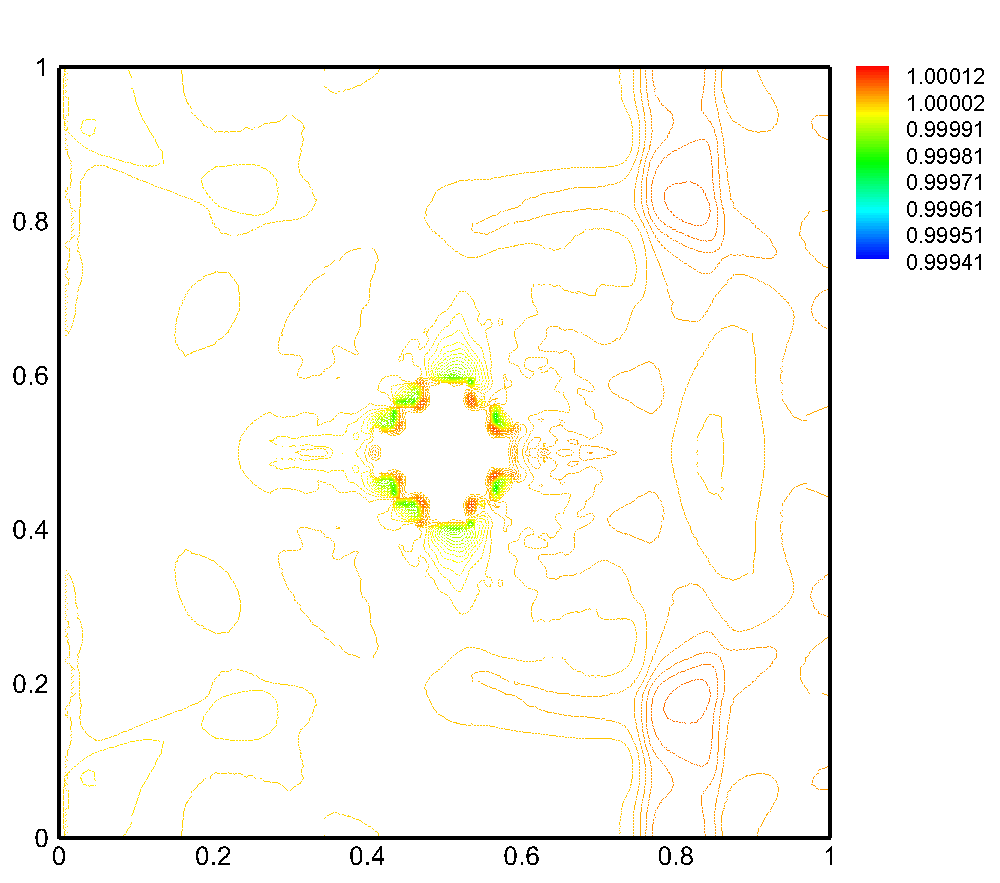}\hspace*{1cm}\includegraphics[height=5.0cm]{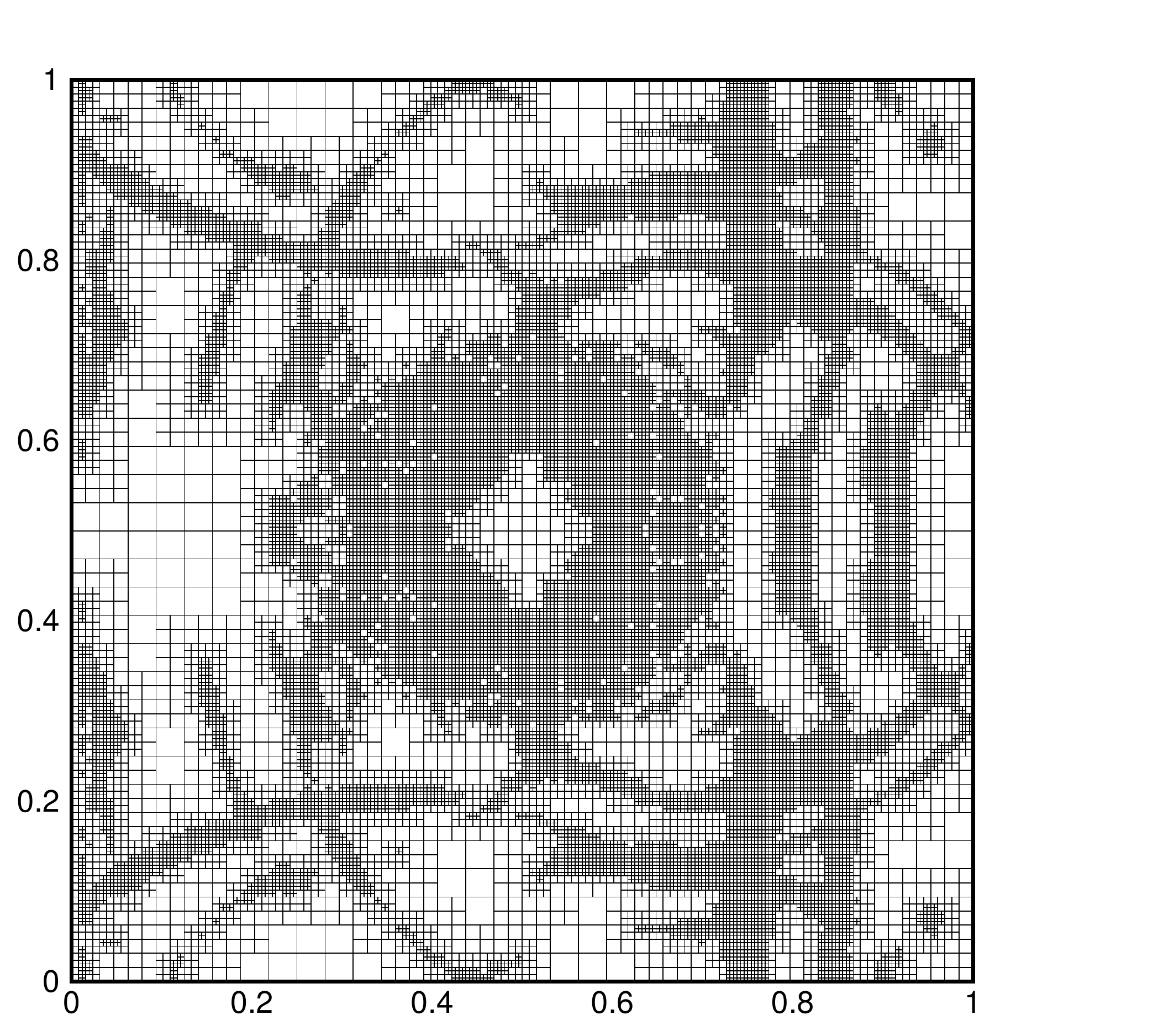}}
\caption{\sf Same as Figure \ref{fig:10.5}, but for the non-well-balanced central-upwind quadtree scheme.\label{fig:10.5.2}}
\end{figure}

We compute both well-balanced and non-well-balanced solutions with $m=8$ and $C_{\rm seed}=0.0002$. The obtained $w$ (left column) and the
corresponding quadtree grids (right column) at $t=0.2, 0.35, 0.5$, and $0.65$ are shown in Figures \ref{fig:10.5} and \ref{fig:10.5.2}. We
note that the number of cells in the well-balanced computation varies from 3712 to 13384, while in the non-well-balanced one it goes up to a
much larger maximum of 34126 cells. However, this level of refinement is apparently not enough to suppress the non-physical parasitic waves,
which propagate all over the computational domain; see Figure \ref{fig:10.5.2}. On the contrary, the well-balanced solution is
oscillation-free as one can clearly see in Figure \ref{fig:10.5}.

\subsection*{Example 5 --- Cylindrical dam break over a step}
\begin{figure}[ht!]
\centerline{\includegraphics[height=5.0cm]{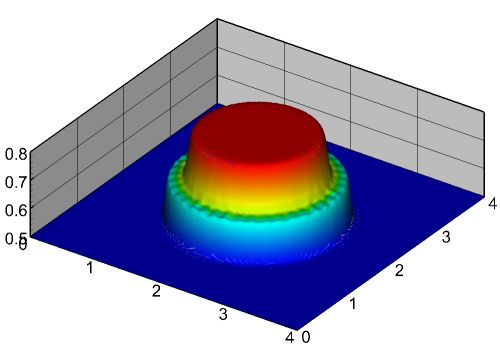}\hspace*{1cm}\includegraphics[height=5.0cm]{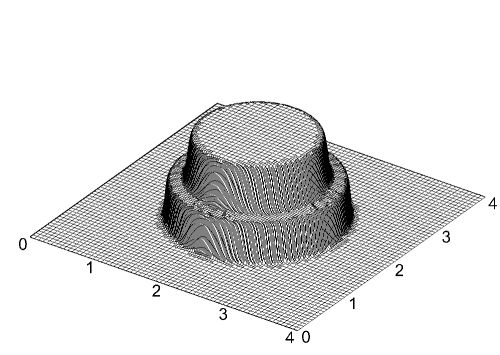}}
\centerline{\includegraphics[height=5.0cm]{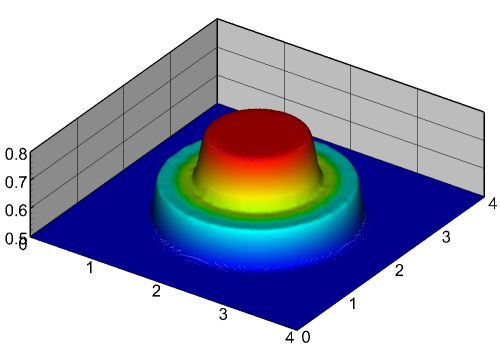}\hspace*{1cm}\includegraphics[height=5.0cm]{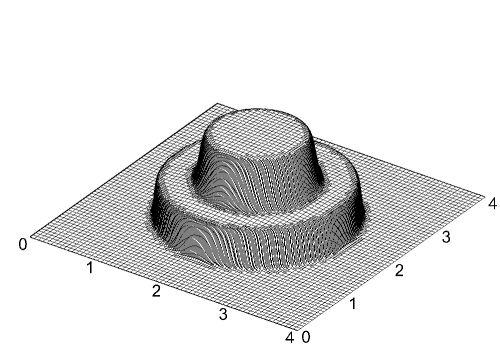}}
\centerline{\includegraphics[height=5.0cm]{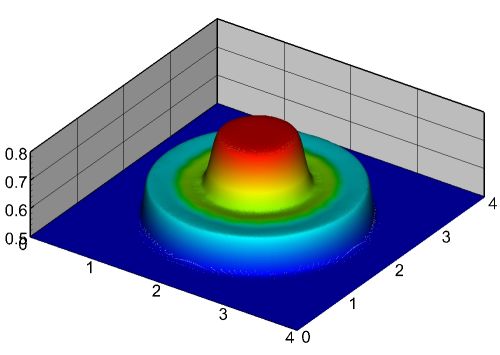}\hspace*{1cm}\includegraphics[height=5.0cm]{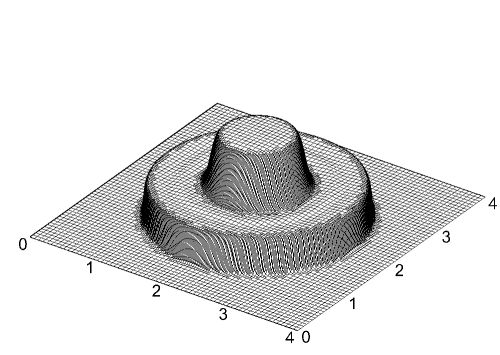}}
\centerline{\includegraphics[height=5.0cm]{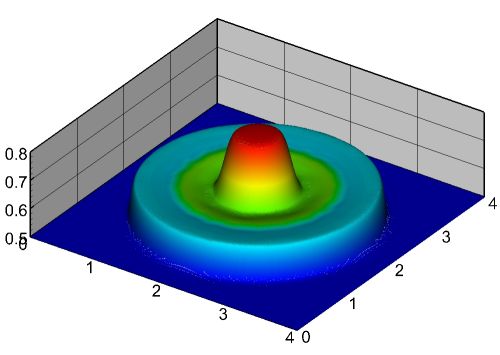}\hspace*{1cm}\includegraphics[height=5.0cm]{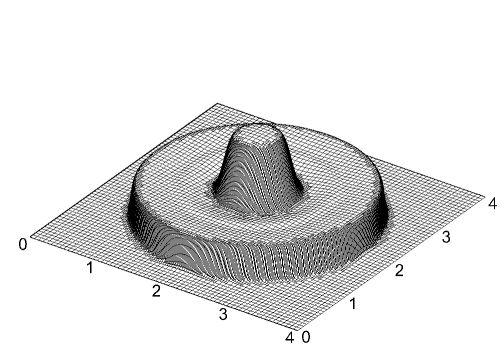}}
\caption{\sf Example 5: Computed water surface $w(x,y,t)$ (left column) and the corresponding quadtree grids (right column) for $t=0.05$,
$0.1$, $0.15$, and $0.2$ (from top to down).\label{fig:7.5}}
\end{figure}

In this example taken from \cite{DUMBSER20138057}, we demonstrate the capability of the proposed scheme to solve problems with discontinuous
bottom topography. We consider the following initial conditions and bottom topography:
$$
h(x,y,0)=\left\{\begin{aligned}&1,&&r\le1,\\&0.5,&&r>1,\end{aligned}\right.\quad u(x,y,0)=v(x,y,0)\equiv0,\quad
B(x,y)=\left\{\begin{aligned}&-0.2,&&r\le1,\\&0,&&r>1,\end{aligned}\right.
$$
where $r=\sqrt{(x-2)^2+(y-2)^2}$. The computational domain is $[0,4]\times[0,4]$, $g=9.8$ and the point values of $B$ at the vertices of
$C_{j,k}$ are obtained using the bottom topography function with $\ell:=6$ (\S\ref{S:4.5}). We compute the solution until the final time
$t=0.2$ using $C_{seed}=0.1$ and $m=8$. Figure \ref{fig:7.5} shows the obtained $w$ at times $t=0.05$, 0.1, 0.15 and 0.2, and the
corresponding quadtree grids. The number of cells at time $t=0$ is 6172 and it reaches a maximum of 12508 cells at later times. As one can
see, the proposed scheme is capable of accurately capturing the solution in the case of discontinuous bottom topography. 

\subsection*{Example 6 --- Sudden contraction}
The last example is a modification of the benchmark in \cite{Hubbard2001}; also see \cite{Bryson2005,Bryson2011}. The purpose of this
example is twofold: to show the ability of the central-upwind quadtree scheme to capture shocks and sharp waves in supercritical flows and to
demonstrate the positivity preserving property of the proposed scheme.

We consider an open channel with a sudden contraction. The geometry of the channel is established on its contraction, where
$$
y_b(x)=\left\{\begin{array}{lc}0.5,&x\le1,\\0.4,&\mbox{otherwise}.\end{array}\right.
$$
The computational domain is $[0,3]\times[0.5-y_b(x),{0.5+y}_b(x)]$. Solid wall boundary conditions are imposed at all of the boundaries
except for the left (inflow boundary with $u(0,y,t)\equiv2$) and right (zero-order extrapolation) ones. The following initial conditions are
prescribed:
$$
w(x,y,0)\equiv1,\quad u(x,y,0)\equiv2,\quad v(x,y,0)\equiv0. 
$$
In this example, we take $m=8$ and $m=9$ refinement levels of the quadtree grid and set $C_{\rm seed}=2$ in \eqref{3.28}. This value of
$C_{\rm seed}$ is greater than the ones used in Examples 1 and 2 since this numerical experiment focuses on capturing sharp waves and thus
choosing small values of $C_{\rm seed}$ would have increased the computational cost as the local gradients are relatively large in most
parts of the computational domain.

\begin{figure}[H]
\centerline{\includegraphics[width=8.5cm]{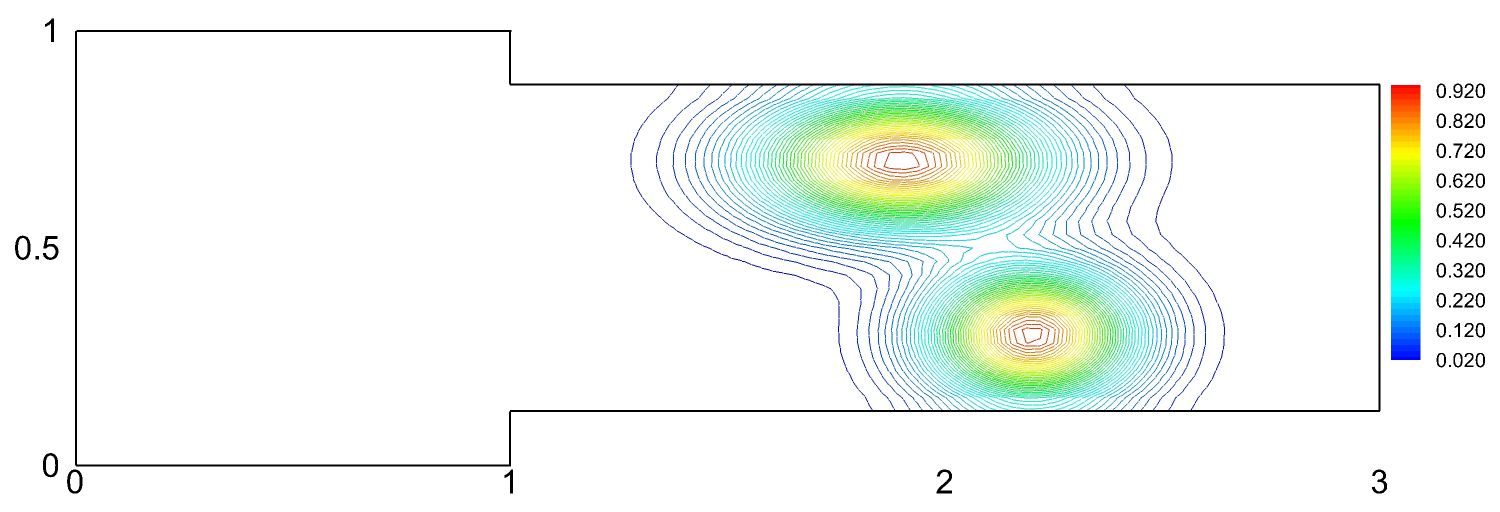}\hspace*{0cm}\includegraphics[width=8.5cm]{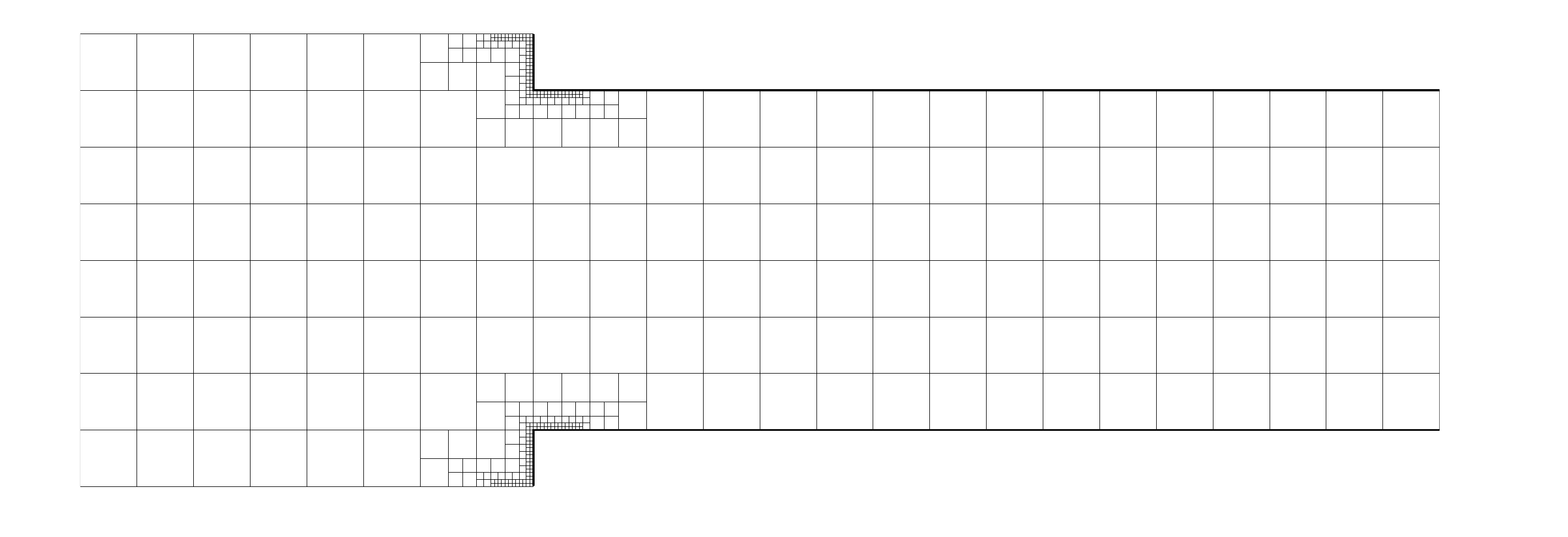}}
\caption{\sf Example 6: Bottom topography (left) and initial quadtree grid with $m=9$ (right).\label{fig:8}}
\end{figure}

We compute the solution twice: first, we use the flat bottom topography $B(x,y)\equiv0$ in order to demonstrate the ability of the scheme
to capture hydraulic jumps and sharp waves, and second, we use the bottom topography given by
$$
B(x,y)=0.95\left[e^{-10(x-1.9)^2-50(y-0.7)^2}+e^{-20(x-2.2)^2-50(y-0.3)^2}\right]
$$
and shown in Figure \ref{fig:8} together with the initial quadtree grid for $m=9$ (notice that the grid is refined near the boundaries at
the contraction to improve accuracy). In the nonflat bottom topography case, the water at the top of the humps is quite shallow (that is why
this is a good example to test the positivity preserving property) and the Froude number there is initially about 2.

\begin{figure}[ht!]
\centerline{\includegraphics[width=8.5cm]{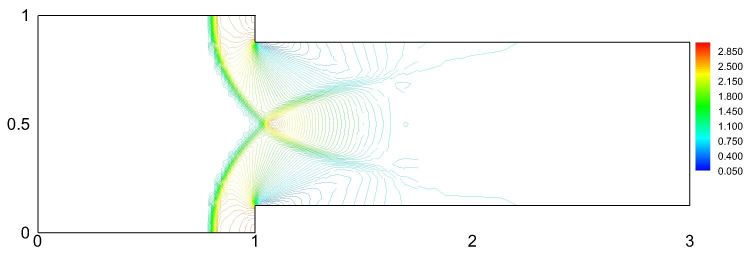}\hspace*{0.2cm}\includegraphics[width=8.5cm]{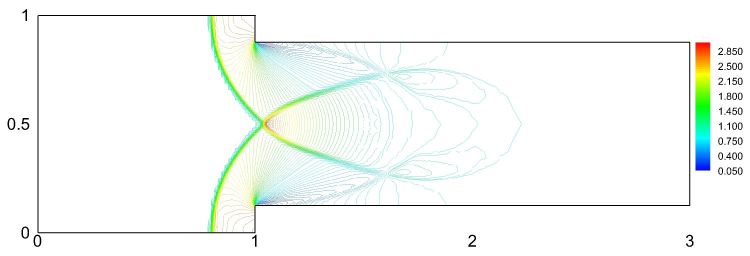}}
\vspace*{0.25cm}
\centerline{\includegraphics[width=8.5cm]{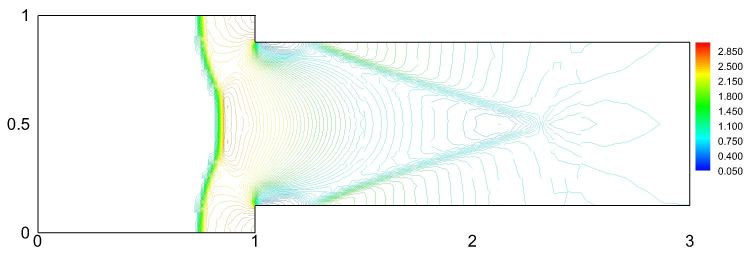}\hspace*{0.2cm}\includegraphics[width=8.5cm]{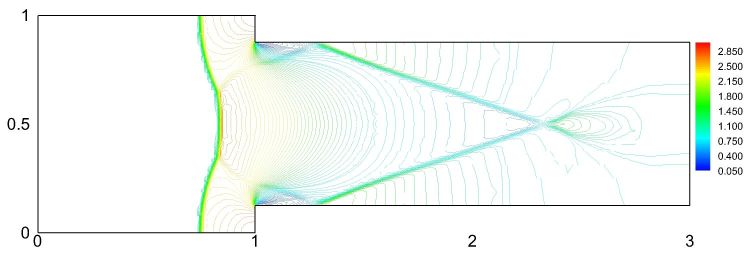}}
\vspace*{0.25cm}
\centerline{\includegraphics[width=8.5cm]{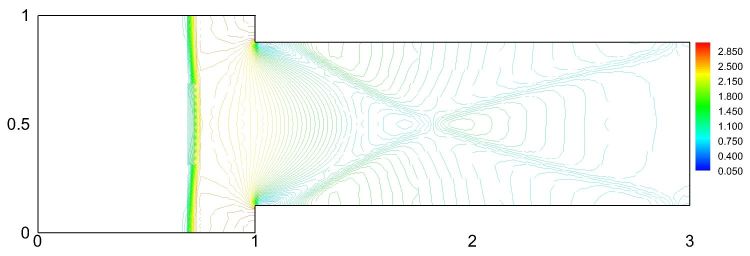}\hspace*{0.2cm}\includegraphics[width=8.5cm]{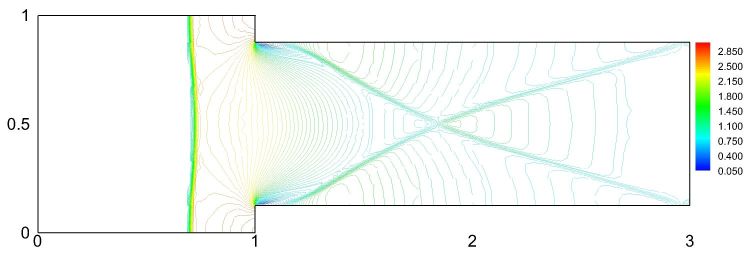}}
\vspace*{0.25cm}
\centerline{\includegraphics[width=8.5cm]{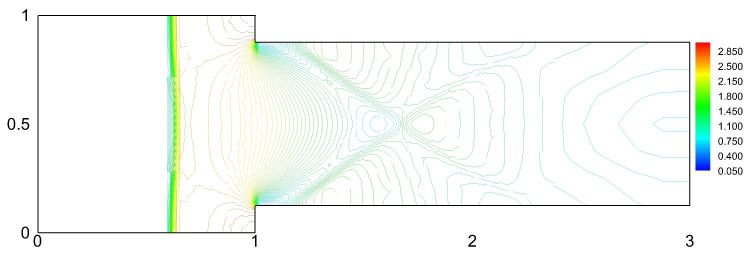}\hspace*{0.2cm}\includegraphics[width=8.5cm]{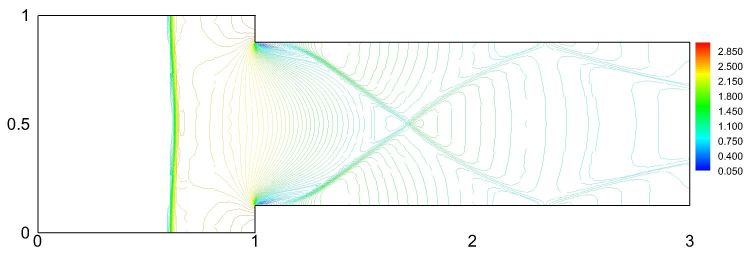}}
\caption{\sf Example 6: Computed water surface $w(x,y,t)$ for $t=0.5$, 1, $1.5$ and 2 (from top to down) obtained using the flat bottom
topography for $m=8$ (left column) and $m=9$ (right column).\label{fig:10}}
\end{figure}

\begin{figure}[ht!]
\centerline{\includegraphics[width=8.5cm]{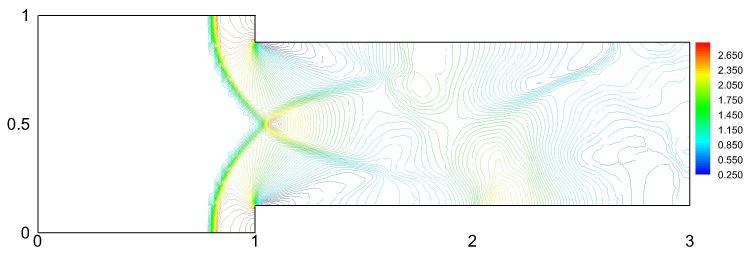}\hspace*{0.2cm}\includegraphics[width=8.5cm]{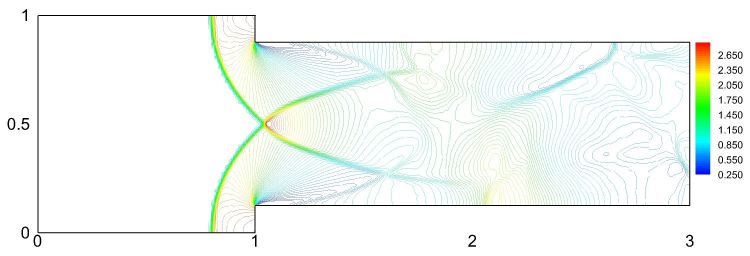}}
\vspace*{0.25cm}
\centerline{\includegraphics[width=8.5cm]{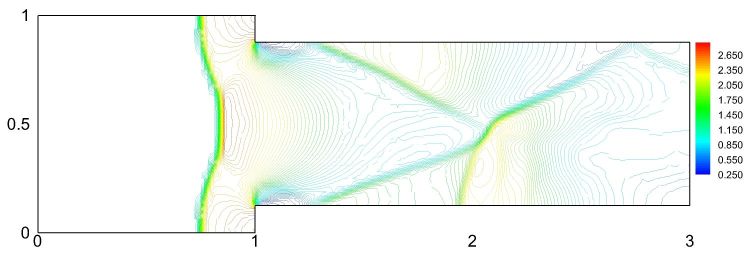}\hspace*{0.2cm}\includegraphics[width=8.5cm]{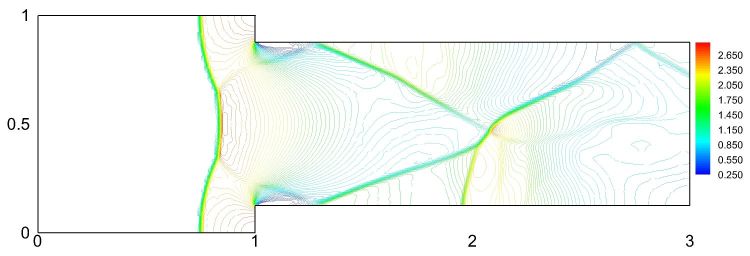}}
\vspace*{0.25cm}
\centerline{\includegraphics[width=8.5cm]{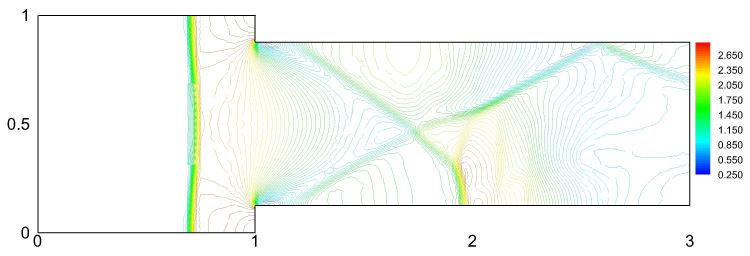}\hspace*{0.2cm}\includegraphics[width=8.5cm]{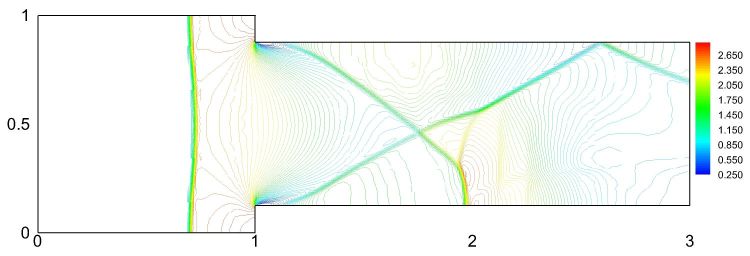}}
\vspace*{0.25cm}
\centerline{\includegraphics[width=8.5cm]{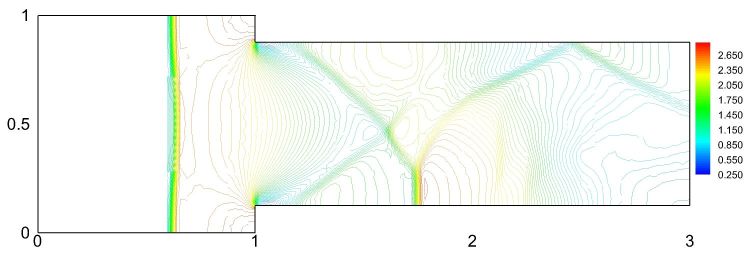}\hspace*{0.2cm}\includegraphics[width=8.5cm]{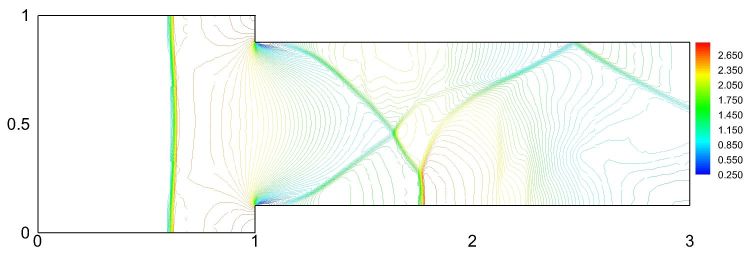}}
\caption{\sf Same as Figure \ref{fig:10}, but for the nonflat bottom topography.\label{fig:11}}
\end{figure}

We compute the solution until the final time $t=2$ in order to simulate a transient flow state. We plot the snapshots of $w$ at times
$t=0.5$, 1, $1.5$ and 2 in Figures \ref{fig:10} and \ref{fig:11} for the flat and nonflat bottom topographies, respectively. As one can
see, the proposed central-upwind quadtree scheme preserves positivity of the computed water depth and is able to capture hydraulic jumps.
Increasing $m$ from 8 to 9 clearly improves the accuracy and resolution of the hydraulic jumps. Finally, in Table \ref{Tbl:Tbl1}, we present
the minimum and maximum number of cells during the time evolution for different quadtree levels and topographies.
\begin{table}[ht!]
  \caption{\sf Example 6: Minimum and maximum number of cells for each of the four solutions.\label{Tbl:Tbl2}}
\centering
\begin{tabular}{*5c}
\toprule
Quadtree level&\multicolumn{2}{c}{$m=8$}&\multicolumn{2}{c}{$m=9$}\\
\midrule
{}&min&max&min&max\\
$B(x,y)\equiv0$&298&3154&436&10954\\
$B(x,y)\neq0$&298&5140&436&21340\\
\bottomrule
\end{tabular}
\end{table}

\section{Conclusion}\label{S:6}
An adaptive, well-balanced, positivity preserving central-upwind scheme over quadtree grids for the shallow water equations over irregular
bottom topography has been presented. Six numerical experiments have been performed in order to verify the accuracy and robustness of the
proposed scheme. The first numerical benchmark test has addressed the accuracy of the scheme. The second numerical example has focused
on the positivity and symmetry preserving as well as adaptability of the scheme. The third and fourth numerical examples have demonstrated
the well-balanced property, symmetry preserving and adaptability features of the proposed method. The fifth test has focused on the
capability of the scheme to model flows over a discontinuous bottom topography. The last numerical example has demonstrated the positivity
preserving and shock-capturing features of the method. The obtained results show that the proposed central-upwind quadtree scheme can
improve the performance and efficiency of calculations compared with regular Cartesian grids.


\acknowledgments{The work of A. Kurganov was supported in part by NSFC grant 11771201 and by the fund of the Guangdong Provincial Key
Laboratory of Computational Science and Material Design (No. 2019B030301001).}

\bibliographystyle{elsarticle-num-names-alphsort}
\bibliography{References}

\end{document}